\documentclass[12pt]{amsart}
\usepackage{div-alg-chars,amscd}
\usepackage[hyperindex]{hyperref}

\usepackage{amsrefs}
\newcommand{\xrhardcode}[2]{\expandafter\def\csname#1\endcsname{#2}}
\xrhardcode{J-rem:equiv-bdd}{2.9}
\xrhardcode{J-subsec:coefficient}{2.3}
\xrhardcode{J-defn:top-F-ss-unip}{2.15}
\xrhardcode{J-lem:bounded-and-top-F-unip}{2.18}
\xrhardcode{J-lem:top-F-unip-Jordan}{2.19}
\xrhardcode{J-lem:top-unip}{2.21}
\xrhardcode{J-rem:top-F-ss-unip-ascent}{2.22}
\xrhardcode{J-defn:top-F-Jordan}{2.23}
\xrhardcode{J-prop:unique-top-F-Jordan}{2.24}
\xrhardcode{J-lem:top-F-Jordan-central}{2.25}
\xrhardcode{J-subsec:existence}{2.8}
\xrhardcode{J-prop:x-depth}{2.33}
\xrhardcode{J-prop:tJd=cpct}{2.36}
\xrhardcode{J-cor:abs-F-ss-tame}{2.37}
\xrhardcode{J-prop:top-F-Jordan-defn}{2.42}
\xrhardcode{J-prop:tame-tunip}{2.43}

\renewcommand\tame{\textsup{tr}}
\author[Adler]{Jeffrey D. Adler}
\address{The University of Akron \\ Akron, OH \ \ 44325-4002}
\curraddr{American University \\ Washington, DC \ \ 20016-8050}
\email{jadler@american.edu}
\author[Spice]{Loren Spice}
\address{The University of Michigan \\ Ann Arbor, MI\ \ 48109-1043}
\email{lspice@umich.edu}
\thanks{The first-named author was partially supported by
the National Security Agency (H98230-05-1-0251 and H98230-07-1-002)
and by a University of Akron Faculty Research Grant (\#1604).
The second-named author was supported by
National Science Foundation Postdoctoral Fellowship award
DMS-0503107.}

\title[Good product expansions]%
{Good product expansions for tame elements of $p$-adic groups}
\date{30 May, 2008}

\begin{document}
\begin{abstract}
We show that, under fairly general conditions,
many elements of a $p$-adic group can be
well approximated by a product
whose factors have properties
that are helpful in performing explicit character computations.
\end{abstract}
\maketitle
\setcounter{tocdepth}{1}
\tableofcontents

\section{Introduction}

Suppose $F$ is a non-discrete, discretely valued, locally
compact field
of residual characteristic $p$
with valuation $\ord$,
$\varpi$ is a uniformizer of $F$,
and $\ff$ is the residue field of $F$.
(We will weaken these hypotheses on $F$ below.  See
especially \S\ref{sec:buildings}.)
Each element of $\ff$ is a coset in the ring of integers of $F$
of the unique maximal ideal in that ring.
The set \mc S consisting
of $0$ and the roots of unity in $F$ of order coprime to $p$
is a set of coset representatives for \ff.
Then every element $\gamma \in F\cross$ has a unique
expression of the form
$$
\gamma = \sum_{i = \ord(\gamma)}^\infty \delta_i\varpi^i
$$
with $\delta_i \in \mc S$ for $i \ge \ord(\gamma)$.
For some purposes,
it turns out to be more convenient to work with expressions
of the form
$$
\gamma = \varepsilon_0\varpi^{\ord(\gamma)}
\prod_{i=1}^\infty (1 + \varepsilon_i \varpi^i) ,
$$
where $\varepsilon_i \in \mc S$ for $i = 0, 1, 2, \dotsc$.
Of course, the terms in this product all commute
with each other.

More generally, suppose $D$ is a
central division algebra of degree $n$ over $F$.
From Proposition I.4.5 of \cite{weil:basic-2nd},
$D$ has a uniformizer $\varpi_D$ such that $\varpi_D^n=\varpi$,
and
there is a set \mc S of coset representatives for the
quotient of the ring of integers of $D$ by its unique maximal
ideal such that every element $\gamma\in D\cross$
has an expansion of the form
\begin{equation}
\tag{$*$}
\gamma
= \varepsilon_0\varpi_D^{\ord_D(\gamma)}
	\prod_{i = 1}^\infty (1 + \varepsilon_i\varpi_D^i),
\end{equation}
where $\ord_D$ is the unique valuation on $D$ extending
$n\dotm\ord$ on $F$,
and
$\varepsilon_i \in \mc S$ for $i = 0, 1, 2, \dotsc$.
When the degree $n$ is coprime to $p$, then
as observed in 
\cite{corwin-moy-sally:gll},
every element $\gamma\in D\cross$ is conjugate to one
for which the factors in ($*$)
all commute with each other.
(Such elements are called \emph{normal} in
\cite{adler-corwin-sally:division-formulas}.  This is the
motivation for our term `\emph{normal approximation}' in
Definition \ref{defn:r-approx} below.)
Then, for $r = 1, 2, \dotsc$, it is natural to consider the
truncation
$$
\gamma_{< r}
:= \varepsilon_0\varpi_D^{\ord_D(\gamma)}
	\prod_{i = 1}^{r - 1} (1 + \varepsilon_i\varpi_D^i).
$$
In general, $\gamma_{< r}$ is ``more singular'' than
$\gamma$ (i.e., has larger centralizer).
Although the specific element $\gamma_{< r}$ depends on the
choice of uniformizer $\varpi_D$, it turns out that the
division algebra $C_D(\gamma_{< r})$ is independent of this
choice.

If the degree $n$ of $D$ is prime and $\pi$ is a
representation of $D\cross$ which is trivial on the $r$th
filtration subgroup of $D\cross$, then $\gamma_{< r}$ is either
central or regular, a dichotomy that is reflected in the
differing formul\ae\ for the character $\Theta_\pi(\gamma)$ in each case.
Whether or not $n$ is prime,
the division algebras which arise as centralizers
of truncations of $\gamma$ play an important role in the character
formul\ae\ of \cite{adler-corwin-sally:division-formulas}
(specifically, see \S4 of loc.\ cit.).

As the character theory of reductive $p$-adic groups
is still in its early stages, there are
not many explicit supercuspidal character formul\ae\
available, but a number of qualitative results have
been proven.  Supercuspidal characters have been computed
for $\SL_2(F)$ in \cite{sally-shalika:characters};
for $\PGL_2(F)$ in \cite{silberger:pgl2};
and for $\GL_2(F)$ in \cite{shimizu:gl2}.  In these cases, the
character tables are complete, but it is necessary to place
some restriction on the residual characteristic of the underlying field $F$,
usually that it be odd.  Under similar conditions, complete
character tables are available for $\GL_3(F)$ in
\cite{debacker:thesis}.
Further explicit information is known for supercuspidal
characters of $\GL_\ell(F)$ and $\SL_\ell(F)$ with $\ell$ prime (see
\cite{corwin-moy-sally:gll},
\cite{debacker:thesis}, \cite{kutzko:supercuspidal-gl2-1},
and \cite{spice:thesis});
$\GL_n(F)$ (see
\cite{adler-corwin-sally:division-formulas}
and \cite{murnaghan:chars-gln})
and $\SL_n(F)$ (see \cite{murnaghan:chars-sln})
with $n$ not necessarily prime; and $\Sp_4(F)$ (see
\cite{boller:thesis}).
In all of these latter cases, we must place a
restriction on the residual characteristic $p$ of $F$, and
the actual character tables are not yet complete.  The most
detailed information that is available is for the so-called
``unramified supercuspidal representations'' of
$\GL_\ell(F)$.

The calculations in \cite{gerardin:gln} (see Proposition 7.4),
\cite{corwin-moy-sally:gll},
\cite{debacker:thesis}, and \cite{spice:thesis}
have suggested that,
since many supercuspidal representations are induced from
compact modulo center groups,
precise control over the terms appearing in (some analogue of)
the Frobenius formula
for induced characters
will be an important tool in the computation of characters
of such representations.
To this end, it is necessary to have a precise understanding
of the effect of conjugation on a regular semisimple
element.

Implicit in the above-mentioned papers, and explicit in the
calculations underlying
\cite{adler-corwin-sally:division-formulas}, is the idea of
a truncation of an element of a general or special linear group,
or of a division algebra (as described above).
In order to compute character formul\ae\
for other groups,
we need to generalize this notion,
hence to express an element
as an infinite product in a convenient way.
This is the main goal of the present paper.

In order for a product expansion to converge,
the terms should lie in smaller
and smaller subgroups in some filtration of $G$.
There is a canonical filtration when
$G=D\cross$, but not in general.
In \cite{moy-prasad:k-types} and
\cite{moy-prasad:jacquet}, Moy and Prasad
defined a collection of filtrations of $G$,
which will play the same role for us.
In fact, we need a slight generalization of the
Moy--Prasad filtrations; so, following Yu in
\cite{yu:supercuspidal}, in \S\ref{sec:concave} we define
filtration subgroups associated to certain functions on the
root system of our group.

In our forthcoming paper
\cite{adler-spice:explicit-chars}, we will apply the
structure theory results of this paper to compute the values
at elements $\gamma$ as above of the characters of many of the
supercuspidal representations constructed by Yu in
\cite{yu:supercuspidal}.

We now outline the content of this paper.
We no longer suppose that $F$ is locally compact.
We will assume for simplicity that it is complete and has
perfect residue field,
but this is not the most general possible condition under
which our results hold (see \S\ref{sec:buildings}).
Let \bG denote a connected, reductive $F$-group,
and $G$ its group of $F$-rational points.

In \S\ref{sec:assumptions}, we describe the hypotheses under
which we work in this paper.
Not all of our results require all hypotheses,
so we describe precisely where in the
document each assumption is used.

In \S\ref{sec:notation}, we recall some basic
definitions and notation, especially regarding Moy--Prasad
filtrations and the Bruhat--Tits building.  Since we are
working over very general fields, we also re-prove some
results about Moy--Prasad filtrations which are familiar in
the case where $F$ is locally compact.

In \S\ref{sec:tame-and-compatible}, we introduce the notion of
\emph{compatibly filtered} subgroups.
These are $F$-subgroups of \bG for which both the
Bruhat--Tits buildings and the Moy--Prasad filtrations sit
nicely inside the associated objects for \bG.
This notion is a mild generalization of that of a tame Levi
subgroup of \bG (see Lemmata \ref{lem:field-descent} and
\ref{lem:levi-descent}).
The principal results of this section, Proposition
\ref{prop:compatibly-filtered-tame-rank}
and Corollary \ref{cor:compatibly-filtered-tame-rank}, show
that the centralizer of a truncation of a tame-modulo-center element
of $G$ is compatibly filtered.
(This centralizer need not, in general, be a Levi subgroup of \bG.)

In \cite{yu:supercuspidal}*{\S\S 1--2}, Yu defined filtration subgroups
of $G$ associated to certain functions on the root
system of \bG.
However, Yu's subgroups are best behaved only when the
ambient group \bG is tame.
Since we wish to avoid such a strong assumption on \bG,
we present a slightly different definition in
\S\ref{sec:concave}, which coincides with Yu's definition
when \bG is tame.
We examine commutators of these subgroups (see, for example,
Lemma \ref{lem:shallow-comm}) and show that they support an
analogue of the Moy--Prasad map (see Lemma
\ref{lem:gen-iwahori-factorization}).
In \S\ref{sec:concave_tame}, we explore further the descent
properties of these subgroups under mild tameness hypotheses
on \bG.  Of particular note in this subsection is Proposition
\ref{prop:heres-a-gp}, which gives a very concrete picture
of our filtration subgroups in many cases.

Now we turn to analyzing the analogue in $G$ of the
product decomposition ($*$).
If $\gamma \in D$ and $\ord_D(\gamma) = 0$, then the first term,
$\varepsilon_0\varpi_D^{\ord_D(\gamma)} = \varepsilon_0$,
in the product decomposition is well defined
(i.e., independent of the choice of $\varpi_D$).  
The analogue in general of this decomposition of $\gamma \in D$
into the ``first term'' and the ``remaining terms'' of ($*$)
is the notion of a topological Jordan decomposition of a
compact element, which was introduced in \cite{hales:simple-defn}*{\S 3},
and Lemma 2 of \cite{kazhdan:lifting}*{\S 3}.
In a separate paper, one of us (L.~S.)\ defines the analogous notion of a
topological Jordan decomposition of an element of a
topological group modulo a subgroup
(see Definition \xref{J-defn:top-F-Jordan} of
\cite{spice:jordan}), and discusses when such
a decomposition exists
(see Proposition \xref{J-prop:tJd=cpct} of loc.\ cit.).

In \S\ref{sec:normal}, we begin to discuss the full product
decomposition for an element of $G$.
In fact, for our purposes, it is most useful
not to consider the entire infinite product at once, but
rather just finite approximations.  Thus we introduce the
idea of a \emph{(normal) $r$-approximation} (see Definition
\ref{defn:r-approx}), which will underlie the remainder of
the article.
We define some useful notation (see Definitions
\ref{defn:funny-centralizer} and
\ref{defn:bracket}) and observe a few basic but fundamental
facts (see Remarks \ref{rem:bracket-facts} and
\ref{rem:approx-facts}).

Before proving our main results about normal approximations,
we need to have a detailed understanding of the properties
of \emph{good elements} (see Definition \ref{defn:good}).
In \S\ref{sec:good-facts}, we study the effect of taking
commutators with such elements.
We omit a detailed description of the results of this
section, since most of them are very technical and intended
only for use in the next section.
The Lie algebra analogues of most of these
results are already known (see \cite{adler:thesis}*{\S 2}
and \cite{jkim-murnaghan:charexp}*{\S 2}).   

In \S\ref{sec:approx-facts}, we are in a position
to answer questions about the existence and uniqueness of
normal approximations.  Under mild hypotheses on \bG, all
tame and bounded modulo center elements of $G$ admit
normal approximations (see Lemma \ref{lem:simult-approx}).
Although normal approximations are by no means unique, they
are ``unique enough'', in a precise sense.
In particular, the
centralizer of a truncation of an element is uniquely determined
(see Proposition \ref{prop:unique-approx}).
Finally, in \S\ref{sec:normal-conj}, we address the question
of how to recognize when an element is close to a tame Levi subgroup.
This will be important for
our forthcoming character computations (see
\cite{adler-spice:explicit-chars}).
Lemma \ref{lem:rigidity} describes the elements which
conjugate something close to a tame Levi subgroup, into
something close to the tame Levi subgroup.
Proposition \ref{prop:aniso-Levi} characterizes the elements
which are close to a tame Levi subgroup
that is $F$-anisotropic modulo its center.

\textsc{Acknowledgements:}
This paper was motivated in large part by the notes of
the late
Lawrence Corwin on his computation of characters of
division algebras
and has benefited from our conversations with Paul Sally,
Gopal Prasad, Stephen DeBacker, Brian Conrad, and Jiu-Kang Yu.
It is a pleasure to thank all of these people.
	%%Theople.

\section{Some assumptions}
\label{sec:assumptions}

Suppose $F$ is a discretely valued field as in \S\ref{sec:buildings},
$F\tame$ is the maximal tamely ramified extension
of $F$ in a fixed algebraic closure,
and $\bG$ is a connected reductive $F$-group.

For easy reference, we collect here all the hypotheses
on $F$ and $\bG$
which
we will need throughout this document.
Many of our results remain valid if we only assume some subset
of these hypotheses.  See Remark \ref{rem:which-hyps} for a
more detailed discussion.
Although the hypotheses include some terms which have not
yet been defined, every term in each hypothesis will be
defined before that hypothesis is used.

\begin{hyps}
\hfill
\begin{enumerate}[\bf(A)]
\item
\label{hyp:reduced}
The relative root system $\lsub{F\tame}\Phi(\bG)$ is reduced.
\item
\label{hyp:conn-cent}
If $E/F$ is a discretely valued, separable extension and
$\gamma \in \bG(E)_{0+}$ is semisimple,
then $C_\bG(\gamma)\conn$ is a Levi subgroup of \bG.
\item
\label{hyp:good-weight-lattice}
There is a connected reductive $F$-group \mo J such that
	\begin{itemize}
	\item the absolute ranks of \bG and \mo J are the same,
	\item \bG is a compatibly filtered $F$-subgroup of \mo J
(in the sense of Definition \ref{defn:compatibly-filtered}),
and
	\item for some (hence every) maximal torus
$\bT'$ in the derived group $\mo J\semi$,
the order of $P(\mo J\semi, \bT')/\bX^*(\bT')$
is not divisible by $\chr\ff$.
	\end{itemize}
\item
\label{hyp:torus-H1-triv}
For any tower of discretely valued separable extensions $E/L/F$ with
$E/L$ Galois and unramified,
and any maximal $L$-torus \bT in \bG with the same
$L\tame$-split rank as \bG,
we have
$H^1(E/L, \bT(E)_{c:c{+}}) = \sset 0$ for $c \in \R_{> 0}$.
\end{enumerate}
\end{hyps}

\begin{rk}
\label{rem:which-hyps}
Hypothesis \eqref{hyp:reduced} applies in all of
\S\ref{sec:concave_tame}, and therefore is needed 
in the part of the document which relies on that section,
namely, everything from Lemma \ref{lem:iso-quotients}
onward.
Hypothesis \eqref{hyp:conn-cent} is used only in the proof of
Lemma \ref{lem:x-depth}.
Hypothesis \eqref{hyp:good-weight-lattice} is used only in the
proof of Lemma \ref{lem:GE2}.
Hypothesis \eqref{hyp:torus-H1-triv} is used in the proof of
Corollary \ref{cor:quotient-H1-triv}, and therefore again is
needed in the part of the document which relies on
\S\ref{sec:concave_tame}.
\end{rk}

\begin{rk}
\label{rem:when-hyps-hold}
There are several situations in which these hypotheses are known to hold.
In particular, we will see below that they all hold if
$\chr \ff = 0$.

Hypothesis \eqref{hyp:reduced} is satisfied if and only if,
for some finite tame Galois extension $L/F$,
$\lsub L\Phi(\bG)$ is reduced.
By \cite{springer:lag}*{\S 17}, Hypothesis \eqref{hyp:reduced}
is satisfied if the (absolute) root system of \bG contains no
irreducible factor of type $\ms A_{2n}$.
Hypothesis \eqref{hyp:reduced} obviously holds when
$\bG/Z(\bG)\conn$ (which we will later denote by \tbG)
contains an $F\tame$-split maximal $F$-torus;
in particular, when $\chr \ff = 0$.
Lemma \ref{lem:rel-root-red} shows that
Hypothesis \eqref{hyp:reduced}
holds when $\chr \ff \neq 2$.

Any finite-order root value of a semisimple
element of $\bG(E)_{0+}$ (where $E/F$ is a
discretely valued, separable
extension) has $p$-power order, where $p := \chr \ff$.
(The term ``root value'' is defined in Definition
\ref{defn:root-value}.)
Thus Proposition \ref{prop:levi} shows that Hypothesis
\eqref{hyp:conn-cent} is satisfied if $p$ is not a
bad prime for \bG (in the sense of Definition
\ref{defn:bad-prime}); in particular, if $p = 0$.

Each of the following conditions implies
Hypothesis \eqref{hyp:good-weight-lattice}:
\begin{itemize}
\item $\chr \ff = 0$;
\item \bG has simply connected derived group;
\item
$\chr \ff > 3$, and \bG contains no factor of type
$\ms A_n$ with $n$ divisible by $\chr \ff$.
\end{itemize}

By Proposition 5.5 of \cite{yu:models}, Hypothesis
\eqref{hyp:torus-H1-triv} is satisfied if all the relevant
tori satisfy condition (\textbf T) of \cite{yu:models}*{\S 4.7.1}.
This is true, in turn, if \bG is tame;
in particular, if $\chr \ff = 0$.
By Proposition 4.4.16 of \cite{bruhat-tits:reductive-groups-2},
it is also true if
\bG is adjoint or simply connected.
	%%To apply this result, one has to go up to an unramified
	%%extension, so that \bG is quasisplit; but that's
	%%OK, because (T) considers only the behaviour of
	%%tori after tame base change.
\end{rk}

\begin{rk}
Clearly,
Hypothesis \eqref{hyp:reduced}
is inherited by full-tame-rank closed
reductive subgroups of \bG;
Hypothesis
\eqref{hyp:good-weight-lattice} is inherited by
full-rank compatibly filtered $F$-subgroups of \bG;
and Hypothesis \eqref{hyp:torus-H1-triv} is inherited by
full-rank, full-tame-rank $F$-subgroups of \bG.

Hypothesis \eqref{hyp:conn-cent} is inherited by
connected full-rank compatibly filtered $F$-subgroups of \bG.
To see this,
suppose Hypothesis \eqref{hyp:conn-cent} holds for \bG,
let \bH be such a subgroup,
and suppose $\gamma \in \bH(E)_{0+}$ is semisimple
for some separable extension $E/F$.
Then $\gamma \in \bG(E)_{0+}$\,, so
$\bM := C_\bG(\gamma)\conn$ is a Levi subgroup of \bG; that
is, $\bM = C_\bG(\bS)$, where $\bS = Z(\bM)\conn$.
By hypothesis, there is a maximal $F$-torus $\bT \subseteq \bH$
such that $\gamma \in T$.  Then \bT is a maximal torus in \bM, so
$\bS \subseteq \bT \subseteq \bH$.  Thus
$\bM \cap \bH = C_\bG(\bS) \cap \bH = C_\bH(\bS)$
is a Levi subgroup of \bH.
Clearly
$C_\bH(\gamma)\conn \subseteq \bM \cap \bH \subseteq C_\bH(\gamma)$.
Since $\bM \cap \bH$ is connected, we have
$\bM \cap \bH = C_\bH(\gamma)\conn$.
That is, Hypothesis \eqref{hyp:conn-cent} also holds for \bH.

Further, all hypotheses
are preserved under base change
(to a discretely valued separable extension of $F$).
As a consequence, whenever we have proved a result for the
group $G = \bG(F)$, we will feel free to use it also for a
group $\bH(E)$, where \bH is a full-rank, full-tame-rank
reductive (necessarily closed and compatibly filtered)
$F$-subgroup of \bG and $E/F$ is a discretely valued
separable extension.
	%%We do this (cite results about $G$ as if they were
	%%about $\bH(E)$).  Since we explicitly mention
	%%the same behaviour as it applies to notation,
	%%we might as well do it here, too.
\end{rk}

\section{Preliminaries}
\label{sec:notation}

\subsection{Generalities on linear reductive groups}
\label{sec:notation-generalities}
For an abstract group $G$, let $Z(G)$ denote the center of $G$.
For a field $F$ and a linear algebraic $F$-group $\bG$,
let $Z(\bG)$ denote the center of $\bG$.
(By this, we do not mean the scheme-theoretic center of $\bG$,
but rather the underlying reduced scheme.
Note that, for example, if $\chr F = p$ and $\bG = \SL_p$,
then $Z(\bG)$ is the trivial variety, not the scheme $\mu_p$
whose ring of regular functions is $F[X]/(X^p-1)$.
This makes a difference in Remark \ref{rem:G-image}.)
By Theorem 18.2(ii) of \cite{borel:linear}, if \bG is connected and
reductive, then $Z(\bG)$ is defined over $F$
and
$Z(\bG)(F) = Z(\bG(F))$.
\indexmem{\tbG}
Denote by $\bG\conn$ the identity component of \bG,
and by \tbG the quotient $\bG/Z(\bG)\conn$.
If \bG is connected and reductive, then, by
Proposition 22.4 of \cite{borel:linear}, the
natural $F$-homomorphism $\bG \to \tbG$ is central (in the
sense of \cite{borel:linear}*{\S 22.3}).
The image of $g \in \bG(\ol F)$ under this map will be
denoted by \ol g.

We will use bold letters to denote algebraic groups,
and the corresponding normal letters to denote their groups
of $F$-rational points.
Thus, for example, $G=\bG(F)$.
Note, however, that $\tG = \tbG(F)$ need not equal $G/Z(G)\conn$.
For convenience, we will often define notation such
as, say, $G_{x, r}$ (see \S\ref{sec:filtrations-and-depth}), and
omit the analogous definition of the notation
$\bG(E)_{x, r}$\,, for $E/F$ a discretely valued algebraic
extension.

For $g,h\in \bG(\ol F)$, let $\Int(g)$ denote the inner automorphism
of $\bG$ given by $x\mapsto gxg\inv$;
$[g,h]$ denote $ghg\inv h\inv$;
and $C_\bG(g)$ denote the centralizer of $g$ in $\bG$.
%%%% Why the centralizer is defined over F:
% By Proposition 13.19 of \cite{borel:linear},
% if $g$ is semisimple, then $C_\bG(g)$ is reductive.
% By Proposition 9.1(1) of \cite{borel:linear},
% if also $g \in G$, then $C_\bG(g)$ is defined over $F$.
We will write $\lsup g h = \Int(g)h$
and $h^g = \Int(g\inv)h$.
For subsets $S, S' \subseteq \bG(\ol F)$,
we define $\lsup S h$, $\lsup h S'$, $\lsup S S'$, etc., in
the obvious way.
(Use similar notations for an abstract group $G$.)
Put $\Ad(g) = d(\Int(g))$ and $\Ad^*(g) = \Ad(g)^*$.
If $X \in \Lie(\bG)(\ol F)$ and $X^* \in \Lie(\bG)^*(\ol F)$,
we will write
\begin{align*}
\lsup g X & = \Ad(g)X, & X^g & = \Ad(g\inv)X, \\
\lsup g X^* & = \Ad^*(g)X^*, & (X^*)^g & = \Ad^*(g\inv)X.
\end{align*}
For a subset $S' \subseteq \Lie(\bG)(\ol F)$ or
$S' \subseteq \Lie(\bG)^*(\ol F)$, define
$\lsup g S'$, etc., as above.

Let $\bX^*(\bG) = \Hom\alg(\bG,\bGL_1)$
and $\bX_*(\bG) = \Hom\alg(\bGL_1, \bG)$.
For any extension $E/F$, 
denote by $\bX^*_E(\bG)$ the set of $\chi \in \bX^*(\bG)$
defined over $E$, and similarly for $\bX_*^E(\bG)$.

If $g\in \bG(\ol F)$, then let $g\semi$ and $g\unip$
denote the semisimple and unipotent parts, respectively,
of the Jordan decomposition of $g$.
Note that, if $g$ is $F$-rational, then $g\semi$ and $g\unip$
are defined over some finite, totally inseparable extension of $F$.

If a torus \bS acts on \bG, we denote by $\wtilde\Phi(\bG, \bS)$
and by $\Phi(\bG, \bS)$ the collections of weights and of
non-zero weights for the corresponding action on $\Lie(\bG)$.

From now on, \bG is a linear reductive $F$-group.
If \bS is a maximal $F$-split torus in \bG, then
$\Phi(\bG, \bS)$ is a root system in
$\bX^*(\bS') \otimes_\Z \Q$
by Theorem 21.6 of \cite{borel:linear}.
(Here, $\bS'$ is the identity component of the intersection
of \bS with the derived group of \bG.)
We will sometimes write $\lsub F\Phi(\bG)$ instead of
$\Phi(\bG, \bS)$ if we only care about the isomorphism
type of the root system.
We will call the elements of $\lsub F\Phi(\bG)$ \emph{$F$-roots}.
If $\bS = \bS'$ (for example, if \bG is semisimple), then
we define the \emph{weight lattice} $P(\bG, \bS)$
to be the space of those
$\chi \in \bX^*(\bS) \otimes_\Z \Q$
such that
$\langle \chi,\alpha\spcheck\rangle \in\Z$
for all $\alpha\spcheck \in \Phi\spcheck(\bG,\bS)$,
where $\Phi\spcheck(\bG, \bS)$ is the coroot system dual to
$\Phi(\bG, \bS)$ and
$\langle\cdot, \cdot\rangle$ is the natural pairing
of $\bX^*(\bS) \otimes_\Z \Q$ with $\bX_*(\bS) \otimes_\Z \Q$.

By Grothendieck's theorem (see
Theorem 18.2(i) of \cite{borel:linear}),
a torus in \bG is a maximal $F$-torus if and only if it is a
maximal torus which is defined over $F$.  We will often use
this fact without further remark.

Put $\tR := \R \sqcup  \set{ r{+} }{r \in \R} \sqcup \sset\infty$,
and extend the ordering on $\R$ to one on $\tR$ as follows:
for all $r,s\in \R$,
\begin{align*}
r < s{+} & \quad\text{if and only if}\quad r\leq s; \\
r{+} < s{+}&\quad\text{if and only if}\quad r< s; \\
r{+} < s{\phantom{+}}& \quad\text{if and only if}\quad r< s; \\
\end{align*}
and $r, r{+} < \infty$.
For $r\in\R$ and $\lambda \in \R_{> 0}$, define $(r{+}){+} := r{+}$
and $\lambda(r{+}) := (\lambda r){+}$.
Define also $\lambda\infty := \infty$ and $\infty{+} := \infty$.
Extend the additive structure
on \R to an additive structure on \tR
in the natural way.
Let $\tR_{> 0} := \set{r\in\smash{\tR}}{r > 0}$
and $\tR_{\ge 0} := \tR_{> 0} \cup \sset 0$, and define
$\R_{> 0}$ and $\R_{\ge 0}$ similarly.

\subsection{Buildings and affine root groups}
\label{sec:buildings}

Assume from now on that $F$
has a non-trivial discrete valuation $\ord$,
and that there is a subfield of $F$,
over which $F$ is algebraic,
which is complete and has perfect residue field.
For any algebraic extension $E/F$, denote again by $\ord$
the unique extension of $\ord$ to a valuation on $E$.
If this extended valuation $\ord$ remains discrete, we will
say that $E/F$ is a discretely valued algebraic extension.
Note that this happens precisely when $E/F$ has finite ramification
degree.
Fix an algebraic closure \ol F of $F$.
\indexmem{F\unram}\indexmem{F\tame}\indexmem{F\sep}%
Let $F\unram/F$, $F\tame/F$ and $F\sep/F$ be the maximal
unramified, tame, and separable subextensions of $\ol F /F$, respectively.
% In \cite{serre:local-fields}, Serre says that
% an extension of discretely valued fields is \emph{unramified}
% if its ramification degree is $1$, \emph{and} the corresponding
% extension of residue fields is separable.
% However, we don't need to worry about this, since
% we're always assuming that our residue fields are perfect.
For $r \in \tR$, put $F_r = \set{t \in F}{\ord(t) \ge r}$.
\indexmem{\ff_F}
Put $\ff_F = F_0/F_{0+}$, the residue field of $F$.
When $F$ is understood, we will often just write \ff.
Put also $F\cross_0 = \set{t \in F}{\ord(t) = 0}$
and (if $r > 0$) $F\cross_r = 1 + F_r$\,.

By Proposition 16.4.9 of \cite{springer:lag}, 
some twist of
\bG by an element of the cohomology set
$H^1(F\sep/F\unram, \bG\subad(F\sep))$
is $F\unram$-quasisplit,
where $\bG\subad$ is the adjoint group of \bG.
By Theorem 12 of \cite{lang:quasi-alg-closure}, $F\unram$
is $C_1$ (as in the second definition of
\cite{lang:quasi-alg-closure}*{p.~374}).
By Corollary II.3.2 of \cite{serre:galois},
$\dim(F\unram) \le 1$, so,
by \cite{borel-springer:reductive-groups-2}*{\S 8.6},
$H^1(F\sep/F\unram, \bG\subad(F\sep)) = \{0\}$.
Thus \bG is $F\unram$-quasisplit.  We will use this fact
frequently without mention.

Let $\BB(\bG,F)$ denote the (enlarged) Bruhat--Tits
building, and $\BB\red(\bG,F)$ the
reduced building, of \bG over $F$.
Note that $\BB\red(\bG, F)$ and $\BB(\tbG, F)$ are
canonically isomorphic, and
$$
\BB(\bG,F) = \BB\red(\bG,F) \times V_F(Z(\bG))
$$
\indexmem{V_F(Z(\bG))}
(where $V_F(Z(\bG))$ is an affine space under
$\bX_*^F(\bG) \otimes_\Z \R = \bX_*^F(Z(\bG)) \otimes_\Z \R$).
Denote by $\pr$ or $\pr_F$ the natural projection
$\BB(\bG, F) \to \BB\red(\bG, F)$,
and by \ox the image of a typical element $x \in \BB(\bG, F)$
under $\pr_F$.
If $E/F$ is a discretely valued algebraic extension,
then there are canonical embeddings
$\BB(\bG,F) \hookrightarrow \BB(\bG,E)$
and $\rBB(\bG, F) \hookrightarrow \rBB(\bG, E)$,
and we have that $\pr_E\bigr|_{\BB(\bG, F)} = \pr_F$.
If $E/F$ is Galois, then $\Gal(E/F)$ acts on $\BB(\bG,E)$,
and we have that
$\BB(\bG,F) \subseteq \BB(\bG,E)^{\Gal(E/F)}$,
with equality when $E/F$ is tamely ramified
(see Proposition 5.1.1 of \cite{rousseau:thesis}).

Suppose that \bS is a maximal $F$-split torus in \bG
(but not necessarily a maximal torus).
For $a \in \Phi(\bG, \bS)$,
Proposition 21.9(i) of
\cite{borel:linear} shows that there is a unique connected
$F$-subgroup $\bU_a \subseteq \bG$, normalized by
$C_\bG(\bS)$, such that $\Lie(\bU_a)$ is the direct sum
of the root subspaces of $\Lie(\bG)$ corresponding to
positive integer multiples of $a$.
If $a = 0$, then put $\bU_a = C_\bG(\bS)$.
If $a \not\in \wtilde\Phi(\bG, \bS)$, then put
$\bU_a = \sset 1$.

\indexmem{\AA(\bS, F)}
We will denote by $\AA(\bS, F)$ the ``empty'' apartment
associated to \bS over $F$.  (As in
\cite{landvogt:compactification}, we use the word ``empty''
to indicate that there is no polysimplicial structure.
However, for convenience, we \emph{will} regard it as
equipped with its natural metric as an affine Euclidean
space.%
	%%See \cite{bruhat-tits:reductive-groups-1}*{\S 2.5}
	%%for the definition of the metric on the building.
)
If we regard $\AA(\bS, F)$ as a subset of $\BB(\bG, F)$, and
wish to emphasize its polysimplicial structure, we will
write $\lsub\bG\AA(\bS, F)$ instead.
Recall that the space $\bX_*(\bS)\otimes_\Z \R$
acts transitively on $\AA(\bS,F)$, giving the latter
the structure of an affine space.
The gradient $\dot\varphi$ of
an affine function $\varphi$ on $\AA(\bS,F)$
is thus a linear function on $\bX_*(\bS) \otimes_\Z \R$.
In particular, via the natural pairing
$\bX^*(\bS) \times \bX_*(\bS) \to \Z$,
we may identify $\dot\varphi$ with an element of
$\bX^*(\bS) \otimes_\Z \R$.

For any affine function $\varphi$ on $\AA(\bS, F)$, define
$\varphi{+}$ to be the function on $\AA(\bS, F)$ given by
$(\varphi{+})(x) = (\varphi(x)){+}$ for $x \in \AA(\bS, F)$,
and put $(\varphi{+})\spdot := \dot\varphi$.
We will also refer to $\varphi{+}$ as an affine function.
\indexmem{U_\varphi}
For every root $a \in \Phi(\bG,\bS)$,
there are filtrations
$(U_\varphi)_\varphi$
and
$(\mf u_\varphi)_\varphi$
of the corresponding root group
$U_a$
and root space
$\Lie(U_a)$,
respectively,
both indexed by affine functions $\varphi$ on
$\AA(\bS,F)$ such that $\dot\varphi=a$.
(The indexing of the filtration depends on our choice of
valuation $\ord$.)
\indexmem{\lsub F U_\varphi}
If necessary, we will write
$\lsub F U_\varphi$
instead of just $U_\varphi$
to indicate the dependence on the field $F$.
We will not define these objects here (but see
\S\ref{sec:unip-exp} and
the proof of Proposition \ref{prop:compatibly-filtered-tame-rank}).

\indexmem{\PPsiGTF{}{\bG}{\bS}{E}}
\indexmem{\PPsiGTF{\protect\wtilde}{\bG}{\bS}{E}}
Denote by \PPsiGTF{}\bG\bS F the set consisting of the
affine functions
$\varphi$
on $\AA(\bS, F)$
such that $\dot\varphi \in \Phi(\bG, \bS)$,
together with the constant function
$\varphi = \infty$ on $\AA(\bS, F)$.
Denote by \PPsiGTF\wtilde\bG\bS F
the union of \PPsiGTF{}\bG\bS F
with the collection of $\tR_{\ge 0}$-valued
constant functions on $\AA(\bS, F)$.
The elements of $\PPsiGTF\wtilde\bG\bS F$ will be called
\emph{affine roots} (or \emph{affine $F$-roots}).
\indexme{affine root (unusual usage)}
Note that this is contrary
to the usual usage, which calls an element
$\varphi \in \PPsiGTF\wtilde\bG\bS F$ an affine root only if
$U_\varphi \ne U_{\varphi+}$
(or $U_\varphi \ne U_{\varphi{+}}\dotm U_{2\varphi}$\,,
if $\dot\varphi$ is multipliable).

Suppose that $E/F$ is a discretely valued algebraic extension such
that \bG contains a maximal $E$-split torus $\bS^\sharp$
defined over $F$.
For any affine $E$-root $\psi$,
$\lsup\sigma(\lsub E U_\psi) = \lsub E U_{\psi \circ \sigma}$ for
$\sigma \in \Aut(E/F)$;
and, if $\bS^\sharp = \bS$, then
$\lsub F U_\psi = \lsub E U_\psi \cap U_\dpsi$\,.

A proof of the following lemma appears in the proof of
Proposition 1.4.1 of \cite{adler:thesis}.
\begin{lm}
\label{lem:filter-descent}
Suppose that \bG is $F$-quasisplit and
\bS is a maximal $F$-split torus in \bG.
For any affine function $\varphi$ on $\AA(\bS, F)$,
$\lsub F U_\varphi
= \bigl(\prod \lsub E U_\psi\bigr) \cap \bU_{\dot\varphi}(F)$,
where $E$ is the splitting field of $C_\bG(\bS)$ and the
product runs over all affine $E$-roots $\psi$ whose
restriction to $\AA(\bS, F)$ is $\varphi$.
\end{lm}

Because of the following lemma,
the fact that our field $F$ need not be complete
does not cause any serious difficulties.
The statement is lengthy only so that the result can cover
all necessary applications; the proof itself is quite easy.

\begin{lm}
\label{lem:complete-subfield}
Suppose that
\begin{itemize}
\item $M'$ is a subfield of $F$ such that
	\begin{itemize}
	\item $F/M'$ is algebraic,
	\item $M'$ is complete with respect to the restriction of
$\ord$,
	and
	\item $\ff_{M'}$ is perfect;
	\end{itemize}
and
\item $\bG^1, \dotsc, \bG^n$ are finitely many
$F$-groups (not necessarily connected or reductive).
\end{itemize}
Suppose that we are given also, for each $1 \le i \le n$,
\begin{itemize}
\item a finite subset $\mc S^i$ of $G^i$;
\item an $F$-split torus $\bT^i$ in $\bG^i$;
and
\item a finite subset $\mc F^i$ of $\BB(\bG^i, F)$
(whenever the building makes sense).
\end{itemize}
Then there is a finite subextension $F'/M'$ of $F/M'$ such
that $F'$ is complete with respect to $\ord$, $F/F'$ is
unramified, and, for each $1 \le i \le n$,
\begin{itemize}
\item
$\bG^i$ and $\bT^i$ are defined over $F'$,
\item
$\mc S^i \subseteq \bG^i(F')$;
\item
$\mc F^i \subseteq \BB(\bG^i, F')$ (whenever the
building makes sense);
and
\item
$\bT^i$ splits over $F'$.
\end{itemize}
\end{lm}

\begin{proof}
Let $\varpi$ be a uniformizer for $F$ and put
$M = M'[\varpi]$.  Then the value groups of
$M$ and $F$ are the same.  Since the residue field of
$M'$, hence also of $M$, is perfect,
$\ff_F/\ff_M$ is separable.
Thus $F/M$ is unramified (in particular, separable).
% Here is why $F/M$ is separable.
%
% If $F/M$ is not separable,
% then there is a tower of fields $F/A/M$ such that $F/A$ is
% purely inseparable.  Without loss of generality, $M = A$.
% Choose $t \in F$.
% By the definition of pure inseparability, there is some
% natural number $N$ such that $t^{p^N} \in M$.
% (Here, $p = \chr M$, which is non-zero since $M$ is not
% perfect.)
% Since
% $\ord(t^{p^N}) = p^N\ord(t) \in p^N\ord(F) = p^N\ord(M)$,
% we may, and hence do, assume, after multiplying $t$ by an
% element of $M$, that $\ord(t) = 0$.
% Further, since $\ff_M$ is perfect, there is $s \in \ord(M)$
% such that $s^{p^N} \equiv t^{p^N} \pmod{M_{0+}}$.
% Thus we may, and hence do, assume that $t \in M\cross_{0+}$\,.
% Then $t - 1$ is a root of the Eisenstein polynomial
% $(X + 1)^{p^N} - t^{p^N} \in M[X]$,
% so $\ord(M) = \ord(F) \supseteq \ord(M[t - 1]) = p^{-N}\ord(M)$.
% Thus, $N = 0$, so $t \in M$.
% Since $t \in F$ was arbitrary, we have $F = M$, a
% contradiction.

We may, and hence do, assume that each
$\bT^i$ is some $\bG^j$.
Since each $\bG^i$ is an affine $F$-variety, there
exist integers $m$ and $r$ and polynomials
$f_{11}(\vec x), \dotsc, f_{nr}(\vec x) \in F[x_1, \dotsc, x_m]$
such that the ring of regular functions on $\bG^i$ is of
the form
$$
\ol F[x_1, \dotsc, x_m]
	/
\langle f_{ij}(\vec x) : j = 1, \dotsc, r\rangle
$$
for $i = 1, \dotsc, n$.
Denote by $\mc X^i$ an integral basis of
$\bX_*(\bT^i) = \bX_*^F(\bT^i)$.
Now let $N$ be the fixed field in $M\sep = F\sep$
of the stabilizer in $\Gal(M\sep/M)$ of
$\set{f_{ij}(\vec x)}{i = 1, \dotsc, n; j = 1, \dotsc r}$,
and $F'$ the fixed field in $N\sep = F\sep$ of the
common stabilizer in $\Gal(N\sep/N)$ of
the various $\mc X^i$, $\mc F^i$, and $\mc S^i$.
Since all of the sets in question are finite, the relevant
stabilizers are open, so $N/M$ and $F'/N$,
hence $F'/M$, are finite.  Since $M/M'$ is finite, so is
$F'/M'$.
Further, since $\Gal(F\sep/F)$ fixes the various
$\mc X^i$, $\mc F^i$, $\mc S^i$, and
$f_{ij}(\vec x)$, we have that $F' \subseteq F$.
\end{proof}

\subsection{Les \'epinglages}
\label{sec:unip-exp}

In this subsection only,
let \bT be a maximal $F$-torus in \bG, and $E/F$ a
discretely valued
Galois extension over which \bT splits.
For each $\alpha \in \Phi(\bG, \bT)$, denote by $F_\alpha$
the fixed field in $E$ of $\stab_{\Gal(E/F)} \alpha$.
By Theorem 18.7 of \cite{borel:linear}, there is an
$F_\alpha\sep$-isomorphism (in fact, an $E$-isomorphism)
$\Add \to \bU_\alpha$ such that, for any $t \in \bT(\ol F)$,
the pullback of the conjugation action of $t$ on $\bU_\alpha$ is
scalar multiplication by $\alpha(t)$ on $\Add$.
Since $\Add$, \bT, and $\bU_\alpha$ are defined over
$F_\alpha$, the set of such isomorphisms is a $\GL_1$-torsor
over $F_\alpha$.  Since
$H^1(F_\alpha\sep/F_\alpha, \GL_1(F_\alpha\sep)) = \sset 0$,
the torsor has a $\Gal(F_\alpha\sep/F_\alpha)$-fixed point.
Choose such a point, i.e., isomorphism, and call it $\mexp_\alpha$.
Since
$\sigma \circ \mexp_\alpha \circ \sigma\inv = \mexp_\alpha$
for $\alpha \in \Phi(\bG, \bT)$ and
$\sigma \in \stab_{\Gal(E/F)} \alpha$,
we may, and hence do,
make our choices in such a way that
$\sigma \circ \mexp_\alpha \circ \sigma\inv
= \mexp_{\sigma\alpha}$
for $\alpha \in \Phi(\bG, \bT)$ and
$\sigma \in \Gal(E/F)$.

Now suppose that $L/F$ is a Galois
subextension of $E/F$ such that
\begin{itemize}
\item \bG is $L$-quasisplit,
\item \bT contains a maximal $L$-split torus $\bS^\sharp$
in \bG which is defined over $F$,
and
\item $\Phi(\bG, \bS^\sharp)$ is reduced.
\end{itemize}
Define fields $L_\alpha$ for $\alpha \in \Phi(\bG, \bT)$ by
analogy with the fields $F_\alpha$ above,
and let $\mc O_L$ and $\mc O_{L_\alpha}$ be the rings of integers
of $L$ and $L_\alpha$, respectively.
As in \cite{bruhat-tits:reductive-groups-2}*{\S 4.1.5},
the $E$-\'epinglage $(\mexp_\alpha)_{\alpha \in \Phi(\bG, \bT)}$
of $\bG \otimes_F E$ gives rise to
an $L$-\'epinglage $(\mexp_a)_{a \in \Phi(\bG, \bS^\sharp)}$
of $\bG \otimes_F L$,
from which in turn we deduce, as in
\cite{bruhat-tits:reductive-groups-2}*{\S 4.3},
for each affine $L$-root $\varphi$ an
$\mc O_L$-scheme $\mf U_\varphi$ such that
\begin{itemize}
\item the generic fiber of $\mf U_\varphi$ is
$\bU_{\dot\varphi}$;
\item if $\alpha \in \Phi(\bG, \bT)$ restricts to
$\dot\varphi$, then there is an $\mc O_L$-isomorphism
$R_{\mc O_{L_\alpha}/\mc O_L}\Add \to \mf U_\varphi$
which induces the map $\mexp_{\dot\varphi}$ on generic fibers;
and
\item $\mf U_\varphi(\mc O_L) = U_\varphi$.
\end{itemize}
Moreover, if $\varphi_1 \ge \varphi_2$ and
$\dot\varphi_1 = \dot\varphi_2$, then there is an
inclusion
$\mf U_{\varphi_1} \hookrightarrow \mf U_{\varphi_2}$
which induces an isomorphism on generic fibers.

Fix an element $a \in \Phi(\bG, \bS^\sharp)$.
For any affine $L$-root $\varphi$ with $\dot\varphi = a$,
$\Lie(\bU_{\dot\varphi})$ is the generic fiber of
$\Lie(\mf U_\varphi)$.
By \cite{yu:models}*{\S 8.7},
$\mf u_\varphi = \Lie(\mf U_\varphi)(\mc O_L)$.
Thus
the natural $\mc O_L$-isomorphism
$\Add \cong \Lie(\Add)$
furnishes an isomorphism
$\mexp_\varphi : \mf u_\varphi \to U_\varphi$.
These maps are compatible, in the sense that, if
$\varphi_1 \ge \varphi_2$ and $\dot\varphi_j = a$ for $j = 1, 2$,
then
$\mexp_{\varphi_2}\bigr|_{\mf u_{\varphi_1}} = \mexp_{\varphi_1}$.
Thus they may be pieced together to form an isomorphism
$\Lie(\bU_a)(L) \to \bU_a(L)$, which (by abuse of notation)
we will denote by $\mexp_a$.
By following the details of the construction, one verifies
easily that
$\sigma \circ \mexp_a \circ \sigma\inv = \mexp_{\sigma a}$
for $\sigma \in \Gal(L/F)$.

\subsection{Filtrations and depth}
\label{sec:filtrations-and-depth}
The group $G$ acts on $\BB(\bG,F)$ by isometries.  We call a
subgroup $K$ of $G$ \emph{bounded} if its image in the isometry
group of $\BB(\bG, F)$ is bounded in the sense of
Exemple 3.1.2(b) of
\cite{bruhat-tits:reductive-groups-1}, i.e., if and only if
the orbit under $K$ of any bounded subset of $\BB(\bG, F)$
remains bounded.
Clearly, $\stab_G(x)$ is bounded for any $x \in \BB(\bG, F)$.
Conversely, suppose that $K$ is a bounded subgroup of $G$,
and fix $\ol y \in \rBB(\bG, F)$ (arbitrarily).
By Proposition 3.2.4 of loc.\ cit.,
$K$ fixes a point \ox in the closure of the convex hull of
$K\dota\ol y$, hence acts by translations on any lift
$x \in \BB(\bG, F)$ of \ox.
Since the orbit of $x$ under $K$ is bounded,
we have that $K \subseteq \stab_G(x)$.
If \mo N is a closed normal subgroup of \bG, then we call a
subgroup $K$ of $G$ \emph{bounded modulo \mo N{}} if its image
in $(\bG/\mo N)(F)$ is bounded.
Then a subgroup of $G$ is bounded modulo center if it fixes
a point in $\rBB(\bG, F)$.
Call an element of $G$ \emph{bounded}
(respectively, \emph{bounded modulo \mo N{}})
if it belongs to a subgroup of $G$ that
is bounded
(respectively, bounded modulo \mo N).
When $F$ is locally compact (equivalently, \ff is finite),
a subgroup is bounded if and only if it is pre-compact.
When $F$ is an algebraic extension of a locally
compact field, an element is bounded if and only
if it belongs to a compact subgroup.
For every point $x\in \BB(\bG,F)$, the stabilizer
$\stab_G(x)$ is open and (as we have observed) bounded,
and contains a normal and finite-index subgroup
$G_x$\,, the \emph{parahoric subgroup} associated to $x$,
which depends only on the image \ox of $x$ in $\rBB(\bG, F)$.
(See D\'efinition 5.2.6 of
\cite{bruhat-tits:reductive-groups-2}.)

\begin{lm}
\label{lem:move-within-apt}
Suppose that
\begin{itemize}
\item \bG is $F$-quasisplit,
\item $x, y \in \AA(\bS, F)$,
and
\item $g \in G$ satisfies
$g\dota x = y$.
\end{itemize}
Then $g \in N_G(S)\dotm G_x$\,.
\end{lm}

\begin{proof}
We have that $x \in \AA(\bS, F) \cap \AA(\bS^g, F)$.
By Proposition 4.6.28(iii) of
\cite{bruhat-tits:reductive-groups-2}, there is
$h \in G_x$ such that $\bS^{g h} = \bS$.  That is,
$g h \in N_G(S)$, so $g \in N_G(S)\dotm G_x$\,.
\end{proof}

The set $T$ of $F$-rational points of an $F$-torus \bT
comes equipped with a natural filtration.
First, note that $T$ has a unique maximal bounded subgroup
$T\subb$\,, and a unique parahoric subgroup $T_0$\,.
These two groups are equal if \bT is $F$-split,
or, more generally, if, over some unramified extension
of $F$, $\bT$ is a product of induced tori;
but, in general, we can only say that
$T_0$ is a finite-index subgroup of $T\subb$\,.
For $r\in\tR_{\geq0}$, we put
$$
T_r = \sett{t\in T_0}{$\ord(\chi(t) - 1) \geq r$
	for all $\chi\in \bX^*(\bT)$}.
$$
The filtration
$(\Lie(T)_r)_{r \in \tR}$
on $\Lie(T)$ is defined similarly
(except that we do not have to worry about passing to
finite-index subgroups).

\indexmem{\Lie(G)_{x, r}}
\indexmem{G_{x,r}}
For $(x, r) \in \BB(\bG, F) \times \tR$
with $r < \infty$,
Moy and Prasad (see \cite{moy-prasad:k-types} and
\cite{moy-prasad:jacquet}) define
a lattice
$\Lie(G)_{x, r}$
in $\Lie(G)$
and, if $r\geq 0$,
a bounded open subgroup
$G_{x,r}$ of $G$
as follows.
By Corollaire 5.1.12 of
\cite{bruhat-tits:reductive-groups-2}, there exists a
maximal $F\unram$-split torus which is defined
over $F$.  Fix such a torus $\bS^\sharp$.
Since \bG is $F\unram$-quasisplit,
$\bT := C_\bG(\bS^\sharp)$ is a maximal torus in \bG.
Put
\begin{gather*}
\bG(F\unram)_{x,r} =
	\langle
	\bT(F\unram)_r,
	U_\psi : \psi\in \PPsiGTF{}\bG{\bS^\sharp}{F\unram},
			\psi(x) \ge r
	\rangle \quad
\text{if $r \ge 0$}
\\
\intertext{and}
\Lie(\bG)(F\unram)_{x, r}  =
	\Lie(\bT)(F\unram)_r
	\oplus
	\sum_{ \substack{\psi\in \PPsiGTF{}\bG{\bS^\sharp}{F\unram} \\
			\psi(x) \ge r} }
	\mf u_\psi\,.
\end{gather*}
(In both cases, we may restrict the indexing set to
those $\psi \in \PPsiGTF{}\bG{\bS^\sharp}{F\unram}$
for which $\psi(x) = r$.)
These groups and lattices are $\Gal(F\unram/F)$-stable.
Put
$$
G_{x,r} = \bG(F\unram)_{x,r} \cap G
\quad\text{and}\quad
\Lie(G)_{x, r}  = \Lie(\bG)(F\unram)_{x, r} \cap \Lie(G).
$$
We have that $G_{x,0} = G_x$\,.
These definitions extend in an obvious fashion to the case
$r = \infty$, giving $G_{x, r} = \sset 1$ and
$\Lie(G)_{x, r} = \sset 0$.

\begin{rk}
\label{rem:actually-product}
By Proposition 6.4.48 of \cite{bruhat-tits:reductive-groups-1}
and Lemma \ref{lem:complete-subfield},
the multiplication map
$$
\bT(F\unram)_r
\times
\prod_{ \substack{\psi\in \PPsiGTF{}\bG{\bS^\sharp}{F\unram} \\
		\psi(x) = r} }
\lsub{F\unram}U_\psi 
\to \bG(F\unram)_{x, r}
$$
(the product taken in any order) is a bijection when $r>0$.
\end{rk}

For a fixed $x \in \BB(\bG, F)$,
$(\Lie(G)_{x, r})_{\substack{r \in \tR_{\ge 0} \\ r < \infty}}$
and $(G_{x,r})_{\substack{r\in\tR_{\geq 0} \\ r < \infty}}$
are filtrations of
$\Lie(G)_{x, 0}$ and $G_{x, 0}$
by normal lattices
and open normal subgroups, respectively.
(The indexings of these filtrations depend
on our choice of valuation $\ord$.)
\indexmem{G_r}
Put $G_r = \bigcup_{x\in\BB(\bG,F)} G_{x,r}$
and
$\Lie(G)_r
= \bigcup_{x \in \BB(\bG, F)} \Lie(G)_{x, r}$
for $r \in \tR_{\ge 0}$ or $r \in \tR$, as appropriate, with
$r < \infty$.
Put
$G_\infty = \bigcap_{\substack{r \in \tR_{\ge 0} \\ r < \infty}} G_r$
and
$\Lie(G)_\infty = \bigcap_{\substack{r \in \tR \\ r < \infty}} \Lie(G)_r$\,.
(These sets are related to the sets of unipotent elements in
$G$ and nilpotent elements in $\Lie(G)$, respectively.  See
\cite{adler-debacker:bt-lie}*{\S\S 2.5, 3.7.1}.)

\indexmem{\depth}\indexmem{\depth_x}
For $x \in \BB(\bG, F)$,
there are $(\R_{\ge 0} \cup \sset\infty)$-valued functions
$\depth_x$ on $G_{x, 0}$
and $\depth$ on $G_0$\,,
given by
$\depth_x(g)
= \max\set{r \in \R_{\ge 0} \cup \sset\infty}
	{g \in G_{x, r}}$
and
$\depth(g) = \max\set{r \in \R_{\ge 0} \cup \sset\infty}{g \in G_r}$
for all appropriate $g$.
Similar (but $(\R \cup \sset\infty)$-valued) functions are
defined on $\Lie(G)$.
If necessary, we will denote these functions by
$\depth_{G,x}$ and $\depth_G$ to indicate the dependence on $G$.

If a group \mc G has a filtration $(\mc G_i)_{i \in I}$, then
we shall frequently write $\mc G_{i:j}$ in place of
$\mc G_i/\mc G_j$ when $\mc G_j \subseteq \mc G_i$
(even if the quotient is not a group).  For example, we put
$F_{r:t} = F_r/F_t$,
$U_{\varphi_1:\varphi_2} = U_{\varphi_1}/U_{\varphi_2}$, and
$G_{x, r:t} = G_{x, r}/G_{x, t}$
for $r \le t$ (and $r \ge 0$, in the last case) and for
affine $F$-roots $\varphi_1$ and $\varphi_2$ such that
$\dot\varphi_1 = \dot\varphi_2$ and $\varphi_1 \le \varphi_2$.

\begin{dn}
\label{defn:G_x}
\indexmem{\ms G_x}
There is a (not necessarily connected) reductive \ff-group
$\ms G_x$ such that
$\stab_{\bG(E)}(x)/\bG(E)_{x, 0+}
= \ms G_x(\ff_E)$ whenever $E/F$ is an unramified extension.
%%Why it's a reductive group.
%We have that $\stab_G(x) = \mf G_x(\mc O)$ for some smooth
%\mc O-group scheme $\mf G_x = \mf G^\dag_x$,
%and $(\mf G_x)\conn(\mc O) = \mf G_x\conn(\mc O)$ maps onto
%the set of \ff-rational points of the
%special fiber of $\mf G_x\conn$.
%The unipotent radical of the special fiber is the
%intersection of the unipotent radicals of its Borel
%subgroups.  A Borel subgroup of the special fiber is the
%image there of an Iwahori subgroup corresponding to a
%chamber with $x$ in its closure.  The unipotent radical of
%such a subgroup is the image of the pro-unipotent radical of
%the Iwahori.  Since the intersection of the pro-unipotent
%radicals of Iwahoris corresponding to chambers with $x$ in
%their closure is just $G_{x, 0+}$ (by Remark
%\ref{rem:actually-product}, say), we have that
%$\stab_G(x)/G_{x, 0+}$ maps onto the set of \ff-rational
%points of a reductive group (the maximal reductive quotient
%of the special fiber of $\hat\mf G_x$).
If necessary, we will write $\ms G_x^F$ instead of just
$\ms G_x$ to indicate the dependence on the field $F$.
\end{dn}

In the notation of Definition \ref{defn:G_x},
$G_{x, 0:0+}$ is the group of \ff-rational points of $\ms G_x\conn$.
(In the notation of \cite{tits:corvallis}*{\S 3},
$\ms G_x = \ol G_x/R\unip(\ol G_x)$ and
$\ms G_x\conn = \ol G{}_x\red$.
Most authors use $\ms G_x$ to denote what we are calling
$\ms G_x\conn$.)

If \tbG is $F$-anisotropic, then
$\Lie(G)_{x, r} = \Lie(G)_r$
and
$G_{x,r} = G_r$
for all $x\in \BB(\bG,F)$ and $r \in \tR$
(respectively, $r \in \tR_{\ge 0}$).
In particular,
$G$ has a canonical filtration.
We have already seen what this filtration looks like
when $\bG$ is a torus.
In any case, if \bS is a maximal $F$-split torus in \bG,
then $C_\bG(\bS)$ is $F$-anisotropic modulo its center.
For any constant
\tR-valued function $\varphi = r$ on $\AA(\bS, F)$,
we put
$\mf u_\varphi = \Lie(C_\bG(\bS))(E)_r$
and (if $r \ge 0$) $U_\varphi = C_\bG(\bS)(E)_r$\,.
If $\varphi$ is an affine function on $\AA(\bS, F)$ whose
gradient does not belong to $\wtilde\Phi(\bG, \bS)$, then we put
$U_\varphi = U_{\varphi+} = \sset 1$ and
$\mf u_\varphi = \mf u_{\varphi+} = \sset 0$.

\subsection{Filtrations and descent}
\label{sec:notation_filt-and-descent}

In this section, we gather together some results,
most of which are well known.
Suppose
\begin{itemize}
\item $(x, r) \in \BB(\bG, F) \times \tR_{\ge 0}$,
\item $E/F$ is
a discretely valued algebraic extension,
and
\item \bM is an $F$-Levi subgroup of \bG.
\end{itemize}
(%
\indexme{Levi@$F$-Levi subgroup}%
\indexme{Levi@$E$-Levi $F$-subgroup}%
\label{defn:levi}%
%	%%Just for the page number.
Throughout, we will call a subgroup of \bG an
\emph{$F$-Levi subgroup} if it is a Levi component of some
parabolic $F$-subgroup of \bG.
An $F$-subgroup of \bG which is a Levi component of some
parabolic $E$-subgroup of \bG (for some extension $E/F$)
will be called an \emph{$E$-Levi $F$-subgroup}.)

Recall that we have defined functions $\depth$ and $\depth_x$
on certain subsets of $G$.
Denote by $\depth_M$ and $\depth_{M, x}$ the corresponding
functions defined on the corresponding subsets of $M$.
For any discretely valued algebraic extension $K/F$,
denote by $\depth^K$ and $\depth_x^K$
the corresponding functions defined on the corresponding subsets of $\bG(K)$.

\begin{lemma}
\label{lem:torus-field-descent}
Let \bT be an $F$-torus.
When $E/F$ is tame and $r > 0$,
$T \cap \bT(E)_r = T_r$\,.
When $E/F$ is separable, $T_0 \cap \bT(E)_r = T_r$\,.
In general, there is some finite Galois extension $K/E$ such
that
$T_0 \cap \bT(K)_r = T_r$\,.
\end{lemma}

\begin{proof}
For any extension $K/F$, we have by the definition of the
filtration that
$T_0 \cap \bT(K)_r \subseteq T_r$\,.
If $T_0 \subseteq \bT(K)_0$\,, then also
$T_r \subseteq \bT(K)_r$\,.
If $E/F$ is separable, then fix an element $t \in T_0$\,.
By Lemma \ref{lem:complete-subfield}, there is a complete
subfield $F'$ of $F$ such that
\begin{itemize}
\item $F/F'$ is unramified,
\item \bT is defined over $F'$,
and
\item $t \in \bT(F')$.
\end{itemize}
In particular,
$t \in \bT(F)_0^{\Gal(F/F')}$, which, by definition,
is $\bT(F')_0$\,.
By another application of Lemma \ref{lem:complete-subfield},
there is a finite subextension $E'/F'$ of $E/F'$ such that
$E/E'$ is unramified.  By Lemma 2.1.2 of
\cite{adler-debacker:mk-theory}, $t \in \bT(E')_0$\,.
Since $E/E'$ is unramified,
$\bT(E')_0 \subseteq \bT(E)_0$ by definition.
Thus, $t \in \bT(E)_0$\,; so $T_0 \subseteq \bT(E)_0$\,.
If \bT is $K$-split, then
$\bT(K)_0 = \bT(K)\subb$\,, hence contains $T_0$\,.
The third statement follows.

If $\chr \ff = 0$, then \bT and $E/F$ are tame, so
Proposition 4.7.2 of \cite{yu:models} gives the first statement
in this case.
The first statement will follow in general once we know that
$T \cap \bT(E)_{0+} \subseteq T_0$
when $E/F$ is tame and $p := \chr \ff > 0$.
Since $\bT(E)_{0+} \subseteq \bT(\wtilde E)_{0+}$\,,
where $\wtilde E/F$ is the Galois closure of $E/F$, we may,
and hence do, assume that $E/F$ is Galois.
Since $T_0 = \bT(F\unram)_0^{\Gal(F\unram/F)}$ by definition,
it suffices to assume that
$F = F\unram$.
In the notation of
\cite{kottwitz:isocrystals-2}*{\S 7.3},
with $L = F$ and $L' = E$
(so, in particular, $\beta$ is
the inclusion of $T$ in $\bT(E)$),
we have by (7.3.2) of loc.\ cit.\ that
$$
\alpha(w_{\bT(E)}(\beta(t)))
= \alpha(N(w_T(t)))
= [E : F]w_T(t)
$$
for $t \in T$.
(\cite{kottwitz:isocrystals-2} works over the completion of the
maximal unramified extension of a $p$-adic field; but, as in
\cite{rapoport:T1-is-T0}, the same reasoning works for any
Henselian field with algebraically closed residue field,
such as $F$.)
By Lemma 2.3 of \cite{rapoport:T1-is-T0},
$T_0 = \ker w_T$
and $\bT(E)_0 = \ker w_{\bT(E)}$,
so $\gamma^{\indx E F} \in T_0$
whenever $\gamma \in T \cap \bT(E)_0$\,.
If further $\gamma \in \bT(E)_{0+}$\,, then, since
the sequence $(\gamma^{p^n})_{n \in \Z_{> 0}}$ in $T$ tends
to the identity element, and
$T_0$ is an open subgroup of $T$, there is
some positive integer $N$ such that $\gamma^{p^N} \in T_0$\,.
Since $\indx E F$ and $p$ are coprime, we have that
$\gamma \in T_0$\,.
\end{proof}

\begin{lemma}
\label{lem:field-descent}
If $E/F$ is separable, then
$\bG(E)_{x,r} \cap G\supseteq G_{x,r}$\,.
If $E/F$ is unramified,
or $r > 0$ and $E/F$ is tame,
then we have equality.
\end{lemma}

\begin{proof}
By Lemma \ref{lem:complete-subfield} and the definition of
the filtration,
$\bG(E)_{x, r} \cap G$ is the union of all the subgroups of the
form $\bG(L)_{x, r} \cap \bG(F')$\,, where
$F'$ is a complete subfield of $F$,
$L$ is a finite extension of $F'$ contained in $E$,
and
$E/L$ is unramified.
Thus, we may, and hence do, assume that $E/F$ is
finite, and that $F$ (hence $E$) is complete.
For $r = 0$, the containment is Lemma 2.1.2 of
\cite{adler-debacker:mk-theory}.
For $r > 0$, Proposition 1.4.1 of
\cite{adler:thesis} shows that the containment always holds.
For $E/F$ unramified, equality follows from the
definition of the filtration.
For $r > 0$ and $E/F$ tame, equality is proven in
Proposition 1.4.2 of \cite{adler:thesis} 
(using the conclusion of Lemma
\ref{lem:torus-field-descent}).
\end{proof}

\begin{lemma}
\label{lem:levi-descent}
$G_r \cap M = M_r$\,,
and, if
$x \in \BB(\bM, F)$,
then $G_{x, r} \cap M = M_{x, r}$\,.
\end{lemma}

In the statement of the lemma, we have regarded $\BB(\bM, F)$ 
as a subset of $\BB(\bG, F)$, as in \cite{debacker:bt-group}
(or Definition \ref{defn:compatibly-filtered} below).

\begin{proof}
Once we know the results for all $r < \infty$, the results
for $r = \infty$ will follow.  Therefore, we may, and hence
do, assume that $r < \infty$.

It is clear (from the second part of the lemma) that $M_r \subseteq G_r \cap M$.
Choose an element $g \in G_r \cap M$.
Then there is some $y \in \BB(\bG, F)$ such that
$g \in G_{y, r}$\,.
By Lemma \ref{lem:complete-subfield}, there is a complete
subfield $F'$ of $F$ such that
\begin{itemize}
\item \bG and $g$ are defined over $F'$,
\item \bM is an $F'$-Levi subgroup of \bG,
\item $y \in \BB(\bG, F')$,
and
\item $F/F'$ is unramified.
\end{itemize}
Then, by Lemma \ref{lem:field-descent},
$g \in \bG(F')_{y, r} \cap \bM(F')$.
By Theorem 4.1.5 of \cite{debacker:bt-group},
$g \in \bM(F')_{y, r}$\,.
By another application of Lemma \ref{lem:field-descent},
$g \in M_r$\,.

By a similar argument, we may, and hence do, assume for the
second statement that $F$ is complete.
For $r = 0$, the statement is Lemma 4.2.2 of
\cite{debacker:bt-group}.
For $r > 0$, it is Theorem 4.2 of
\cite{moy-prasad:jacquet}.
\end{proof}

\begin{lemma}
\label{lem:domain-field-ascent}
When $E/F$ is separable,
$\bG(E)_r \cap G \supseteq G_r$\,.
If $E/F$ is unramified, or $r > 0$ and $E/F$ is tame, then
we have equality.
\end{lemma}

\begin{proof}
Once we know the result for all $r < \infty$, it follows
for $r = \infty$; so we assume throughout that $r < \infty$.
In this case, 
the containment follows from Lemma \ref{lem:field-descent}.

Suppose that $\gamma \in \bG(E)_r \cap G$.
Then $\gamma \in \bG(\wtilde E)_r \cap G$, where
$\wtilde E$ is the Galois closure of $E$ over $F$.
Since $\wtilde E/F$ is unramified (respectively, tame) if
$E/F$ is, we may, and hence do, assume that $E/F$ is Galois.
Choose $x \in \BB(\bG, E)$ with $\gamma \in \bG(E)_{x, r}$\,.
Then $\gamma \in \bG(E)_{\sigma x, r}$ for
$\sigma \in \Gal(E/F)$, so $\gamma \in \bG(E)_{\tilde x, r}$\,,
where $\tilde x$ is the center of mass of the
$\Gal(E/F)$-orbit of $x$.
Suppose that $E/F$ is tame.
Then $\tilde x \in \BB(\bG, E)^{\Gal(E/F)} = \BB(\bG, F)$.
If $E/F$ is unramified or $r > 0$, then we have by another
application of Lemma \ref{lem:field-descent} that
$\gamma \in G_{\tilde x, r}$\,.
That is, $\bG(E)_r \cap G \subseteq G_r$\,.
The reverse containment being obvious, we have equality, as
desired.
\end{proof}

\begin{lemma}
\label{lem:unipotent}
If \bS is a maximal $F$-split torus in \bG and
$x, y \in \AA(\bS, F)$, then, for any
$g \in G_{x, r} \cap G_{y, r+}$\,,
we have that
$g G_{x, r+} \subseteq G_{z, r+}$
for $z \in (x, y)$ sufficiently close to $x$;
and there is a
parabolic $F$-subgroup \bP of \bG containing \bS
(depending on $g$)
such that $g \in R\unip(P)\cdot G_{x, r+}$\,,
where $R\unip(\bP)$ is the unipotent radical of \bP.
\end{lemma}

Here, $(x, y)$ denotes the open line segment
between $x$ and $y$ in $\AA(\bS,F)$.

\begin{proof}
By Lemma \ref{lem:complete-subfield}, there is a complete
subfield $F'$ of $F$ such that
\begin{itemize}
\item \bG and \bS are defined over $F'$,
\item \bS is $F'$-split,
\item $g$ is defined over $F'$,
and
\item $F/F'$ is unramified.
\end{itemize}
By Lemma \ref{lem:field-descent},
$g \in \bG(F')_{x, r} \cap \bG(F')_{y, r+}$\,.
Thus, we may, and hence do, assume that $F$ is complete.

The first statement is proved in the proof of Corollary 3.7.10 of
\cite{adler-debacker:bt-lie}.
The second is proved in the ``$\subset$'' part of the proof of
Lemma 4.3.2 of \cite{debacker:bt-group} for $r = 0$, and
in the proof of Lemma 3.7.6 of
\cite{adler-debacker:bt-lie} for $r > 0$.
\end{proof}

\begin{lm}
\label{lem:stab-deep}
$\stab_G(\ox) \cap G_0 = G_{x, 0}$\,.
\end{lm}

\begin{proof}
Let $g\in \stab_G(\ox) \cap G_0$\,.
By definition, $G$ acts by translations on the factor
$V_F(Z(\bG))$ in
$\BB(\bG, F) = \BB\red(\bG, F) \times V_F(Z(\bG))$.
Since $G_0$ consists of bounded elements, it acts trivially
on $V_F(Z(\bG))$,
so in fact $g \in \stab_G(x)$.
Further, $g\in G_{y,0}$ for some $y\in \BB(\bG,F)$.
From Lemma \ref{lem:complete-subfield},
we may pick a complete subfield $F'$ of $F$
such that
\begin{itemize}
\item \bG is defined over $F'$,
\item $F/F'$ is unramified,
\item $x,y\in \BB(\bG,F')$,
and
\item $g\in \bG(F')$.
\end{itemize}
Then $g\in \stab_{\bG(F')}(x)$, and,
from Lemma \ref{lem:field-descent},
$g\in \bG(F')_{y,0}\subseteq \bG(F')_0$\,.
By Lemma 4.2.1 of \cite{debacker:bt-group},
$\bG(F')_0 \cap \stab_{\bG(F')}(x) = \bG(F')_{x, 0}
\subseteq G_{x,0}$\,.
The reverse containment, hence equality, is obvious.
\end{proof}

\subsection{The effect of the center on depth and degeneracy}

Fix $(x, r) \in \BB(\bG, F) \times \tR_{\ge 0}$.

\begin{lem}
\label{lem:depth-in-center}
$Z(G) \cap G_{x, r} = Z(G) \cap G_{y, r}$
for any $y \in \BB(\bG, F)$.
\end{lem}

\begin{proof}
By Th\'eor\`eme 7.4.18 of
\cite{bruhat-tits:reductive-groups-1}, there exists a
maximal $F$-split torus \bS such that
$x, y \in \AA(\bS, F)$.
By Lemma \ref{lem:levi-descent},
$$
Z(G) \cap G_{x, r}
= Z(G) \cap (C_G(S) \cap G_{x, r})
= Z(G) \cap C_G(S)_r\,,
$$
and similarly for $Z(G) \cap G_{y, r}$\,.
\end{proof}

\begin{lm}
\label{lem:depth-mod-center}
If $z \in Z(G)$, then
$z G_{x, r} \cap G_r \ne \emptyset$ if and only if
$z \in G_{x, r}$\,.
\end{lm}

\begin{proof}
The `if' part is clear.

For the `only if' part, suppose that
$z G_{x, r} \cap G_r \ne \emptyset$.
Since Lemma \ref{lem:stab-deep} gives
$z G_{x, r} \cap G_0 \subseteq \stab_G(\ox) \cap G_0 = G_{x, 0}$\,,
we have $z \in G_{x, 0}$\,.
Put $t = \depth_x(z)$.
By Lemma \ref{lem:depth-in-center},
$Z(G)\cap G_{y,t+}$ is independent
of $y \in \BB(\bG, F)$,
so we have that $t=\depth(z)$.
If $t \geq r$, then we are done,
so assume $t<r$.
From
Lemma \ref{lem:unipotent},
all elements of $z G_{x, r} \subseteq zG_{x,t+}$ have depth $t$,
so
$zG_{x,r} \cap G_r$ is empty, which is a contradiction.
\end{proof}

\begin{cor}
\label{cor:depth-mod-center}
$Z(G)G_{x, r} \cap G_r = G_{x, r}$\,.
\end{cor}

\begin{lm}
\label{lem:degen-mod-center}
If $g \in G_{x, r}$
and $g G_{x, r+} \cap Z(G)G_{r+} \ne \emptyset$,
then $g \in Z(G)G_{r+}$\,.
\end{lm}

\begin{proof}
Suppose that $z \in Z(G)$ is such that
$g G_{x, r+} \cap z G_{r+} \ne \emptyset$.
Then
$\emptyset \neq z\inv g G_{x,r+} \cap G_{r+}
\subseteq z\inv G_{x,r} \cap G_r$\,,
so, by Lemma \ref{lem:depth-mod-center},
$z\in G_{x,r}$\,.
Therefore, $z\inv g \in G_{x,r}$\,, so
Lemma \ref{lem:unipotent}
implies that $z^{-1}g \in G_{r+}$\,.
\end{proof}

\section{Tameness and compatible filtration}
\label{sec:tame-and-compatible}

\begin{rk}
\label{rem:G-image}
Let \mo N be a closed normal $F$-subgroup of \bG.
By Theorem AG.17.3 and Proposition 6.5 of \cite{borel:linear}, the map
$\bG \to \bG/\mo N$ is a submersion of varieties.
Fix an open subset $U$ of $G$.
By Lemma \ref{lem:complete-subfield},
$U = \varinjlim (U \cap \bG(F'))$, the limit taken over all
complete subfields $F'$ of $F$ over which \bG is defined.
For such a field $F'$,
the induced map $\bG(F') \to (\bG/\mo N)(F')$ is a submersion of
analytic manifolds, hence, by Theorem II.III.10(2)(2) of
\cite{serre:lie-alg+lie-gp}*{p.~85},
an open map.
Thus the image $\ol{U \cap \bG(F')}$ of $U \cap \bG(F')$
in $(\bG/\mo N)(F')$ is open there;
so the image $\ol U = \varinjlim \ol{U \cap \bG(F')}$ of $U$
in $(\bG/\mo N)(F) = \varinjlim (\bG/\mo N)(F')$ is open there.
That is, $G \to (\bG/\mo N)(F)$ is an open map.

Denote by \ol G the image of $G$ in $(\bG/\mo N)(F)$, with the subspace
topology.  Then \ol G is
\begin{itemize}
\item open in $(\bG/\mo N)(F)$;
\item closed in $(\bG/\mo N)(F)$;
and
\item homeomorphic to $G/N$.
\end{itemize}
If $F$ is complete (respectively, locally compact),
then so are $G$ and (\bG/\mo N)(F), hence also \ol G.
\end{rk}

\begin{defn}
\label{defn:tame}
Suppose that \mo N is a closed normal $F$-subgroup of \bG.
\indexme{tame@$F$-tame!group}
\indexme{tame@$F$-tame!modulo a subgroup!group}
Say that \bG is \emph{$F$-tame}
(respectively, \emph{$F$-tame modulo \mo N{}})
if it contains a maximal $F$-torus \bT such that
\bT (respectively, $\bT/(\mo N \cap \bT)$) splits over a
tame extension of $F$.
If \bS is an $F$-torus in \bG, we will often say that \bS is
\emph{$F$-tame modulo center} instead of $F$-tame modulo
$Z(\bG)$ if \bG is understood from the context.

\indexme{tame@$F$-tame!element}
\indexme{tame@$F$-tame!modulo a subgroup!element}
Say that an element $\gamma \in G$ is \emph{$F$-tame in \bG{}}
(respectively, \emph{$F$-tame in \bG modulo \mo N{}})
if it is semisimple and $C_\bG(\gamma)\conn$ is $F$-tame
(respectively, $F$-tame modulo \mo N).
We will omit ``in \bG'',
and say that $\gamma$ is \emph{$F$-tame modulo center}
instead of $F$-tame modulo $Z(\bG)\conn$,
if \bG is understood from the context.

We will frequently say \emph{tame} instead of $F$-tame if
$F$ is understood.
\end{defn}

Note that a torus which is $F$-tame modulo a group \mo N need not
actually contain \mo N.
If $\mo N'$ is another normal subgroup such that
$\mo N' \supseteq \mo N$ and
$\bG/\mo N' \to \bG/\mo N$ is a central isogeny, then, by
Corollaire 1.9(a) of \cite{borel-tits:reductive-groups}, we
have that an $F$-torus or element of $G$
is $F$-tame modulo \mo N if and only if it
is $F$-tame modulo $\mo N'$.
In particular, $F$-tameness modulo $Z(\bG)$ is the same as
$F$-tameness modulo $Z(\bG)\conn$.

\begin{defn}
\label{defn:compatibly-filtered}
Suppose that $i : \bH \to \bG$ is a closed embedding of
reductive $F$-groups.
Given a discretely valued algebraic extension $E/F$,
we say that a map $i_E : \BB(\bH, E) \to \BB(\bG, E)$
\indexme{preserves filtrations}
\emph{preserves filtrations over $F$} if
\begin{enumerate}
\item
$i_E$ is an $\Aut(E/F)$-equivariant isometry;
\item
$i_E(h x) = i(h)i_E(x)$ for $h \in \bH(E)$ and
$x \in \BB(\bH, E)$;
\item\label{defn:compatibly-filtered_affine}
for every maximal $E\tame$-split $E$-torus
$\bS_\bH \subseteq \bH$,
there is a maximal $E\tame$-split $E$-torus
$\bS_\bG \subseteq \bG$
such that $i(\bS_\bH) \subseteq \bS_\bG$ and
$i_E$ restricts to an affine injection of 
$\AA(\bS_\bH, E\tame)^{\Gal(E\tame/E)} \subseteq \BB(\bH, E)$
into
$\AA(\bS_\bG, E\tame)^{\Gal(E\tame/E)} \subseteq \BB(\bG, E)$;
and
\item
for all $(x, r) \in \BB(\bH,E) \times \tR_{> 0}$,
$\bH(E)_{x,r}$ is the preimage in $\bH(E)$ of $\bG(E)_{i_E(x),r}$\,.
\end{enumerate}

We say that an $F$-embedding
$i : \bH \to \bG$ of reductive $F$-groups is
\indexme{filtration-preserving map}
\emph{filtration preserving over $F$} if there is a system
$$
(i_E \colon  \BB(\bH, E) \to \BB(\bG, E))
	_{\text{$E/F$ a discretely valued
		separable
		extension}}
$$
of filtration-preserving embeddings of buildings such that, for all
pairs $E/F$ and $E'/F$ of discretely valued separable extensions
with $E' \subseteq E$,
$i_E\bigr|_{\BB(\bH, E')} = i_{E'}$.
If \bH is a reductive $F$-subgroup of \bG such that the inclusion
$\bH \hookrightarrow \bG$ is filtration preserving over $F$,
then we say that \bH is a \emph{compatibly filtered $F$-subgroup} of \bG.
\indexme{compatibly filtered $F$-subgroup}
\end{defn}

When such a system
exists, we will often use it
to identify the buildings of $\bH$ (over various fields)
with subsets of those of $\bG$,
using the same letter for an element of $\BB(\bH, F)$
and its image in $\BB(\bG, F)$.

For every discretely valued separable extension $E/F$,
there is a canonical $\bH(E)$- and $\Aut(E/F)$-equivariant isomorphism
$i_E : \BB(\bH\conn, E) \to \BB(\bH, E)$.
We have
$\bH(E)_{i_E(x), r} := \bH\conn(E)_{x, r}$ for
$(x, r) \in \BB(\bH\conn, E) \times \tR_{\ge 0}$,
and the maps $i_E$ are compatible in the sense of Definition
\ref{defn:compatibly-filtered}.
Thus $\bH\conn$ is a compatibly filtered $F$-subgroup of \bG if
\bH is.

Note that, if \bH is a compatibly filtered $F$-subgroup of \bG
and \bL is a compatibly filtered $F$-subgroup of \bH,
then \bL is a compatibly filtered $F$-subgroup of \bG.

\begin{lemma}
\label{lem:group-descent}
Suppose that \bH is a compatibly filtered $F$-subgroup of \bG
and $x \in \BB(\bH, F)$.
Then $H_{x, 0} \subseteq G_{x, 0}$\,.
\end{lemma}

\begin{proof}
By Definition \ref{defn:compatibly-filtered},
$\bH(F\unram)_{x, 0+}
= \stab_{\bH(F\unram)}(x) \cap \bG(F\unram)_{x, 0+}$\,.
Thus
\begin{multline*}
\ms H_x(\ol\ff)
= \stab_{\bH(F\unram)}(x)/\bH(F\unram)_{x, 0+} \\
\subseteq
\stab_{\bG(F\unram)}(x)/\bG(F\unram)_{x, 0+}
= \ms G_x(\ol\ff),
\end{multline*}
so
$$
\bH(F\unram)_{x, 0:0+} = \ms H_x\conn(\ol\ff)
\subseteq \ms G_x\conn(\ol\ff) = \bG(F\unram)_{x, 0:0+}\,.
$$
In particular, $\bH(F\unram)_{x, 0} \subseteq \bG(F\unram)_{x, 0}$\,.
Since $H_{x, 0} = \bH(F\unram)_{x, 0}^{\Gal(F\unram/F)}$\,,
and similarly for $G_{x, 0}$\,, we have the desired
containment.
\end{proof}

\begin{lm}
\label{lem:ratl-maxl-torus}
If $E/F\unram$ is an algebraic extension,
then any $E$-split $F$-torus in \bG is contained in some
maximal $E$-split torus in \bG which is defined over $F$.
\end{lm}

We emphasize that the torus in question is maximal among
$E$-split tori, not just among $E$-split $F$-tori.

\begin{proof}
Let $E'/F$ be the maximal Galois subextension of $E/F$.
Since $F\unram/F$ is Galois, $E'$ contains $F\unram$.
Since an $F$-torus is $E$-split if and only if it is
$E'$-split, we may, and hence do, assume, upon replacing $E$
by $E'$, that $E/F$ is Galois.

Fix an $E$-split $F$-torus $\bS'$.
Since the $E$-split ranks of \bG and $C_\bG(\bS')$ are the
same, upon replacing \bG by $C_\bG(\bS')$,
we may, and hence do, assume
that $\bS'$ is central in \bG.
Denote by $\bS'\lsplit$ the maximal $F$-split subtorus of $\bS'$.
By Corollaire 5.1.12 of
\cite{bruhat-tits:reductive-groups-2},
there exists a maximal $F\unram$-split $F$-torus
$\bS^{\prime\,\sharp}$ in \bG containing $\bS'\lsplit$.
Since $\bS^{\prime\,\sharp}$ is $E$-split, it is contained in a
maximal $E$-split torus $\bS^\sharp$ in \bG.  In particular,
$\bS^\sharp$ is a maximal $E$-split torus in
$C_\bG(\bS^{\prime\,\sharp})$.
Since \bG is $F\unram$-quasisplit,
$C_\bG(\bS^{\prime\,\sharp})$ is a maximal $F$-torus.
In particular, it has a unique maximal $E$-split
subtorus, namely $\bS^\sharp$; so $\bS^\sharp$ is $\Gal(E/F)$-stable,
hence an $F$-torus.
Finally, since $\bS'$ is $E$-split, there is an element
$g \in \bG(E)$ such that $\lsup g\bS' \subseteq \bS^\sharp$.  Since
$\bS'$ is central in \bG, we have that $\lsup g\bS' = \bS'$.
\end{proof}

The following proposition is an example of a situation in which
a subgroup of \bG is compatibly filtered over $F$.

\begin{pn}
\label{prop:compatibly-filtered-tame-rank}
Suppose that \bH is a connected reductive $F$-subgroup of \bG such
that the absolute ranks of \bH and \bG are the same, and
there is some tame extension
$L/F\unram$ such that
\bH and \bG have the same $L$-split rank.
Then \bH is a
compatibly filtered $F$-subgroup of \bG.
Moreover, for any discretely valued separable extension $E/F$,
the image of $\BB(\bH, E)$ in $\BB(\bG, E)$ is
independent of the choice of filtration-preserving
embeddings over $F$.
\end{pn}

\begin{rk}
We will actually prove a stronger uniqueness statement;
namely, that, if $(i_E)$
and $(i'_E)$ are two compatible systems as in Definition
\ref{defn:compatibly-filtered}, then there is some
$\lambda \in \bX_*^F(Z(\bH)) \otimes_\Z \R$ such that, for
all discretely valued separable extensions $E/F$,
$i_E(x + \lambda) = i'_E(x)$ for $x \in \BB(\bH, E)$.
\end{rk}

\begin{proof}
By Lemma \ref{lem:ratl-maxl-torus}, there exists a maximal
$L$-split torus in \bH which is defined over $F$.  Upon
replacing $L$ by the compositum with $F\unram$ of the
splitting field of such a torus, we may, and hence do,
assume that $L$ is discretely valued and strictly Henselian,
hence that \bG and \bH are $L$-quasisplit.

Let $E/F$ be a discretely valued separable extension.
Assume first that $E \supseteq L$ (hence that $E$ is strictly Henselian);
we will show how to handle arbitrary $E$ later.
Let \bS be a maximal $E$-split torus in \bH.
By Lemma \ref{lem:rank-ascent}, \bS is also a maximal
$E$-split torus in \bG.
Since they have the same underlying affine space,
there is a natural affine isometry $i_{E, \bS}$ of
$\lsub\bH\AA(\bS, E)$ with $\lsub\bG\AA(\bS, E)$.
Since $N_\bH(\bS)(E) \subseteq N_\bG(\bS)(E)$
and $\Aut(E/F)$ act on $\AA(\bS, E)$ by
affine transformations, this identification is
$N_\bH(\bS)(E)$- and $\Aut(E/F)$-equivariant.
If $K/E$ is a further discretely valued separable extension
and $\bS^\sharp$ is a maximal $K$-split torus containing
\bS, then $i_{K, \bS^\sharp}$ restricts to $i_{E, \bS}$.

We recall from \cite{tits:corvallis}*{\S 1.4} (where they
are denoted by $X_\varphi$) the definition of the affine
root subgroups $U_\varphi$.
This will be useful later in the proof.
Suppose
$a \in \Phi(\bH, \bS) \subseteq \Phi(\bG, \bS)$.
From Lemma \ref{lem:no-root-ambience},
the root subgroups of \bG and \bH associated to $a$ are the
same, so we may write $\bU_a$ without ambiguity.
Suppose $\varphi$ is an affine function on $\AA(\bS, E)$
such that $\dot\varphi = a$.
If $a = 0$, then recall that $U_\varphi = \bT(E)_r$\,,
where \bT is the unique maximal torus in \bH (or \bG)
containing \bS and
$r$ is the value of the constant function $\varphi$.
This is obviously independent of whether the ambient group is
\bH or \bG.
If $a \in \Phi(\bG, \bS)$, then,
for $u\in \bU_a(E) \smallsetminus \sset 1$,
the set
$\bU_{-a}(\ol F)u\bU_{-a}(\ol F)
\cap N_\bG(\bS)(\ol F)$
has only one element, say $n_{\bG, \bS}(u)$
(which actually lies in
$\bU_{-a}(E)u\bU_{-a}(E) \cap N_\bG(\bS)(E)$).
This element uniquely determines an affine function
$\varphi_{\bG, \bS}(u)$ on $\AA(\bS, E)$
such that $\dot\varphi_{\bG, \bS}(u) = a$ and
$\varphi_{\bG, \bS}(u)$ vanishes on the hyperplane
fixed by $n_{\bG, \bS}(u)$.
(Note that, since the action of $N_\bG(\bS)(E)$ on
$\AA(\bS, E)$ depends on our choice of valuation $\ord$,
so does this hyperplane.)
Then we put
$U_\varphi^\bG =
\sset 1
\cup
\set{u \in \bU_a(E) \smallsetminus \sset 1}
	{\varphi_{\bG, \bS}(u) \geq \varphi}$.
We define $U_\varphi^\bH$ (and the associated notation)
similarly.  Since
$n_{\bH, \bS}(u) \in \bU_{-a}(\ol F)u\bU_{-a}(\ol F)
			\cap N_\bG(\bS)(\ol F)$,
in fact $n_{\bH, \bS}(u) = n_{\bG, \bS}(u)$,
so $\varphi_{\bH, \bS}(u) = \varphi_{\bG, \bS}(u)$
for $u \in \bU_a(E) \smallsetminus \sset 1$.  Thus,
$U_\varphi^\bG = U_\varphi^\bH$, and we may use the notation
$U_\varphi$ without ambiguity.

Now choose an element $x_\bH \in \BB(\bH, E)$.
Then there is some $h \in \bH(E)$ such that
$y_\bH := h\inv x_\bH \in \lsub\bH\AA(\bS, E)$.
Put $y_\bG := i_{E, \bS}(y_\bH)$ and
$x_\bG = h y_\bG$.

We claim that $x_\bG$ is independent of the choice of
$h \in \bH(E)$ as above.  Indeed, if also $h' \in \bH(E)$
satisfies
$y'_\bH := h'^{-1}x_\bH \in \lsub\bH\AA(\bS, E)$,
then $h\inv h'$ carries $y'_\bH$ to $y_\bH$.
By Lemma \ref{lem:move-within-apt}, there are
$n \in N_\bH(\bS)(E)$ and
$h_0 \in \bH(E)_{y'_\bH, 0}$ such that
$h\inv h' = n h_0$.
In particular, $n y'_\bH = y_\bH$;
so also $n y'_\bG = y_\bG$, where $y'_\bG = i_{E, \bS}(y'_\bH)$.

Now $\bH(E)_{y'_\bH, 0}$ is generated by those
$U_\varphi = U_\varphi^\bH$
for which
$\varphi \in \PPsiGTF\wtilde\bH\bS E
	\subseteq \PPsiGTF\wtilde\bG\bS E$
and $\varphi(y'_\bH) \ge 0$.
Since $y'_\bH \in \lsub\bH\AA(\bS, E)$ and
$y'_\bG \in \lsub\bG\AA(\bS, E)$ correspond to the same point of
the underlying affine space $\AA(\bS^\sharp, E)$,
every such $\varphi$ also satisfies $\varphi(y'_\bG) \ge 0$,
so $U_\varphi \subseteq \bG(E)_{y'_\bG, 0}$\,.
That is,
$\bH(E)_{y'_\bH, 0} \subseteq \bG(E)_{y'_\bG, 0}$\,.
In particular,
$(h\inv h')y'_\bG = n y'_\bG = y_\bG$, i.e.,
$h'y'_\bG = h y_\bG = x_\bG$.

Thus we may unambiguously put $i_E(x_\bH) = x_\bG$.
If $K/E$ is a further discretely valued separable extension
and $\bS^\sharp$ a maximal $K$-split torus containing \bS,
then
$y_\bH \in \AA(\bS, E) \subseteq \AA(\bS^\sharp, K)$
and $i_{K, \bS^\sharp}(y_\bH) = i_{E, \bS}(y_\bH) = y_\bG$,
so the same process puts $i_K(x_\bH) = h y_\bG = i_E(x_\bH)$.

By Remark \ref{rem:actually-product}, for $r \in \tR_{> 0}$,
$$
\bG(E)_{y_\bG, r}
= \bH(E)_{y_\bH, r}
\times
\prod U_\varphi\,,
$$
the latter product taken over those
$\varphi \in \PPsiGTF\wtilde\bG\bS E
	\smallsetminus \PPsiGTF\wtilde\bH\bS E$
such that $\varphi(y_\bH) = r$.
In particular, clearly
$\bH(E)_{y_\bH, r}
\subseteq \bG(E)_{y_\bG, r} \cap \bH(E)$.
On the other hand, if
$g := h\prod u_\varphi \in \bG(E)_{y_\bG, r}$
belongs to $\bH(E)$, then so does $g\inv h$; and the image
(under the multiplication map $\mult_\bH$ of
Lemma \ref{lem:Zariski}) of
$(g\inv h, (u_\varphi))$ is $1$.  Since the map in question
is injective, we have that $g = h \in \bH(E)_{y_\bH, r}$\,.
Thus $\bH(E)_{y_\bH, r} = \bG(E)_{y_\bG, r} \cap \bH(E)$,
so
\begin{equation}
\tag{$*$}
\bH(E)_{x_\bH, r} = \lsup h\bH(E)_{y_\bH, r}
= \lsup h\bG(E)_{y_\bG, r} \cap \bH(E)
= \bG(E)_{i_E(x_\bH), r} \cap \bH(E).
\end{equation}

By construction, $i_E$ is $\bH(E)$-equivariant,
hence independent of the choice of \bS.  It is easy to see
that it is also $\Aut(E/F)$-equivariant.
Note that $i_E$ is an isometry when restricted to
$\lsub\bH\AA(\bS, E)$.
By Proposition 2.3.1 of
\cite{bruhat-tits:reductive-groups-1},
any pair of points of $\BB(\bH, E)$ may be (simultaneously)
conjugated into $\lsub\bH\AA(\bS, E)$ by a point of
$\bH(E)$.
Thus, the $\bH(E)$-equivariance of $i_E$ implies that it is
an isometry on all of $\BB(\bH, E)$.

Now drop the assumption that $E$ contains $L$,
and assume only that $E/F$ is some discretely valued
separable extension.
Pick a discretely valued tame extension $K/E L$
such that $K/F$ is Galois.
Then we have a $\Gal(K/F)$-equivariant isometry
$i_K : \BB(\bH, K) \to \BB(\bG, K)$,
which restricts to an $\Aut(E/F)$-equivariant isometry
$\BB(\bH, E) = \BB(\bH, K)^{\Gal(K/E)}
\to \BB(\bG, K)^{\Gal(K/E)} = \BB(\bG, E)$,
independent of the choice of $K$.
We define $i_E$ to be this latter map.
By ($*$) and Lemma \ref{lem:field-descent}, for
$(x, r) \in \BB(\bH, E) \times \tR_{> 0}$, we have
$\bH(E)_{x, r} = \bG(E)_{i_E(x), r} \cap \bH(E)$.

Now suppose that $\bS'$ is a maximal $E\tame$-split
$E$-torus in \bH.  By Lemma \ref{lem:ratl-maxl-torus},
$\bS'$ is a maximal $E\tame$-split torus in \bH, hence, by
Lemma \ref{lem:rank-ascent}, in \bG.
Denote by $K/E$ the splitting field of $\bS'$.
Then $i_E$ is a restriction of $i_K$, and the restriction of
$i_K$ to
$\AA(\bS', K)^{\Gal(K/E)} = \AA(\bS', E\tame)^{\Gal(E\tame/E)}$
agrees with the restriction of the affine injection (indeed,
isomorphism) $i_{K, \bS'}$\,.

We have constructed one compatible system $(i_E)$ of
filtration-preserving embeddings over $F$.  Suppose that
$(i'_E)$ is another such system.
Fix a discretely valued separable extension $E/F$.
As above,
we may find a discretely valued tame Galois extension $K/E$ and a $K$-split
$E$-torus $\bS'$ such that $\bS'$ is an $E\tame$-split torus
maximal in \bH and \bG.
We may, and hence do, suppose that $K$ is strictly Henselian.
By Definition
\ref{defn:compatibly-filtered}\eqref{defn:compatibly-filtered_affine},
$i_K(x + \lambda) = i_K(x) + \lambda$ and
$i'_K(x + \lambda) = i'_K(x) + \lambda$ for
$x \in \lsub\bH\AA(\bS', K)$ and
$\lambda \in \bX_*^K(\bS') \otimes_\Z \R$;
and, since the unique $E\tame$-split $E$-torus in
\bG containing $\bS'$ is $\bS'$ itself,
the image under each of $i_K$ and $i'_K$ of
$\lsub\bH\AA(\bS', K)$ is $\lsub\bG\AA(\bS', K)$.
Since $\lsub\bH\AA(\bS', K)$ and $\lsub\bG\AA(\bS', K)$
are affine spaces under
$\bX_*^K(\bS') \otimes_\Z \R$, there is some
$\lambda_E \in \bX_*^K(\bS') \otimes_\Z \R$ such that
$i_K(x + \lambda_E) = i'_K(x)$ for all
$x \in \lsub\bH\AA(\bS', K)$.
Fix, arbitrarily, a point $o \in \lsub\bH\AA(\bS', K)$.
Since $K$ is strictly Henselian and
$\bH(K)_{o + \lambda_E, r}
= \bG(K)_{i_K(o + \lambda_E), r} \cap \bH(K)
= \bG(K)_{i'_K(o), r} \cap \bH(K)
= \bH(K)_{o, r}$
for $r \in \tR_{> 0}$,
Remark \ref{rem:actually-product} shows that
$\varphi(o) = \varphi(o + \lambda_E)$ for all affine
$K$-roots $\varphi$.  That is,
$\lambda_E \in \bX_*^K(Z(\bH)) \otimes_\Z \R$.
Since $i_K$ and $i'_K$ are $\Aut(K/F)$-equivariant,
$\lambda_E \in (\bX_*^K(Z(\bH)) \otimes_\Z \R)^{\Aut(K/F)}
	= \bX_*^F(Z(\bH)) \otimes_\Z \R$.
Since $i_K$ and $i'_K$ are $\bH(K)$-equivariant,
$i_K(x + \lambda_E) = i'_K(x)$ for $x \in \BB(\bH, K)$,
so $i_E(x + \lambda_E) = i'_E(x)$ for $x \in \BB(\bH, E)$.
For $x \in \BB(\bH, F)$, we have that
$\lambda_F = i_F(x) - i'_F(x) = i_E(x) - i'_E(x) = \lambda_E$.
Thus, $\lambda_E$ does not depend on $E$.
\end{proof}

\begin{cor}
\label{cor:compatibly-filtered-tame-rank}
Let \mo D be a subvariety, defined over $F$, of an
$F\tame$-split modulo center torus in \bG, and put
$\bH = C_\bG(\mo D)$.  Then \bH and $\bH\conn$ are
compatibly filtered $F$-subgroups of \bG.
\end{cor}

\begin{proof}
By Corollary 9.2 of \cite{borel:linear}, \bH (hence
%by Proposition 1.2 of loc.\ cit.
also $\bH\conn$) is defined over $F$.
Clearly, $\bH\conn$ and \bG have the same absolute rank.
By Lemma \ref{lem:torus-quotient-product}, they also have
the same $F\tame$-split rank.
Thus $\bH\conn$ satisfies the conditions of
Proposition \ref{prop:compatibly-filtered-tame-rank}, so that 
$\bH\conn$ is a compatibly filtered $F$-subgroup of \bG.
Now it suffices to show that, for a discretely valued
separable extension $E/F$, a point $x \in \BB(\bH, E)$,
and $r \in \tR_{>0}$,
we have that
$\bG(E)_{x, r} \cap \bH(E) \subseteq \bH\conn(E)$.
By Lemma \ref{lem:field-descent}, enlarging $E$ only makes
this statement stronger, so we will do so as necessary.

First, we may, and hence do, assume that \bH is $E$-split.
If $x \in \BB(\bH, E)$, then let \bT be an $E$-split maximal
torus in \bH such that $x \in \AA(\bT, E)$.
Then \bT contains $Z(\bH\conn)$, hence, in particular,
contains \mo D.
By \cite{springer-steinberg:conj}*{\S II.4.1(a)},
the component group of \bH has a set of representatives in
$N_\bH(\bT)(E)$.
Suppose that $h \in N_\bH(\bT)(E)$ is such that
$\bG(E)_{x, r} \cap h\bH\conn(E) \ne \emptyset$.
Let $\bB_\bH$ be a Borel $E$-subgroup of $\bH\conn$ containing \bT.
Then $\lsup h\bB_\bH$ contains $\lsup h\bT = \bT$, so there
is some $n \in N_{\bH\conn}(\bT)(E)$ such that
$\lsup h\bB_\bH = \lsup n\bB_\bH$.
We have that $n\inv h \in N_\bH(\bT)$ and
$(n\inv h)\bH\conn = h\bH\conn$; so, upon replacing $h$ by
$n\inv h$, we may, and hence do, assume that
$\lsup h\bB_\bH = \bB_\bH$.

Now let \bB be a Borel $E$-subgroup of \bG containing
$\bB_\bH$, say with opposite Borel $\bB'$ (with respect to
\bT).
By Remark \ref{rem:actually-product},
$\bG(E)_{x, r} \subseteq \bB(E)\bB'(E)$.
Thus $\bB(E)\bB'(E) \cap h\bH\conn(E) \ne \emptyset$;
hence, \emph{a fortiori},
$\bB\bB' \cap h\bH\conn \ne \emptyset$.
Since $\bB\bB' \cap h\bH\conn$ is a non-empty Zariski open subset
of $h\bH\conn$, it intersects any other such subset; in
particular,
$\bB\bB' \cap h\bB_\bH\bB_\bH' \ne \emptyset$,
where $\bB_\bH'$ is the opposite Borel to $\bB_\bH$ (with
respect to \bT).
Since $h$ normalizes $\bB_\bH$, it follows that
$\bB h\bB' = \bB\bB'$.
Now let $n_0$ be a representative in $N_\bG(\bT)$ of the
long element of the Weyl group of \bG.
Then $\bB' = n_0\bB n_0\inv$, so
$\bB h n_0\bB = \bB n_0\bB$.
By the Bruhat decomposition, this means that $h n_0$
and $n_0$ project to the same element of the Weyl group of
\bG; i.e., $h \in \bT(E)$.
In particular, $h \in \bH\conn(E)$.
That is, the only connected component of \bH intersecting
$\bG(E)_{x, r}$ is the identity component, as desired.
\end{proof}

\begin{rem}
There are other situations where a subgroup $\bH$
of a group $\bG$ is compatibly filtered.
(For example, consider the case where
$\chr \ff \neq 2$,
and $\bH$ is a classical group
embedded in a general linear group $\bG$
in the usual way.)
However, we will not need this fact.
\end{rem}

\begin{lm}
\label{lem:levi-building-descent}
Suppose that
%\mo D is a subvariety, defined over $F$, of an
%$F\tame$-split torus in \bG, and that
%$\bM := C_\bG(\mo D)\conn$
	%%This hypothesis is just uselessly restrictive.
\bM is an $F\tame$-Levi $F$-subgroup of \bG,
and \bH is a compatibly filtered $F$-subgroup of \bG
containing the maximal $F\tame$-split torus in the center of
\bM.
Then
$\BB(\bM, E) \cap \BB(\bH, F) = \BB(\bM \cap \bH, F)$ for any
discretely valued separable extension $E/F$.
\end{lm}

\begin{proof}
Let \bS be the maximal $F\tame$-split torus in the center of
\bM, so that $\bM = C_\bG(\bS)$.
We may, and hence do, replace $F$ by the splitting field of
\bS.
Fix a discretely valued separable extension $E/F$.

Now \bM is an $F$-Levi subgroup of \bG, and
$\bM \cap \bH = C_\bH(\bS)$ is an $F$-Levi subgroup of \bH.
Since $S\subb$ is a sufficiently large subgroup of both
$S\subb$ and $\bS(E)\subb$\,, in the sense of
\cite{prasad-yu:actions}*{\S 1.3}, Proposition 1.3 of
loc.\ cit.\ gives that
$$
\BB(\bM, E) \cap \BB(\bH, F)
= \BB(\bG, E)^{S\subb} \cap \BB(\bH, F)
= \BB(\bH, F)^{S\subb}
= \BB(\bM \cap \bH, F),
$$
as desired.
\end{proof}

Because we will need it for the next result, we reproduce
here the definitions
(Definition \xref{J-defn:top-F-ss-unip})
of absolute semisimplicity and topological unipotence
from \cite{spice:jordan}
(where they are called absolute $F$-semisimplicity and
topological $F$-unipotence).
The notations $F'$ and $\mc F(F)$ introduced in the definition
are used \emph{only} to define absolute semisimplicity.
The group $\mc F(F)$ is discussed in much more detail in
\S\xref{J-subsec:coefficient} of loc.\ cit.

\begin{dn}
\label{defn:top-F-ss}
Let $F'$ be a complete subfield of $F$ such that $F/F'$ is
unramified.
If $\chr \ff = 0$, then choose a
$\Gal({F'}\unram/F')$-stable subfield of
${F'}\unram$ which is mapped by the natural projection
isomorphically to $\ol\ff$,
and write $\mc F(F)$ for the multiplicative group of this
subfield.
If $\chr \ff = p > 0$, then put
$\mc F(F) =
\bigcup_{\text{$L/F'$ finite unramified}}
	\bigcap_{n = 0}^\infty (L\cross)^{p^n}$.
\indexme{absolutely semisimple}
In either case, an element $\gamma \in G$ is
\emph{absolutely semisimple} if and only if it is
semisimple, and its character values (in the sense of
Definition \ref{defn:root-value}) lie in $\mc F(F)$.
\end{dn}

\begin{prop}
\label{prop:abs-ss-good-descent}
Suppose $g\in G$ is absolutely semisimple
(in the sense of Definition \ref{defn:top-F-ss}).
Then 
$$
\BB(C_\bG(g), F) = \BB(\bG,F)^g.
$$
\end{prop}

This result also appears as Proposition \xref{J-prop:x-depth}
in \cite{spice:jordan}, with a somewhat different proof.
Note that, when $F$ is an algebraic extension of a locally compact
field, $g$ has finite order coprime to $\chr\ff$,
so the result is just Theorem 1.9 of \cite{prasad-yu:actions}.

\begin{proof}
By Proposition 5.1.1 of \cite{rousseau:thesis}, it suffices
to prove this result over some tame extension of $F$.
By Corollary \xref{J-cor:abs-F-ss-tame} of \cite{spice:jordan},
$g$ is tame.
Thus we may, and hence do, assume that $g$ belongs to some
$F$-split torus.
Let \mc F be the (cyclic) group generated by $g$, and
\bS its Zariski closure in \bG.

We prove the result first in case no character value of $g$
(in the sense of Definition \ref{defn:root-value})
projects to a non-trivial root of unity in \ff.
By Lemma \ref{lem:no-roots-unity-torus}, in this case,
\bS is a torus.
By Proposition 1.3 of \cite{prasad-yu:actions}, the result
will follow in this case once we know that \mc F is a
sufficiently large subgroup of $S\subb$ (in the sense of \S 1.3
of loc.\ cit.).

Trivially, $C_\bG(\mc F) = C_\bG(\bS)$.  Now suppose that \bT
is a maximal $F$-split torus containing \bS, and $x$ a point
of $\AA(\bT, F)$.
Denote by $\ms T_x$ the \ff-split torus in $\ms G_x\conn$
corresponding to \bT (so that $\ms T_x(\ff)$ is the image in
$\ms G_x\conn(\ff)$ of $T\subb$), and by \ms g the image of
$g$ in $\ms G_x\conn$.
Put $\bY = \set{\chi \in \bX^*(\bT)}{\chi(g) = 1}$
and $\ms Y = \set{\chi \in \bX^*(\ms T_x)}{\chi(\ms g) = 1}$,
and let $\ms S_x$ be the Zariski closure of the group
generated by \ms g.
By Proposition 8.2(c) of \cite{borel:linear},
$\bS = \bigcap_{\chi \in \bY} \ker \chi$
and
$\ms S_x = \bigcap_{\chi \in \ms Y} \ker \chi$.
Since \bT is $F$-split, there is an identification $i$ of the
cocharacter lattices of \bT and $\ms T_x$ such that the square
$$
\begin{CD}
\bX_*(\bT) \times F\cross_0 @>\text{eval}>> T\subb \\
@V{i \times \text{proj}}VV @V\text{proj}VV \\
\bX_*(\ms T_x) \times \ff\cross @>\text{eval}>> \ms T_x(\ff)
\end{CD}
$$
commutes.
By duality, we deduce an identification
$i^* : \bX^*(\bT) \cong \bX^*(\ms T_x)$ which commutes with
evaluation, in the natural sense.
In particular, no character value of \ms g is a non-trivial
root of unity.
Since $\mc F(F) \cap F\cross_{0+} = \sset 1$
(in the notation of Definition \ref{defn:top-F-ss}),
we have that $i^*(\bY) = \ms Y$.
By another application of
Lemma \ref{lem:no-roots-unity-torus}, $\ms S_x$
is a torus (necessarily \ff-split),
hence is generated by the images of those cocharacters
$\lambda \in \bX_*(\ms T_x)$ such that
$\lambda \circ \chi = 1$ for $\chi \in \ms Y$.
Such a cocharacter may be lifted (via $i$) to a cocharacter
(again called $\lambda$) of \bT such
that $\lambda \circ \chi = 1$ for $\chi \in \bY$;
i.e., to a cocharacter of \bS.
Thus $\ms S_x(\ff)$ is the image of $S\subb$\,.
In particular, the ``centralizer modulo \mf p'' of $S\subb$
is $C_{\ms G_x\conn}(\ms S_x) = C_{\ms G_x\conn}(\ms g)$,
which is equal to the ``centralizer modulo \mf p'' of \mc F.
This is precisely the definition of a sufficiently large
subgroup.

Now it remains to handle the general case.  Pick any pair
$(\bT, x)$ as above, and preserve the above notation.
By Lemma \ref{lem:power-no-roots-unity},
there exists an integer $n$, not divisible by $\chr \ff$, such that
$\ms g^n$ has no non-trivial roots of unity as character
values.  Since the character values of $\ms g^n$ are the
projections to $\ff\cross$ of those of $g^n$, we have that any
root of unity arising as a character value of $g^n$ must
belong to $\mc F(F) \cap F\cross_{0+} = \sset 1$.
Thus, by what we have already shown,
$\BB(C_\bG(g^n), F) = \BB(\bG, F)^{g^n}$.
By Theorem 1.9 of \cite{prasad-yu:actions},
$\BB(C_\bG(g), F) = \BB(C_\bG(g^n), F)^g$.
We are finished.
\end{proof}

\section{Groups associated to concave functions}
\label{sec:concave}

\subsection{Basic definitions}

\begin{defn}
\label{defn:tame-reductive-sqnc}
\indexme{tame reductive sequence}
Call a sequence $\vbG = (\bG^0\subseteq \cdots \subseteq \bG^d)$
of connected reductive $F$-groups a \emph{tame reductive $F$-sequence}
if all $\bG^i$ have the same absolute and $F\tame$-ranks.
We write
$\Lie(\vbG) = (\Lie(\bG^0), \dotsc, \Lie(\bG^d))$.
A \emph{splitting field} for \vbG is any field over which
all $\bG^i$ for $0 \le i \le d$ are split.
\end{defn}

\begin{rk}
By Proposition \ref{prop:compatibly-filtered-tame-rank},
if $\vbG = (\bG^0 \subseteq \dotso \subseteq \bG^d)$
is a tame reductive $F$-sequence, then,
for all $0 \leq i \leq j \leq d$,
$\bG^i$ is a compatibly filtered $F$-subgroup of $\bG^j$.
We will always choose our embeddings of buildings in such a way
that, for all discretely valued extensions $E/F$
and all $0 \le i \le j \le k \le d$, the composition
$\BB(\bG^i, E) \to \BB(\bG^j, E)
\to \BB(\bG^k, E)$ is the same as the embedding
$\BB(\bG^i, E) \to \BB(\bG^k, E)$.
\end{rk}

\begin{defn}
\label{defn:tame-Levi-sequence}
\indexme{tame Levi sequence}
Call a sequence $\vbG = (\bG^0\subseteq \cdots \subseteq \bG^d)$
of connected reductive $F$-groups a \emph{tame Levi $F$-sequence}
if there is a tame extension $L/F$ such that
each $\bG^i$ is an $L$-Levi subgroup
of $\bG^d$.
\end{defn}

\begin{rem}
Note that Definition \ref{defn:tame-Levi-sequence}
is usually more general than the analogous definition
in \cite{yu:supercuspidal}.  If $\bG^d$ is tame, then these
definitions coincide.

It is clear from the definitions that every tame Levi
$F$-sequence is a tame reductive $F$-sequence.
Moreover, any subsequence
of a tame reductive (respectively, tame Levi) $F$-sequence is also a 
tame reductive (respectively, tame Levi) $F$-sequence.
\end{rem}

\begin{dn}
\label{defn:sufficiently-large}
Suppose that \bT is a maximal torus in \bG and
$x \in \AA(\bT, F)$.
Say that $F$ is
\indexme{sufficiently large field}%
\emph{sufficiently large} (for \bT and $x$) if
\bT is $F$-split,
$F$ is strictly Henselian, and,
whenever $\psi \in \PPsiGTF\wtilde\bG\bT F$ satisfies
$\psi(x) \in \Q\dotm\ord(F\cross)$
and
$\lsub F U_\psi \ne \lsub F U_{\psi+}$\,,
then actually
$\psi(x) \in \ord(F\cross)$.
\end{dn}

\begin{rk}
\label{rem:how-large-is-large}
Suppose that $F$ is sufficiently large and $E/F$ is a
discretely valued algebraic extension.
If $\psi$ is an affine $E$-root satisfying
$\psi(x) \in \Q\dotm\ord(E\cross) = \Q\dotm\ord(F\cross)$
and
$\lsub E U_\psi \ne \lsub E U_{\psi+}$\,,
choose an affine $F$-root $\psi_0$ satisfying
$\dpsi_0 = \dpsi$ and
$\lsub F U_{\psi_0} \ne \lsub F U_{\psi_0+}$\,.
An easy calculation, using
the definition of the affine root groups
$\lsub F U_\psi$ and $\lsub F U_{\psi_0}$
(as in \cite{tits:corvallis}*{\S 1.4}, or the proof of
Proposition \ref{prop:compatibly-filtered-tame-rank}),
shows that
the value of the constant function
$\psi - \psi_0$ lies in
$\ord(E\cross) \subseteq \Q\dotm\ord(F\cross)$,
so $\psi_0(x) \in \Q\dotm\ord(F\cross)$.
Since $F$ is sufficiently large,
$\psi_0(x) \in \ord(F\cross)$,
so $\psi(x) \in \ord(E\cross)$.
That is, $E$ is also sufficiently large.

If $\psi(x) \in \ord(F\cross)$, then,
since $\psi_0(x) \in \ord(F\cross)$, the value of the
constant function $\psi - \psi_0$ lies in $\ord(F\cross)$.
Since $\lsub F U_{\psi_0} \ne \lsub F U_{\psi_0+}$\,,
another computation as above shows
that $\lsub F U_\psi \ne \lsub F U_{\psi+}$\,.
\end{rk}

The following definition appears in
\cite{bruhat-tits:reductive-groups-1}*{\S 6.4.3}.

\begin{dn}
\label{defn:concave}
For any root system $\Phi$, a function
$f : \Phi \cup \sset 0 \to \tR$
is called \emph{concave} if, for any
finite non-empty sequence $(\alpha_i)_i$ in $\Phi \cup \sset 0$
with $\sum_i \alpha_i \in \Phi \cup \sset 0$,
we have
$$
f\Bigl(\sum_i \alpha_i\Bigr) \le \sum_i f(\alpha_i).
$$
\end{dn}

Note that, if $f$ is concave, then $f(0) \ge 0$.

\begin{defn}
\label{defn:admissible}
A sequence $\vec r = (r_0,\ldots, r_d)$
of elements of $\tR_{\ge 0}$
is \emph{admissible}
\indexme{admissible sequence (of numbers)}
if $2r_j \ge r_i$ for all $0\leq i \leq j \leq d$.
(Note that this is more general than the definition
given in \cite{yu:supercuspidal}.)
If $\vec r$ is admissible and
$\vbG = (\bG^0, \dotsc, \bG^d)$ is a tame reductive $F$-sequence,
then denote by $f_{\vbG, \vec r}$ the function on
$\wtilde\Phi(\bG, \bT)$
which is equal to $r_i$ on
$\Phi(\bG^i, \bT)
\smallsetminus \Phi(\bG^{i - 1}, \bT)$
for $0 < i \le d$, and to $r_0$ on
$\wtilde\Phi(\bG^0, \bT)$.
(Here, \bT is any maximal $F$-torus in $\bG^0$.  We will
see in Lemma \ref{lem:torus-what-torus} that the choice of
\bT is immaterial.)
\end{defn}

The next result is the analogue of Lemma 1.2 of
\cite{yu:supercuspidal}.

\begin{lm}
If $\vec r$ is admissible and \vbG is a tame reductive
$F$-sequence, then the
function $f_{\vbG, \vec r}$ is concave.
\end{lm}

\begin{proof}
Write $\vec r = (r_0, \ldots, r_d)$ and $\vbG = (\bG^0, \dotsc, \bG^d)$.
Let \bT be a maximal torus in $\bG^0$.
By Lemma 1.1 of \cite{yu:supercuspidal},
it will be enough to show that
$f_{\vbG, \vec r}(\alpha+\beta)
\leq f_{\vbG, \vec r}(\alpha) + f_{\vbG, \vec r}(\beta)$
whenever
$\alpha, \beta, \alpha + \beta \in \wtilde\Phi(\bG, \bT)$.
Fix $\alpha$ and $\beta$ as above, and
let $0 \le i \le d$ and $0 \le j \le d$
be the smallest indices such that
$\alpha \in \wtilde\Phi(\bG^i, \bT)$ and
$\beta \in \wtilde\Phi(\bG^j, \bT)$.
If $i\neq j$, then, assuming,
as we may without loss of generality,
that $i<j$, we have
$f_{\vbG, \vec r}(\alpha + \beta) = 
r_j
\leq r_i + r_j = f_{\vbG, \vec r}(\alpha) + f_{\vbG, \vec r}(\beta)$.
If $i=j$,
then let $0 \le k \le d$
be the smallest index such that
$\alpha +\beta \in \wtilde\Phi(\bG^k, \bT)$.
Since $k \leq i$,
Definition \ref{defn:admissible} implies that
$f_{\vbG, \vec r}(\alpha + \beta) = r_k
\le 2r_i = f_{\vbG, \vec r}(\alpha) + f_{\vbG, \vec r}(\beta)$.
\end{proof}

\begin{dn}
\label{defn:f-positive}
Suppose that $f : \wtilde\Phi(\bG, \bT) \to \tR$ is any
function and $x \in \AA(\bT, \undefined)$ (where \bT is a maximal torus in \bG).
Say that $x' \in \AA(\bT, \undefined)$ is an
\indexme{positive@$(x, f)$-positive!point}%
\indexme{positive@$(x, f)$-positive!point!strictly}%
\emph{$(x, f)$-positive point} (respectively,
\emph{strictly $(x, f)$-positive point})
if $\psi(x') \ge 0$
(respectively, $\psi(x') > 0$)
whenever $\psi$ is an affine root and $\psi(x) = f(\dpsi)$.
If $y \in \AA(\bT, \undefined)$ and
$h : \wtilde\Phi(\bG, \bT) \to \tR_{\ge 0}$ is a concave
function, say that $(y, h)$ is an
\emph{$(x, f)$-positive pair} if
\indexme{positive@$(x,f)$-positive!pair}%
$h(\dpsi) + \psi(x) = f(\dpsi) + \psi(y)$ for every affine
root $\psi$.
\end{dn}

Notice that, if $f_1 \le f_2$, then an
$(x, f_1)$-positive point is also an $(x, f_2)$-positive point.
 
For the rest of this subsection, fix a function
$f : \wtilde\Phi(\bG, \bT) \to \tR$
and a point $x \in \AA(\bT, \undefined)$
(where \bT is a maximal $F$-torus in \bG).

\begin{lm}
\label{lem:pos-pair-has-pos-point}
If $(y, h)$ is an $(x, f)$-positive pair, then
$y$ is an $(x, f)$-positive point.
\end{lm}

\begin{proof}
Suppose that $\psi$ is an affine root with
$\psi(x) = f(\dpsi)$.
If $f(\dpsi) = \infty$, then $\psi(x) = \infty$,
so $\psi$ is the constant function $\infty$.
In particular, $\psi(y) = \infty$.
Otherwise, let $r \in \R$ be such that
$\psi(x) = f(\dpsi) \in \sset{r, r{+}}$.
Then
$(h(\dpsi) + r){+} = (r + \psi(y)){+}$.
	%%In a sense, we're just ``subtracting $r$ from both
	%%sides'' below; but that operation needs to be justified.
Let $s \in \R_{\ge 0} \cup \sset\infty$ be such that
$h(\dpsi) \in \sset{s, s{+}}$.
Then
$(s + r){+} = (r + \psi(y)){+}$,
so $s + r \le r + \psi(y)$, so $s \le \psi(y)$.
Since $s \ge 0$, we are done.
\end{proof}

\begin{lm}
\label{lem:invariant-pos-pair}
Put $\Gamma_{(x, f)} := \stab_{\Gal(F\sep/F)} (x, f)$.
If $f$ is concave, then
there exists a $\Gamma_{(x, f)}$-fixed
$(x, f)$-positive pair $(y, h)$.
If $f$ takes values in $\R \cup \sset\infty$ and $f(0) > 0$,
then $y$ may be taken to be strictly $(x, f)$-positive.
If $\Gamma_{(x, f)} = \Gal(F\sep/F)$, then
$y \in \BB(\bT, F)$.
\end{lm}

Recall that $\bT$ is a compatibly filtered subgroup of $\bG$,
and that we may thus identify $\BB(\bT,F)$ with a particular
subset of $\BB(\bG, F)$.
The final assertion of the lemma is that $y$ belongs to this particular
subset.

\begin{proof}
By Proposition 6.4.6 of
\cite{bruhat-tits:reductive-groups-1}, there is a linear
form $\lambda$ on $\bX^*(\bT) \otimes_\Z \R$ such that
$(f + \lambda)(\alpha) \ge 0$
(or $(f + \lambda)(\alpha) > 0$,
if $f$ takes values in $\R \cup \sset\infty$ and $f(0) > 0$)
for all $\alpha \in \Phi(\bG, \bT)$.
Certainly, $(f + \lambda)(0) = f(0) \ge 0$
(respectively, $(f + \lambda)(0) = f(0) > 0$)
as well.
Put $\tilde h = f + \lambda$.
We may also regard $\lambda$ as an element of
$\bX_*(\bT) \otimes_\Z \R$.
Then put $\tilde y = x + \lambda$.
Clearly, $(\tilde y, \tilde h)$ is an $(x, f)$-positive pair;
and $\tilde y$ is strictly $(x, f)$-positive if $f$ takes
values in $\R \cup \sset\infty$ and $f(0) > 0$.

Since $x$ and $f$ are $\Gamma_{(x, f)}$-fixed,
we have that
$$
(\tilde h \circ \sigma)(\dpsi) + \psi(x)
= (\tilde h \circ \sigma)(\dpsi) + \psi(\sigma x)
%% Some intermediate steps:
% = \tilde h((\psi \circ \sigma)\spdot) + (\psi \circ \sigma)(x)
% = f((\psi \circ \sigma)\spdot) + (\psi \circ \sigma)(\tilde y)
= (f \circ \sigma)(\dpsi) + \psi(\sigma\tilde y)
= f(\dpsi) + \psi(\sigma\tilde y)
$$
for all affine roots $\psi$ and all
$\sigma \in \Gamma_{(x, f)}$.
Thus, if we let $y$ be the center of mass of the
$\Gamma_{(x, f)}$-orbit of $\tilde y$, and $h$ be the average
of the compositions $\tilde h \circ \sigma$ as
$\sigma$ ranges over $\Gamma_{(f, x)}$,
then we have that
$h(\dpsi) + \psi(x) = f(\dpsi) + \psi(y)$
for all affine roots $\psi$.
Since $h$ evidently takes values in $\tR_{\ge 0}$, we have
that $(y, h)$ is a $\Gamma_{(x, f)}$-invariant $(x, f)$-positive
pair.
Further, if $\tilde y$ is strictly $(x, f)$-positive, then
so is $y$.

Suppose that $\Gamma_{(x, f)} = \Gal(F\sep/F)$.
Let $L/F$ be a discretely valued tame extension such that
the maximal $F\tame$-split subtorus $\bS^\sharp$ of \bT is
$L$-split.
By Proposition \ref{prop:compatibly-filtered-tame-rank},
\bT is a compatibly filtered $L$-subgroup of $\bG$
(in fact, it is an $L$-Levi subgroup);
that is, there is a system of embeddings
$(i_E : \BB(\bT, E) \to \BB(\bG, E))_{E/F}$
as in Definition \ref{defn:compatibly-filtered}.
Let $E$ be the splitting field of \bT over $L$.
Then the map $i_E : \BB(\bT, E) \to \BB(\bG, E)$
is $\Gal(E/L)$-equivariant and restricts to the map
$i_L : \BB(\bT, L) \to \BB(\bG, L)
	\subseteq \BB(\bG, E)$.
In particular, the image in $\BB(\bG, E)$ of
$\BB(\bT, L)$ is contained in that of
$\BB(\bT, E)^{\Gal(E/L)} = \AA(\bT, E)^{\Gal(E/L)}$.
Of course, $\BB(\bT, L) = \AA(\bS^\sharp, L)$ is an affine space
under $\bX_*(\bS^\sharp) \otimes_\Z \R$.
By \cite{tits:corvallis}*{\S 1.10},
% Tits makes a blanket assumption that the field is complete.
% However, that plays no role in this result.
so is $\AA(\bT, E)^{\Gal(E/L)}$.
Thus $\AA(\bT, E)^{\Gal(E/L)} = \BB(\bT, L)$.
Since $y$ is $\Gamma_{(x, f)}$-fixed, hence $\Gal(E/L)$-fixed, we
have that $y \in \AA(\bS^\sharp, L) = \BB(\bT, L)$.
Since $y$ is, further, $\Gal(L/F)$-fixed, we have that
$y \in \BB(\bT, L)^{\Gal(L/F)} = \BB(\bT, F)$.
\end{proof}

\subsection{Groups associated to concave functions}
\label{sec:concave_gen}
Fix, for the remainder of this section,
\begin{itemize}
\item a tame reductive
$F$-sequence $(\bT, \bG^0, \cdots, \bG^d = \bG)$
with \bT a torus,
\item
a point $x$ lying in the image of
$\BB(\bT, F)$ in $\BB(\bG, F)$,
\item
a $\Gal(F\sep/F)$-invariant concave function
$f : \wtilde\Phi(\bG, \bT) \to \tR$,
and
\item
a $\Gal(F\sep/F)$-invariant $(x, f)$-positive pair $(y, h)$
with $y \in \BB(\bT, F)$.
\end{itemize}
By Lemma \ref{lem:invariant-pos-pair},
a pair $(y, h)$ as above exists.
Put $\vbG = (\bG^0, \dotsc, \bG^d)$.

\begin{dn}
\label{defn:vGvr-split}
Suppose that $F$ is sufficiently large
for \bT and $y$.
Put
$$
\lsub\bT G_{x, f}
:= \langle U_\psi
	:
	\psi \in \PPsiGTF\wtilde\bG\bT F\text{ and }
	\: \psi(x) \ge f(\dpsi)\rangle.
$$
If $f = f_{\vbG, \vec r}$
for an admissible sequence $\vec r = (r_0, \ldots, r_d)$,
put
$$
\vG_{x, \vec r} := \lsub\bT G_{x, f}\,.
$$
\end{dn}

Note that the definition still makes sense if we only assume
that \bT is $F$-split (rather than that $F$ is sufficiently
large), but then it might not agree with
Definition \ref{defn:vGvr} below
when $f(0) = 0$.
Note further that, as in \cite{yu:supercuspidal}*{\S 1},
$\vG_{x, \vec r}$ is just an open
subgroup of $G$, not a sequence of such subgroups.

Fix, for the remainder of this section,
a discretely valued Galois extension $E/F$ such that
$E$ is sufficiently large for \bT and $y$.

\begin{dn}
\label{defn:vGvr}
\indexmem{\lsub\bT G_{x, f}}%
\indexmem{\vG_{x, \vec r}}%
Put
$\lsub\bT G_{x, f} := \lsub\bT\bG(E)_{x, f} \cap G_0$\,.
For an admissible sequence $\vec r$, put
$\vG_{x, \vec r} := \vbG(E)_{x, \vec r} \cap G_0$\,.
We will sometimes write
$(G^0, \ldots, G^d)_{x, \vec r}$ instead of
$\vG_{x, \vec r}$\,.
If \bH is a compatibly filtered $F$-subgroup of \bG
containing \bT, then we will use $\lsub\bT H_{x, f}$ as
shorthand for $\lsub\bT H_{x, g}$\,, where
$g = f\bigr|_{\wtilde\Phi(\bH, \bT)}$.
\end{dn}

Note that $\lsub\bT G_{x, f} = \lsub\bT G_{y, h}$\,.

We will show (see Lemma \ref{lem:big-split})
that $\lsub\bT G_{x,f}$ does not depend
on the choice of sufficiently large extension $E$;
and (see Lemma \ref{lem:torus-what-torus})
that, if $\vec r$ is admissible, then
$\vG_{x, \vec r}$ is independent of the choice of \bT
(see Lemma \ref{lem:torus-what-torus}).
For example, under Hypothesis \eqref{hyp:reduced},
if $\vec r = (r, r, \ldots, r)$, then
$\vG_{x, \vec r} = G_{x, r}$
(see Remark \ref{rem:vGvr-facts}).

\begin{rk}
\label{rem:infinite-vec}
We have defined $\vG_{x, \vec r}$ when \vbG is indexed
by a finite set.  For later purposes, it will be convenient
to handle the case where
$\vbG = (\bG^i)_{i \in I}$
and
$\vec r = (r_i)_{i \in I}$
are indexed by any totally ordered set $I$ (such as an
interval).
Suppose that $I \ne \emptyset$.
Since $\bG^i \mapsto \Phi(\bG^i, \bT)$
is an injection from the set
of groups in \vbG to the set of subsets
of $\Phi(\bG, \bT)$, there are actually only finitely many
distinct groups appearing in \vbG.
Denote these finitely many groups by
$\vec\bH = (\bH^0 \subseteq \cdots \subseteq \bH^d)$,
and put $\vec s = (s_0, \ldots, s_d)$, where
$s_j = \inf\set{r_i}{\bG^i = \bH^j}$ for
$0 \le j \le d$.
(Here, the infimum is taken in \tR, not \R.  Thus, for
example, the infima of the intervals $(0, 1)$ and
$(0{+}, 1)$ are both $0{+}$, not $0$.)
Then we define $\vG_{x, \vec r} := \vec H_{x, \vec s}$\,.
If $I = \emptyset$, put $\vG_{x, \vec r} := \sset 1$.
\end{rk}

\begin{lm}
\label{lem:sloppy-vGvr}
If $x' \in \BB(\bT, F)$ is an $(x, f)$-positive point, then,
for any subextension $L/F$ of $E/F$, we have
$\lsub\bT\bG(L)_{x, f}
= \lsub\bT\bG(E)_{x, f} \cap \bG(L)_{x', 0}$\,.
\end{lm}

\begin{proof}
Fix $x'$ and $L/F$ as in the statement of the lemma.
Then clearly
$\lsub\bT\bG(E)_{x, f}$ fixes $x'$,
and $x' \in \BB(\bT, L)$.
By Lemma \ref{lem:stab-deep}, we have
$\lsub\bT\bG(L)_{x, f} =
\lsub\bT\bG(E)_{x, f} \cap \bG(L)_0
\subseteq \lsub\bT\bG(E)_{x, f} \cap \bG(L)_{x', 0}$\,.
The reverse containment being obvious, we have equality.
\end{proof}

\begin{lm}
\label{lem:master-comm}
Let
$f_j : \wtilde\Phi(\bG, \bT) \to \tR$
be a $\Gal(E/F)$-invariant, concave function
for $j = 1, 2$.
Let $f_1 \vee f_2$ be as in Definition
\ref{defn:concave-vee}.
Suppose that
\begin{itemize}
\item
$(f_1 \vee f_2)(0) \ne -\infty$
and
\item
there is some point $x' \in \BB(\bT, F)$
which is both $(x, f_1)$- and $(x, f_2)$-positive.
\end{itemize}
(For example, this occurs when $f_1$ and $f_2$ are both
non-negative, or $f_1 \le f_2$.)
Then
$[\lsub\bT G_{x, f_1}, \lsub\bT G_{x, f_2}]$ is contained in
$\lsub\bT G_{x, f_1 \vee f_2}$\,.
\end{lm}

\begin{proof}
Since $(f_1 \vee f_2)(0) \ne -\infty$,
Proposition 6.4.44 of
\cite{bruhat-tits:reductive-groups-1} 
shows that $f_1 \vee f_2$ is concave and
$$
[\lsub\bT\bG(E)_{x, f_1}, \lsub\bT\bG(E)_{x, f_2}]
\subseteq \lsub\bT\bG(E)_{x, f_1 \vee f_2}\,.
$$
(The cited proposition
depends on condition (Pr) of
\cite{bruhat-tits:reductive-groups-1}, which, by Proposition
6.4.39 of loc.\ cit., is equivalent
to the existence of a ``prolongation of root datum'', in the
sense of \S 6.4.38 of loc.\ cit.
By Lemma 6.1 of \cite{yu:models}, such a prolongation
exists.)
If further there exists a point $x'$ as in the statement of
the lemma, then, by Lemma \ref{lem:sloppy-vGvr},
$$
[\lsub\bT G_{x, f_1}, \lsub\bT G_{x, f_2}]
\subseteq \lsub\bT\bG(E)_{x, f_1 \vee f_2} \cap G_{x', 0}
\subseteq \lsub\bT G_{x, f_1 \vee f_2}\,.\qedhere
$$
\end{proof}

\begin{cor}
\label{cor:master-comm}
Suppose that
$\vec{s}^{(j)} = (s^{(j)}_0, \ldots, s^{(j)}_d)$
are admissible sequences for $j = 1, 2$.
Put $t^{(j)} = \min\sset{s^{(j)}_0, \ldots, s^{(j)}_d}$
for $j = 1, 2$,
$$
r_0 = \min\set{\smash{s^{(1)}_k + s^{(2)}_k}}{0 \le k \le d},
$$
and
$$
r_i = \min\set{\smash{s^{(1)}_i + t^{(2)},
	t^{(1)} + s^{(2)}_i,
	s^{(1)}_k + s^{(2)}_k}}{i < k \le d}
$$
for all $0 < i \le d$.
Then
$[\vG_{x, \vec s^{(1)}}, \vG_{x, \vec s^{(2)}}]
	\subseteq \vG_{x, \vec r}$\,.
\end{cor}

\begin{proof}
Note that $\vec r$ is admissible.
Put $f_j = f_{\vbG, \vec s^{(j)}}$ for $j = 1, 2$,
and $f = f_{\vbG, \vec r}$.
(See Definition \ref{defn:vGvr-split}.)
For notational convenience,
put
$\Phi^0 = \wtilde\Phi(\bG^0, \bT)$,
and
$\Phi^i
= \wtilde\Phi(\bG^i, \bT)
\smallsetminus \wtilde\Phi(\bG^{i - 1}, \bT)$
for $0 < i \le d$.
Note that $0 \in \Phi^0$.
Fix $0 \le i \le d$ and $\alpha \in \Phi^i$.
Write, in any fashion whatsoever,
$\alpha = \sum a_m + \sum b_n$,
where $(a_m)_m$ and $(b_n)_n$ are finite non-empty sequences in
$\wtilde\Phi(\bG, \bT)$.

First suppose that $i = 0$.
Let $k$ be the greatest index for which some
$a_m$ or $b_n$ is in $\Phi^k$.
If there are $a_m, b_n \in \Phi^k$
(in particular, if $k = 0$), then
$$
\sum_m f_1(a_m) + \sum_n f_2(b_n) \ge s^{(1)}_k + s^{(2)}_k
	\ge f(\alpha).
$$
Otherwise, $k > 0$ and there are distinct indices $m \ne m'$,
or $n \ne n'$, such that $a_m, a_{m'} \in \Phi^k$,
or $b_n, b_{n'} \in \Phi^k$.
In the former case, let $0 \le j < k$ be such that
some $b_n$ is in $\Phi^j$.
Then
$$
\sum_m f_1(a_m) + \sum_n f_2(b_n) \ge 2s^{(1)}_k + s^{(2)}_j
	\ge s^{(1)}_j + s^{(2)}_j \ge f(\alpha).
$$
The latter case is handled similarly.
Thus $(f_1 \vee f_2)(\alpha) \ge f(\alpha)$.
(Here, $\vee$ is as in Definition \ref{defn:concave-vee}.)

Now suppose that $i > 0$.
Note that some $a_m$ or $b_n$
must lie in $\Phi^k$ for some $k \ge i$.
If some $a_m \in \Phi^i$,
then
$$
\sum_m f_1(a_m) + \sum_n f_2(b_n) \ge s^{(1)}_i + t^{(2)}
	\ge f(\alpha).
$$
If some $b_n \in \Phi^i$,
then similarly
$$
\sum_m f_1(a_m) + \sum_n f_2(b_n) 
	\ge f(\alpha).
$$
If, for some $k > i$,
there are $a_m, a_{m'} \in \Phi^k$ with $m \ne m'$, then
$$
\sum_m f_1(a_m) + \sum_n f_2(b_n)
	\ge 2s^{(1)}_k + t^{(2)}
	\ge s^{(1)}_i + t^{(2)}
	\ge f(\alpha).
$$
If, for some $k > i$,
there are $b_n, b_{n'} \in \Phi^k$ with $n \ne n'$, then
similarly
$$
\sum_m f_1(a_m) + \sum_n f_2(b_n)
	\ge f(\alpha).
$$
If there are $a_m, b_n \in \Phi^k$ for some $k > i$,
then
$$
\sum_m f_1(a_m) + \sum_n f_2(b_n) \ge s^{(1)}_k + s^{(2)}_k
	\ge f(\alpha).
$$
Thus $(f_1 \vee f_2)(\alpha) \ge f(\alpha)$.

Since $\alpha \in \wtilde\Phi(\bG, \bT)$ was arbitrary, the
result follows from Lemma \ref{lem:master-comm}.
\end{proof}

Recall that the construction of $\vG_{x, \vec r}$
depended on a torus \bT, which we suppressed from the
notation.  We will temporarily indicate the dependence on
\bT by writing $\lsub\bT\vG_{x, \vec r}$\,.

\begin{lm}
\label{lem:stab-norm}
For any admissible sequence $\vec r$
and $g \in G^0 \cap \stab_G(\ox)$,
$\lsup g(\lsub\bT\vG_{x, \vec r})
= \lsub{\lsup g\bT}\vG_{x, \vec r}$\,.
\end{lm}

\begin{proof}
The action of $g$ induces bijections
from
$\AA(\bT, E)$ to $\AA(\lsup g\bT, E)$
and
\PPsiGTF\wtilde\bG\bT E to
\PPsiGTF\wtilde\bG{\lsup g\bT}E
such that, if
$y \mapsto y'$ and
$\psi \mapsto \psi'$,
then
$\lsup g U_\psi = U_{\psi'}$ and
$\psi(y) = \psi'(y')$.
Since $\psi$ and $\psi'$ depend only on the images of their
arguments in $\rBB(\bG, F)$, and since $x$ and $g x$ have
the same image there,
we have in particular that $\psi(x) = \psi'(x)$.
Similarly, the action of $g$ induces a bijection from
$\wtilde\Phi(\bG, \bT)$ to
$\wtilde\Phi(\bG, \lsup g\bT)$ such that
$\wtilde\Phi(\bG^0, \bT)$ is carried to
$\wtilde\Phi(\bG^0, \lsup g\bT)$, and
$\wtilde\Phi(\bG^i, \bT)
\smallsetminus \wtilde\Phi(\bG^{i - 1}, \bT)$
is carried to
$\wtilde\Phi(\bG^i, \lsup g\bT)
\smallsetminus \wtilde\Phi(\bG^{i - 1}, \lsup g\bT)$
for each $0 < i \le d$.
Thus
$\lsup g(\lsub\bT\vbG(E)_{x, \vec r})
= \lsub{\lsup g\bT}\vbG(E)_{x, \vec r}$\,.
Since $\lsub\bT\vG_{x, \vec r}
= \lsub\bT\vbG(E)_{x, \vec r} \cap G_0$\,,
and similarly for $\lsub{\lsup g\bT}\vG_{x, \vec r}$\,,
and since conjugation by $g$ preserves $G_0$\,,
we have the desired equality.
\end{proof}

\begin{lemma}
\label{lem:torus-what-torus}
Suppose that
$\vec r$ is admissible, and
$\bT_j$ is an $E$-split $F$-torus such that
$(\bT_j , \bG^0)$ is a
tame reductive sequence
and $x \in \BB(\bT_j, F)$
for $j = 1, 2$.
Then
$\lsub{\bT_1}\vG_{x, \vec r}
=\lsub{\bT_2}\vG_{x, \vec r}$\,.
\end{lemma}

\begin{proof}
By Proposition 4.6.28(iii) of
\cite{bruhat-tits:reductive-groups-2},
there is $g \in \bG^0(E)_{x, 0}$ such that
$\lsup g\bT_1 = \bT_2$.
By Lemma \ref{lem:stab-norm},
$\lsup g(\lsub{\bT_1}\vbG(E)_{x, \vec r})
= \lsub{\bT_2}\vbG(E)_{x, \vec r}$\,.
By Definition \ref{defn:vGvr-split},
$\bG^0(E)_{x,0}
= \lsub{\bT_1}\vbG(E)_{x,(0,\infty,\dotsc, \infty)}$,
so, by
Corollary \ref{cor:master-comm},
$\bG^0(E)_{x, 0}$
normalizes $\lsub{\bT_1}\vbG(E)_{x, \vec r}$\,.  Thus in fact
$\lsub{\bT_1}\vbG(E)_{x, \vec r}
= \lsub{\bT_2}\vbG(E)_{x, \vec r}$\,.
As in the proof of Lemma \ref{lem:stab-norm}, we conclude that
$\lsub{\bT_1}\vG_{x, \vec r}
= \lsub{\bT_2}\vG_{x, \vec r}$\,.
\end{proof}

\begin{cor}
\label{cor:stab-norm}
$G^0 \cap \stab_G(\ox)$ normalizes $\vG_{x, \vec r}$
for any admissible sequence $\vec r$.
\end{cor}

\begin{proof}
This follows from
Lemmata \ref{lem:stab-norm} and \ref{lem:torus-what-torus}.
\end{proof}

For the remainder of this section,
suppose that $L/F$ is a Galois subextension of $E/F$ such that
\begin{itemize}
\item \bG is $L$-quasisplit,
\item \bT contains a maximal $L$-split torus $\bS^\sharp$
(necessarily defined over $F$),
and
\item
$\lsub L\Phi(\bG)$ is reduced.
\end{itemize}
Note that, by
Proposition 15.5.3(iii) of \cite{springer:lag},
$\wtilde\Phi(\bG, \bS^\sharp)$ is the set of
$\Gal(E/L)$-orbits in $\wtilde\Phi(\bG, \bT)$,
so that $f$ may also be regarded as a function on
$\wtilde\Phi(\bG, \bS^\sharp)$.
For $(a, c) \in \bX_*(\bS^\sharp) \times \tR$
with $c < \infty$,
write $a + c$ as a shorthand for the unique affine function
$\varphi$ on $\AA(\bS^\sharp, L)$ with
$\dot\varphi = a$ and $\varphi(x) = c$.
If $c = \infty$, then write $a + c$ for the constant function
with value $\infty$ on $\AA(\bS^\sharp, L)$.

\begin{lm}
\label{lem:gen-iwahori-factorization}
Suppose that
$f(0) > 0$.
Then the multiplication map
\begin{equation*}
\mult_f :
\prod_{ a \in \wtilde\Phi(\bG, \bS^\sharp)
	}
     U_{a + f(a)} 
\to
\lsub\bT\bG(L)_{x, f}
\end{equation*}
(the product taken in any order) is a continuous bijection.
\end{lm}

\begin{proof}
By Proposition 6.4.48 of
\cite{bruhat-tits:reductive-groups-1}
and Lemma \ref{lem:complete-subfield},
\begin{equation*}
\mult_{f, E} :
\prod_{ \alpha \in \wtilde\Phi(\bG, \bT)
	}
		\lsub E U_{\alpha + f(\alpha)}
\to \lsub\bT\bG(E)_{x, f}
\end{equation*}
is a bijection.
(Here, we have ordered the product so that the factors
corresponding to roots with the same restriction to
$\bS^\sharp$ are contiguous.)
Since $\mult_f$ is a restriction of
$\mult_{f, E}$, it suffices to
show that $\mult_f$ is surjective.
Choose
$g \in \lsub\bT\bG(L)_{x, f} \subseteq \lsub\bT\bG(E)_{x, f}$,
and denote by
$(u_\alpha)_{\alpha \in \wtilde\Phi(\bG, \bT)}$ its preimage
under $\mult_{f, E}$.
Then $u_0 \in \bT(E)_{f(0)}$\,;
and, for $a \in \Phi(\bG, \bS^\sharp)$,
Lemmata \ref{lem:filter-descent} and \ref{lem:Zariski} give
$u_a := \prod_{
		\alpha\bigr|_{\bS^\sharp} = a 
		}
	 u_\alpha \in U_{a + f(a)}$\,. 
In particular, $u_a \in \bG(L)_{y, 0}$ for
$a \in \Phi(\bG, \bS^\sharp)$, so
$u_0 \in \bG(L)_{y, 0}$ also.
By Lemmata \ref{lem:torus-field-descent} and
\ref{lem:levi-descent},
$$
u_0 \in \bT(E)_{f(0)} \cap \bG(L)_{y, 0}
	= \bT(E)_{f(0)} \cap \bT(L)_0
	= \bT(L)_{f(0)} = U_{0 + f(0)}\,.
$$
Then
$(u_a)_{a \in \wtilde\Phi(\bG, \bS^\sharp)}$
is in the domain of $\mult$, and clearly
$g = \mult_f((u_a))$.
\end{proof}

\begin{rk}
\label{rem:little-split}
If $f(0) > 0$, then, by Lemma \ref{lem:gen-iwahori-factorization},
$\lsub\bT\bG(L)_{x, f}$ is
generated by the root subgroups $U_\varphi$ with
$\varphi$ an affine $L$-root satisfying
$\varphi(x) = f(\dot\varphi)$.
In the notation of
\cite{bruhat-tits:reductive-groups-1}*{\S 6.4.2},
$\lsub\bT\bG(L)_{x, f} = X\dotm U_{f^\perp}$,
where $X = U_{0 + f(0)} \subseteq \bT(L)\subb$
and $f^\perp$ is the restriction of $f$ to $\Phi(\bG, \bT)$.
Thus, if \bT is $L$-split, then the group
$\lsub\bT\bG(L)_{x, f}$ constructed here is the same as the
one constructed in Definition \ref{defn:vGvr-split}.
(We will prove later that,
if $L$ is sufficiently large for \bT and $y$, then
this is true even if $f(0) = 0$.
See Lemma \ref{lem:big-split}.)
\end{rk}

\begin{cor}
\label{cor:very-pos-point}
If $f(0) > 0$, then there is a point
$x' \in \BB(\bT, F)$ such that, for any subextension $K/F$
of $L/F$ with $L/K$ tame, we have
$\lsub\bT\bG(K)_{x, f}
= \lsub\bT\bG(L)_{x, f} \cap \bG(K)_{x', 0{+}}$\,.
\end{cor}

\begin{proof}
By Lemma \ref{lem:concave-in-R}, for $\delta$ sufficiently
small, $f_\delta$ is concave.
Since $L$ is discretely valued, we have that, for $\delta$
smaller still,
$U_{a + f(a)} = U_{a + f_\delta(a)}$ for
$a \in \wtilde\Phi(\bG, \bS^\sharp)$, hence, by Remark
\ref{rem:little-split}, that
$\lsub\bT\bG(L)_{x, f} = \lsub\bT\bG(L)_{x, f_\delta}$\,.
By Lemma \ref{lem:invariant-pos-pair}, there is a strictly
$(x, f_\delta)$-positive point $x'$.
In particular, for any $a \in \wtilde\Phi(\bG, \bS^\sharp)$,
we have that $(a + f_\delta(a))(x') > 0$,
so $U_{a + f_\delta(a)} \subseteq \bG(L)_{x', 0{+}}$\,.
Therefore, by another application of Remark
\ref{rem:little-split},
$\lsub\bT\bG(L)_{x, f_\delta} \subseteq \bG(L)_{x', 0{+}}$\,.
It is an easy consequence of Definition \ref{defn:vGvr} and
Lemma \ref{lem:domain-field-ascent} that
$\lsub\bT\bG(K)_{x, f}
= \lsub\bT\bG(L)_{x, f} \cap \bG(K)_0$\,.
Since $L/K$ is tame, we have by Lemma
\ref{lem:field-descent} that
\begin{multline*}
\lsub\bT\bG(K)_{x, f}
= \lsub\bT\bG(L)_{x, f} \cap \bG(K)_0 \\
= \lsub\bT\bG(L)_{x, f} \cap \bG(L)_{x', 0+} \cap \bG(K)
= \lsub\bT\bG(L)_{x, f} \cap \bG(K)_{x', 0+}\,.\qedhere
\end{multline*}
\end{proof}

\begin{cor}
\label{cor:gen-iwahori-factorization-bij}
Suppose that
\begin{itemize}
\item
$f_j$ is a $\Gal(E/F)$-invariant, concave function
with $f_j(0) > 0$
for $j = 1, 2$,
\item $f_2 \ge f_1$,
\item $f_1 \vee f_2 \ge f_2$,
and
\item $(f_1 \vee f_2)(\alpha) > f_2(\alpha)$
whenever $f_2(\alpha) < \infty$.
\end{itemize}
(Here, the operator $\vee$ is as in Definition
\ref{defn:concave-vee}.)
Then the composition
$$
\prod_{a \in \wtilde\Phi(\bG, \bS^\sharp)}
	U_{a + f_1(a)} 
\xrightarrow{\mult_{f_1}} \lsub\bT\bG(L)_{x, f_1}
\to \lsub\bT\bG(L)_{x, f_1:f_2}
$$
induces a continuous bijection
$$
\mult_{f_1:f_2} :
\prod_{a \in \wtilde\Phi(\bG, \bS^\sharp)}
	U_{(a + f_1(a)):(a + f_2(a))}
\to \lsub\bT\bG(L)_{x, f_1:f_2}\,.
$$
\end{cor}

\begin{proof}
By Lemma \ref{lem:concave-in-R}, for $\delta$ and
$\varepsilon$ sufficiently small,
the function
$f_{j, \delta}$
defined in that lemma
is concave for $j = 1, 2$,
and
$f_{1, \delta} \vee f_{2, \delta}
\ge f_{2, \delta} + \varepsilon$
(with notation as in the statement of that lemma).
Since $L$ is discretely valued, we have that, for
$\delta$ smaller still,
$U_{a + f_j(a)} = U_{a + f_{j, \delta}(a)}$
for $a \in \wtilde\Phi(\bG, \bS)$,
hence, by Remark \ref{rem:little-split}, that
$\lsub\bT\bG(L)_{x, f_j}
= \lsub\bT\bG(L)_{x, f_{j, \delta}}$\,,
for $j = 1, 2$.
Let $\delta$ and $\varepsilon$ be so small that all of the
above conditions are satisfied.
Then we may, and hence do, replace $f_j$ by $f_{j, \delta}$
for $j = 1, 2$, so that
$f_1 \vee f_2 \ge f_2 + \varepsilon$.

Define concave functions $f_j$ for $j \in \Z_{> 2}$ by
$f_j = f_2 + (j - 2)\varepsilon$.
Then the groups $\lsub\bT\bG(L)_{x, f_j}$ are a basis of
neighborhoods of the identity in $\lsub\bT\bG(L)_{x, f_2}$.

By normality of $\lsub\bT\bG(L)_{x, f_2}$,
the indicated composition is constant on
cosets of $\prod U_{a + f_2(a)}$.
Thus, there is a map $\mult_{f_1:f_2}$ as indicated.
It is clearly continuous and surjective.

To show injectivity, suppose that
$\vec u_1, \vec u_2
\in \prod U_{a + f_1(a)}$
satisfy
$\mult_{f_1:f_2}(\vec u_1) = \mult_{f_1:f_2}(\vec u_2)$.
By Lemma \ref{lem:complete-subfield},
there is a complete subfield $L'$ of $L$ such that
\begin{itemize}
\item \bG is defined over $L'$,
\item $\bS^\sharp$ is split over $L'$,
\item $x \in \BB(\bG, L')$,
and
\item
$\vec u_1, \vec u_2
\in \prod \lsub{L'}U_{a + f_1(a)}$.
\end{itemize}
Upon replacing $L$ by $L'$, we may, and hence do, assume
that $L$ is complete.
Now suppose that $j \in \Z_{\ge 2}$ and
$\vec u_j
\in \vec u_2\prod U_{a + f_2(a)}$
satisfies
$\mult_{f_1:f_j}(\vec u_1) = \mult_{f_1:f_j}(\vec u_j)$.
Put
$$
\vec w
= \mult_{f_j}\inv(
	\mult_{f_1}(\vec u_j)\inv\mult_{f_1}(\vec u_1)
	).
$$
Since
$f_1 \vee f_j
\ge (f_1 \vee f_2) + (j - 2)\varepsilon
\ge f_2 + (j - 1)\varepsilon
= f_{j + 1}$,
Lemma \ref{lem:master-comm} gives
$[\lsub\bT\bG(L)_{x, f_1}, \lsub\bT\bG(L)_{x, f_j}]
\subseteq \lsub\bT\bG(L)_{x, f_{j + 1}}$\,.
In particular,
\begin{multline*}
\mult_{f_1}(\vec u_1) = \mult_{f_1}(\vec u_j)\mult_{f_j}(\vec w)
	= \Bigl(\prod u_{j, a}\Bigr)\Bigl(\prod w_a\Bigr) \\
	\equiv \prod (u_{j, a}w_a)
	= \mult_{f_1}(\vec u_j\vec w)
	\pmod{\lsub\bT\bG(L)_{x, f_{j + 1}}}.
\end{multline*}
Put $\vec u_{j + 1} := \vec u_j\vec w$.
Then $(\vec u_j)_{j \in \Z_{\ge 2}}$ is a Cauchy sequence,
say with limit $\vec v$.
We have
$\vec v
\in \vec u_2\prod U_{a + f_2(a)}$;
and
$\mult_{f_1}(\vec v)
= \lim_{j \to \infty} \mult_{f_1}(\vec u_j)
= \mult_{f_1}(\vec u_1)$,
so $\vec v = \vec u_1$.
\end{proof}

\begin{cor}
\label{cor:gen-iwahori-factorization-iso}
With notation as in Corollary
\ref{cor:gen-iwahori-factorization-bij}, if
$\lsub\bT\bG(L)_{x, f_2}$ contains
$[\lsub\bT\bG(L)_{x, f_1}, \lsub\bT\bG(L)_{x, f_1}]$,
then $\mult_{f_1:f_2}$ is a $\Gal(L/F)$-equivariant isomorphism.
\end{cor}

\begin{proof}
Since
\begin{multline*}
\mult_{f_1}(\vec u)\mult_{f_1}(\vec u^{\,\prime})
	= \bigl(\prod u_a\bigr)\bigl(\prod u'_a\bigr) \\
	\equiv \prod (u_a u'_a)
		= \mult_{f_1}(\vec u\,\vec u^{\,\prime})
	\pmod{[\lsub\bT\bG(L)_{x, f_1}, \lsub\bT\bG(L)_{x, f_1}]}
\end{multline*}
for
$\vec u, \vec u^{\,\prime} \in \prod U_{a + f_1(a)}$,
we have that $\mult_{f_1:f_2}$ is a homomorphism.
We have already shown that it is a bijection.

Now note that $\lsub\bT\bG(L)_{x, f_1:f_2}$
is stable under the $\Gal(L/F)$-action on $\bG(L)$.
The $\Gal(L/F)$-action on the domain of $\mult_{f_1:f_2}$
is deduced from the action on
$\prod U_{a + f_1(a)}$
defined by
$(\sigma\vec u)_a := \sigma u_{\sigma\inv a}$ for
$\vec u \in \prod U_{a + f_1(a)}$,
$a \in \Phi(\bG, \bS^\sharp)$,
and $\sigma \in \Gal(L/F)$.
Equivariance follows from the fact that
\begin{multline*}
\mult_{f_1}(\sigma\vec u) = \prod \sigma u_{\sigma\inv a}
	\equiv \prod \sigma u_a \\
	= \sigma\bigl(\prod u_a\bigr) = \sigma\mult_{f_1}(\vec u)
	\pmod{[\lsub\bT\bG(L)_{x, f_1}, \lsub\bT\bG(L)_{x, f_1}]}.
\qedhere
\end{multline*}
\end{proof}

\begin{rk}
Although we do not need to do so, one can use
Corollary \ref{cor:gen-iwahori-factorization-iso} to show that
$\mult_{f_1}$ is an open map, hence a homeomorphism.
\end{rk}

\begin{lm}
\label{lem:big-split}
The group $\lsub\bT G_{x, f}$ is independent
of the separable extension of $F$ chosen as a
sufficiently large field for \bT and $y$.
\end{lm}

\begin{proof}
To obviate some confusion, we will temporarily use
the notation $\lsupb *\bT\bG(\cdot)_{x, f}$ for the groups
constructed in Definition \ref{defn:vGvr-split},
where the ground field was assumed to be sufficiently large.
It is enough to show that, if $F$ is sufficiently large, then
$\lsupb *\bT G_{x, f} = \lsub\bT G_{x, f}$\,,
i.e., that
$\lsupb *\bT G_{x, f} = \lsupb *\bT\bG(E)_{x, f} \cap G_0$\,.
Upon replacing $x$ by $y$ and $f$ by $h$, we may, and hence
do, assume that $f$ is non-negative and $F$ is sufficiently
large for \bT and $x$.
Since the containment
$\lsupb *\bT\bG(E)_{x, f} \cap G_0
\supseteq \lsupb *\bT G_{x, f}$
is obvious, we consider only the reverse containment.
By Lemma \ref{lem:sloppy-vGvr},
$\lsupb *\bT\bG(E)_{x, f} \cap G_0 \subseteq G_{x, 0}$\,,
so it suffices to show
$\lsupb *\bT\bG(E)_{x, f} \cap G_{x, 0}
\subseteq \lsupb *\bT G_{x, f}$\,.

Put $f_+ = \max\sset{0{+}, f}$.
Clearly, $f \vee f_+ \ge f \vee f \ge f$.
On the other hand, fix $\alpha \in \wtilde\Phi(\bG, \bT)$
and let $(a_m)_m$ and $(b_n)_n$ be any finite
non-empty sequences in $\wtilde\Phi(\bG, \bT)$ such that
$\sum_m a_m + \sum_n b_n = \alpha$.
Let $b_{n_0}$ be any term of $(b_n)_n$.  Then
$$
\sum_m f(a_m) + \sum_n f_+(b_n)
\ge f_+(b_{n_0}) \ge 0{+}.
$$
Thus, $(f \vee f_+)(\alpha) \ge 0{+}$.
Since $\alpha \in \wtilde\Phi(\bG, \bT)$ was arbitrary, we
have $f \vee f_+ \ge \max\sset{0{+}, f} = f_+$.
By Lemma \ref{lem:master-comm},
$\lsub\bT\bG(E)_{x, f_+}$ is normal in
$\lsub\bT\bG(E)_{x, f}$\,.
Notice that $\lsupb *\bT\bG(E)_{x, f:f_+}$ is
generated by the subgroups $\lsub E U_{\psi:\psi+}$\,,
where $\psi$ is an affine $E$-root satisfying
$\lsub E U_\psi \ne \lsub E U_{\psi+}$
and
$\psi(x) = 0 = f(\dpsi)$.
Choose such a $\psi$.
Since $\psi(x) \in \ord(F\cross)$
and $F$ is sufficiently large for \bT and $x$,
Remark \ref{rem:how-large-is-large} gives
$\lsub F U_\psi \ne \lsub F U_{\psi+}$\,.
The map
$F_0 \cong \lsub F U_\psi
\to \lsub E U_\psi \cong E_0$\,,
deduced from the isomorphisms
$F \to U_\dpsi$ and $E \to \bU_\dpsi(E)$
of \S\ref{sec:unip-exp},
restricts to a map
$F_{0+} \cong \lsub F U_{\psi+}
\to \lsub E U_{\psi+} \cong E_{0+}$\,.
Since $\ff_E/\ff$ is algebraic, and \ff is separably (hence
algebraically, since it is perfect) closed,
the induced map
$\ff = F\unram_{0:0+} \to E_{0:0+} = \ff_E$
is an isomorphism.
Therefore,
$\lsub F U_{\psi:\psi+} \to \lsub E U_{\psi:\psi+}$
is also an isomorphism;
in particular,
$\lsub E U_\psi
\subseteq \lsub F U_\psi\dotm\lsub E U_{\psi{+}}$\,.

Thus
$\lsupb *\bT\bG(E)_{x, f}
\subseteq \lsupb *\bT G_{x, f}
		\cdot\lsupb *\bT\bG(E)_{x, f_+}$\,,
so
$\lsupb *\bT\bG(E)_{x, f} \cap G_{x, 0}
\subseteq \lsupb *\bT G_{x, f}
		\cdot(\lsupb *\bT\bG(E)_{x, f_+}
			\cap G_{x, 0})$.
By Lemma \ref{lem:sloppy-vGvr},
$\lsupb *\bT\bG(E)_{x, f_+} \cap G_{x, 0}
= \lsub \bT G_{x, f_+}$
which, by Remark \ref{rem:little-split},
equals $\lsupb *\bT G_{x, f_+}$.
We are finished.
\end{proof}

\begin{lm}
\label{lem:more-vGvr-facts}
Suppose that
\begin{itemize}
\item
$\bG'$ is a connected compatibly filtered reductive $F$-subgroup
of \bG containing \bT;
\item
$I$ is some indexing set;
and
\item for each $i \in I$,
we have a $\Gal(E/F)$-invariant, concave function
$f_i$ on $\wtilde\Phi(\bG, \bT)$, constant on
$\wtilde\Phi(\bG', \bT)$,
such that there exists a $\Gal(E/F)$-invariant,
concave function $g_i$ for which
	\begin{itemize}
	\item $g_i(0) > 0$,
	\item $\lsub\bT\bG'(E)_{x, f_i(0)}$ normalizes
$\lsub\bT\bG(E)_{x, g_i}$\,,
	\item $f_i \le g_i$,
and
	\item $f_i = g_i$ off $\wtilde\Phi(\bG', \bT)$.
	\end{itemize}
\end{itemize}
Then
$$
\bigcap_{i \in I} \lsub\bT G_{x, f_i} = 
\lsub\bT G_{x, \max_{i \in I} f_i}\,.
$$
\end{lm}

Note that, if $\set{f_i}{i \in I}$ is any collection of
$\Gal(E/F)$-invariant, concave functions on
$\wtilde\Phi(\bG, \bT)$ such that $f_i(0) > 0$ for
$i \in I$, then we may take $\bG' = \bT$ and
$g_i = f_i$ in the above lemma.

\begin{proof}
By the definition of $\lsub\bT G_{x, (\cdot)}$\,, it suffices
to show that
$$
\bigcap_{i \in I} \lsub\bT\bG(E)_{x, f_i}
= \lsub\bT\bG(E)_{x, \max_{i \in I} f_i}\,.
$$
In fact, one containment being obvious, it suffices to show
that
$$
\bigcap_{i \in I} \lsub\bT\bG(E)_{x, f_i}
\subseteq
\lsub\bT\bG(E)_{x, \max_{i \in I} f_i}\,.
$$

Fix $i \in I$.
By Definition \ref{defn:vGvr-split},
$\lsub\bT\bG(E)_{x, f_i}$ is generated by
$\bG'(E)_{x, f_i(0)}$ and
$\lsub\bT\bG(E)_{x, g_i}$\,.
Since $\bG'(E)_{x, f_i(0)}$ normalizes
$\lsub\bT\bG(E)_{x, g_i}$\,,
we have that
$\lsub\bT\bG(E)_{x, f_i}
= \bG'(E)_{x, f_i(0)}
	\dotm\lsub\bT\bG(E)_{x, g_i}$\,.
By Lemma \ref{lem:gen-iwahori-factorization},
$$
\lsub\bT\bG(E)_{x, f_i}
= \bG'(E)_{x, f_i(0)}
	\dotm
	\prod U_{\alpha + g_i(\alpha)}
= \bG'(E)_{x, f_i(0)}
	\dotm
	\prod U_{\alpha + f_i(\alpha)}\,,
$$
where the unlabelled products, here and for the remainder
of this proof, run over
$\alpha \in \Phi(\bG, \bT)
	\smallsetminus \Phi(\bG', \bT)$.
(Here, the notation $\alpha + c$ is as in Lemma
\ref{lem:gen-iwahori-factorization}.)
In particular, for each $i \in I$,
$\lsub\bT\bG(E)_{x, f_i}$ lies in the image of
the multiplication map $\mult_{\bG'}$ of Lemma
\ref{lem:Zariski}, and, by loc.\ cit.,
its preimage under that map is precisely
$\bG'(E)_{x, f_i(0)}
\times
\prod U_{\alpha + f_i(\alpha)}$\,.
Thus $\bigcap_{i \in I} \lsub\bT\bG(E)_{x, f_i}$ lies in the
image of $\mult_{\bG'}$, and its preimage under that map is
precisely
$$
\bigcap_{i \in I} \bG'(E)_{x, f_i(0)}
\times
\prod \bigcap_{i \in I} U_{\alpha + f_i(\alpha)}
=
\bG'(E)_{x, (\max_{i \in I} f_i)(0)}
\times
\prod U_{\alpha + (\max_{i \in I} f_i)(\alpha)}\,;
$$
so $\bigcap_{i \in I} \lsub\bT\bG(E)_{x, f_i}$ is contained
in $\lsub\bT\bG(E)_{x, \max_{i \in I} f_i}$\,, as desired.
\end{proof}

\begin{lm}
\label{lem:center-comm}
Suppose that
\begin{itemize}
\item $\gamma \in Z(G^0) \cap \stab_G(\ox)$;
\item $\vec s = (s_0, \dotsc, s_d)$ and
$\vec t = (t_0, \dotsc, t_d)$ are admissible
sequences;
\item $\min\limits_{0 \le i \le d} s_i > 0$;
\item $\vec t \ge \vec s$;
\item
for all $0 < i \le d$,
either $t_i = \infty$ or
$t_i
< \min\set{s_i + \smash{\min\limits_{0 \le j \le d}} t_j, s_k + t_k}
	{i < k \le d}$;
\item $h \in \vG_{x, \vec s}\dotm G^0$;
and
\item $[\gamma, h] \in \vG_{x, \vec t}$\,.
\end{itemize}
Put $\vec r = (r_0, r_1, \dotsc, r_d)$, where
$r_0 = \min\set{s_k + t_k}{0 < k \le d}$,
and
$r_i = t_i$ for $0 < i \le d$.
Then $[\gamma, h] \in \vG_{x, \vec r}$\,.
\end{lm}

\begin{proof}
By construction, $\vec r$ is admissible.
We may assume that $r_0>t_0$, since otherwise
there is nothing to prove.
Let $f_{\vbG, \vec s}$, $f_{\vbG, \vec t}$,
and $f_{\vbG, \vec r}$ be the concave
functions associated to $\vec s$, $\vec t$,
and $\vec r$, respectively.
(See Definition \ref{defn:vGvr}.)
Then one checks, as in the proof of Corollary
\ref{cor:master-comm}, that
$f_{\vbG, \vec s} \vee f_{\vbG, \vec t}
\ge f_{\vbG, \vec t}$ and
$(f_{\vbG, \vec s} \vee f_{\vbG, \vec t})(\alpha)
> f_{\vbG, \vec t}(\alpha)$
whenever $f_{\vbG, \vec t}(\alpha) < \infty$;
and similarly for $f_{\vbG, \vec s} \vee f_{\vbG, \vec r}$.

Now write $h = h'h_0$, with $h_0 \in G^0$
and $h' \in \vG_{x, \vec s}$\,.  Then $[\gamma, h] = [\gamma, h']$.
By Lemma \ref{lem:gen-iwahori-factorization}, we may write
$$
h'
= \Bigl(\prod_{\alpha \in \wtilde\Phi(\bG, \bT) \smallsetminus
		\wtilde\Phi(\bG^0, \bT)}
	h_\alpha\Bigr)
\dotm
\Bigl(\prod_{\alpha \in \wtilde\Phi(\bG^0, \bT)}
	h_\alpha\Bigr),
$$
with $h_\alpha \in \lsub E U_{\alpha + f_{\vbG, \vec s}(\alpha)}$
for $\alpha \in \wtilde\Phi(\bG, \bT)$.
(Here, the notation $\alpha + c$ is as in Lemma
\ref{lem:gen-iwahori-factorization}.)
Then the commutator of $\gamma$ with $h'$
is the same as the commutator of $\gamma$ with
$\prod_{\alpha \in \wtilde\Phi(\bG, \bT) \smallsetminus
		\wtilde\Phi(\bG^0, \bT)}
	h_\alpha$.
(For the remainder of this proof, all products
should be understood as products over the same collection of
roots as above.)
In particular, this latter commutator lies in
$\vbG(E)_{x, \vec t}$\,.
That is,
$\prod \Int(\gamma)h_\alpha \equiv \prod h_\alpha
	\pmod{\vbG(E)_{x, \vec t}}$.
By Corollary \ref{cor:gen-iwahori-factorization-bij}, this means
that
$$
\Int(\gamma)h_\alpha \equiv h_\alpha
	\pmod{\lsub E U_{\alpha + f_{\vbG, \vec t}(\alpha)}
		= \lsub E U_{\alpha + f_{\vbG, \vec r}(\alpha)}
	\subseteq \vbG(E)_{x, \vec r}}
$$
for $\alpha$ as above,
hence (since
$\vbG(E)_{x, \vec s}$ normalizes $\vbG(E)_{x, \vec r}$)
that
$$
\Int(\gamma)\Bigl(\prod h_\alpha\Bigr)
= \prod \Int(\gamma)h_\alpha
\equiv
\prod h_\alpha
\pmod{\vbG(E)_{x, \vec r}};
$$
i.e.,
$[\gamma, h] = [\gamma, \prod h_\alpha] \in \vbG(E)_{x, \vec r}$\,.
By Lemma \ref{lem:sloppy-vGvr}, we have
$\vG_{x, \vec t} \subseteq G_{x, 0}$\,.
Since $[\gamma, h] \in \vG_{x, \vec t}$\,, we have
$[\gamma, h] \in \vbG(E)_{x, \vec r} \cap G_{x, 0}
		= \vG_{x, \vec r}$\,,
as desired.
\end{proof}

\subsection{Tame descent}
\label{sec:concave_tame}
We keep the notation of \S\ref{sec:concave_gen}.
In particular, $L/F$ satisfies the hypotheses introduced
before Lemma \ref{lem:gen-iwahori-factorization}.
For the remainder of the section,
suppose in addition that $L/F$ is tame;
that is, that Hypothesis \eqref{hyp:reduced} holds.
Let $N/F$ be the maximal unramified subextension of $L/F$.

Remember that, by Lemma \ref{lem:complete-subfield},
$F$ is an unramified extension of a complete
subfield.  If $F'$ is such a subfield, then
$L/F'$ is tame, hence separable, and the Galois closure
$\wtilde L$ of $L$ over $F'$ is again a
tame, discretely valued extension of $F'$.
Thus we may, and hence do, assume that 
there is some complete subfield $F'$ of $F$ such that
$L/F'$ is Galois.
However, we will not need this assumption until
Lemma \ref{lem:quotient-H1-iso}.

If $\chr \ff = 0$, then we assume in addition that $L$ is
strictly Henselian.  Since the strict Henselization of a
Galois extension is again Galois, this is compatible with
the assumption above.
Note that then $N$ is also strictly Henselian.
Thus, by the proof of Proposition IV.2.8
of \cite{serre:local-fields}, $L/N$ is cyclic.
However, we will not need this assumption until
Lemma \ref{lem:H1-of-pro-p}.

\begin{rk}
\label{rem:vGvr-facts}
Suppose that there is a connected, reductive,
compatibly filtered $F$-subgroup \bH of \bG, containing \bT,
such that $f$ is identically $\infty$ off $\wtilde\Phi(\bH, \bT)$.
Denote by $M$ the maximum value of $f$ on $\wtilde\Phi(\bH, \bT)$.
Clearly, $\lsub\bT\bG(E)_{x, f} \supseteq \bH(E)_{x, M}$\,,
so Lemmata \ref{lem:field-descent} and
\ref{lem:group-descent} give
$$
\lsub\bT G_{x, f} = \lsub\bT\bG(E)_{x, f} \cap G_0
\supseteq \bH(E)_{x, M} \cap H_0 \supseteq H_{x, M}\,.
$$
Now suppose that $f$ is non-negative, and put $m = \min f$.
If $\bH = \bG$ and $m = 0$, then, by Lemma \ref{lem:sloppy-vGvr},
$\lsub\bT G_{x, f} \subseteq G_{x, 0} = H_{x, m}$\,.
If $m > 0$, then let \ol m be the concave function on
$\wtilde\Phi(\bG, \bT)$ which takes the value $m$ on
$\wtilde\Phi(\bH, \bT)$, and is $\infty$ off it.
By Lemma \ref{lem:gen-iwahori-factorization},
$\lsub\bT\bG(L\unram)_{x, \ol m} = \bH(L\unram)_{x, m}$\,,
so Lemma \ref{lem:field-descent} gives
$$
\lsub\bT G_{x, f} \subseteq \lsub\bT G_{x, \ol m}
\subseteq \lsub\bT\bG(L\unram)_{x, \ol m}^{\Gal(L\unram/F)}
= \bH(L\unram)_{x, m}^{\Gal(L\unram/F)}
= H_{x, m}\,.
$$
In particular, if $f$ takes the value $c \in \tR_{\ge 0}$
everywhere on $\wtilde\Phi(\bH, \bT)$,
and $c > 0$ or $\bH = \bG$, then
$\lsub\bT G_{x, f} = H_{x, c}$\,.
\end{rk}

\begin{lm}
\label{lem:shallow-comm}
Suppose that
$s' \in \tR_{\ge 0}$
and
$\vec s$ is an admissible sequence.
Then
$[G^0_{x, s'}, \vG_{x, \vec s}]
	\subseteq \vG_{x, s' + \vec s}$\,.
\end{lm}

\begin{proof}
Put
$\vec s^{(1)} = (s', \infty, \ldots, \infty)$
and $\vec s^{(2)} = \vec s$.
By Remark \ref{rem:vGvr-facts}, we have
$G^0_{x, s'} = \vG_{x, \vec s^{(1)}}$\,.
Now the result follows from
Corollary \ref{cor:master-comm}.
\end{proof}

\begin{lm}
\label{lem:tame-descent}
Suppose that $K/F$ is a discretely valued tame extension.
If $K/F$ is unramified or $f(0) > 0$, then
$\lsub\bT\bG(K)_{x, f} \cap G
= \lsub\bT G_{x, f}$\,.
\end{lm}

\begin{proof}
It is clear that the right-hand side is contained in the
left, so we need only show the reverse containment.
By Lemma \ref{lem:big-split}, we may, and hence do, assume
that $K \subseteq E$.
Note that the composite extension $K L/F$ remains tame, so
we may, and hence do, assume further that $K \subseteq L$.
If $K/F$ is unramified, then let
$x' \in \BB(\bG, F)$ be an $(x, f)$-positive point.
(For example, we could take $x' = y$, where $(y, h)$ is the
pair introduced at the beginning of
\S\ref{sec:concave_gen}.)
By Lemmata \ref{lem:field-descent} and \ref{lem:sloppy-vGvr}
(applied twice), we have that
$$
\lsub\bT\bG(K)_{x, f} \cap G
= \lsub\bT\bG(E)_{x, f} \cap \bG(K)_{x', 0} \cap G
= \lsub\bT\bG(E)_{x, f} \cap G_{x', 0}
= \lsub\bT G_{x, f}\,,
$$
as desired.

If $f(0) > 0$, then, since $L/K$ is tame, we have by
Corollary \ref{cor:very-pos-point} that there is a point
$x' \in \BB(\bG, F)$ with
$\lsub\bT G_{x, f}
= \lsub\bT\bG(L)_{x, f} \cap G_{x', 0+}$
and
$\lsub\bT\bG(K)_{x, f}
= \lsub\bT\bG(L)_{x, f} \cap \bG(K)_{x', 0+}$
Since $K/F$ is tame, we have by Lemma
\ref{lem:field-descent} that
$G_{x', 0+} = \bG(K)_{x', 0+} \cap G$.
Combining these three equations, we find that
$\lsub\bT G_{x, f} = \lsub\bT\bG(K)_{x, f} \cap G$,
as desired.
\end{proof}

\begin{rk}
If there is a tame extension over which \vbG splits,
then Lemma \ref{lem:tame-descent} shows that the definition
of $\lsub\bT G_{x,f}$ here coincides with the one in
\cite{yu:supercuspidal}*{\S 13} when $f(0) > 0$.
\end{rk}

Until the end of the proof of Lemma \ref{lem:H1-of-pro-p}
only, put $p = \chr \ff$ and call a $\Gal(L/F)$-group \mc G
\emph{$p$-filtered}
if it possesses a
filtration $(\mc G_i)_{i \in \Z_{\ge 0}}$ by closed normal
$\Gal(L/F)$-subgroups such that
$\mc G = \mc G_0$,
$\mc G = \varprojlim \mc G/\mc G_i$, and
\begin{itemize}
\item 
$p > 0$ and $\mc G_i/\mc G_{i + 1}$ is an
Abelian $p$-torsion group;
or
\item
$p = 0$ and
$\mc G_i/\mc G_{i + 1}$ is a \Q-vector space on which $\Gal(L/F)$
acts linearly
\end{itemize}
for all $i \in \Z_{\ge 0}$.
Note that, if $p > 0$
(but not necessarily if $p = 0$),
then a $\Gal(L/F)$-subgroup or $\Gal(L/F)$-quotient of a
$p$-filtered group is again $p$-filtered.

\begin{rk}
\label{rem:p-filt}
If $p > 0$, then
$\mc G := \bG(L)_{x', 0{+}}$ is clearly $p$-filtered (with
$\mc G_i = \bG(L)_{x', 2^i\varepsilon}$\,, say, for
$\varepsilon \in \R_{> 0}$ sufficiently small)
for any $x' \in \BB(\bG, F)$;
so, if $f(0) > 0$, then Corollary \ref{cor:very-pos-point}
shows that $\lsub\bT\bG(L)_{x, f}$ is $p$-filtered.
In particular, in case $\bG = \bT$, we see that
$\bT(L)_c$\,, hence $\bT(L)_{c:d}$\,, is $p$-filtered for
any $c, d \in \tR_{\ge 0}$ with $c \le d$.

If $p = 0$, then \bT (indeed, any $L$-torus)
is a tame $L$-torus, so,
by Proposition 5.5 of \cite{yu:models},
$\bT(L)_{c:d}$ is
$p$-filtered for any $c, d \in \tR_{> 0}$ with
$c \le d$.
Then,
by Lemma \ref{lem:concave-in-R} and Corollary
\ref{cor:gen-iwahori-factorization-iso},
$\lsub\bT\bG(L)_{x, f}$ is $p$-filtered if
$(f \vee f)(\alpha) > f(\alpha)$ whenever
$f(\alpha) < \infty$.
\end{rk}

Recall that $N/F$ is the maximal
unramified subextension of $L/F$.

\begin{lm}
\label{lem:H1-of-pro-p}
%%Note that this lem:H1-of-pro-p used to be cor:H1-of-pro-p.
%%The old H1-of-pro-p is no longer needed.
Suppose that $A$ is a $\Gal(L/F)$-group, and $B$ a closed
$\Gal(L/F)$-subgroup.  Suppose further that
$B$ and $A/B$ are $p$-filtered.
Then the natural map
$$
H^1(N/F, A^{\Gal(L/N)}/B^{\Gal(L/N)})
\to
H^1(L/F, A/B)
$$
is a bijection.
\end{lm}

\begin{proof}
Note that it suffices to prove that
$H^1(L/N, C) = \sset 0$ when $C$ is $p$-filtered.
Indeed, once we have done so, we will have the exact
sequences
$$
B^{\Gal(L/N)} \to A^{\Gal(L/N)}
\to (A/B)^{\Gal(L/N)} \to H^1(L/N, B) = 0
$$
(by Proposition I.5.5.38 of \cite{serre:galois}),
so that $(A/B)^{\Gal(L/N)} = A^{\Gal(L/N)}/B^{\Gal(L/N)}$,
and
$$
0 \to H^1(N/F, (A/B)^{\Gal(L/N)}) \to H^1(L/F, A/B)
\to H^1(L/N, A/B) = 0
$$
(by \S 5.8(a) of loc.\ cit.).

By Proposition 2.8 of \cite{yu:supercuspidal}
and
Proposition I.5.5.38 of \cite{serre:galois},
it suffices to prove that $H^1(L/N, C) = \sset 0$ when
$p > 0$ and $C$ is an Abelian $p$-torsion group,
or
$p = 0$ and $C$ is a \Q-vector space on which
$\Gal(L/N)$ acts linearly.

In the former case, note that every element of the
cohomology group $H^1(L/N, C)$ has $p$-power order.
On the other hand, by Proposition I.2.4.9 of
\cite{serre:galois}, every element has order dividing
the cardinality of $\Gal(L/N)$, which is indivisible by $p$.
Thus, $H^1(L/N, C) = \sset 0$.

In the latter case, recall that our assumptions on $L$ imply that $L/N$
is cyclic,
say with generator $\sigma$.
Since $\sigma$ has finite order, say $e$, it acts semisimply on $C$,
so $C = \ker(\sigma - 1) + \im(\sigma - 1)$.
Now $\sum_{i = 0}^{e - 1} \sigma$ vanishes on $\im(\sigma - 1)$, 
and acts on $\ker(\sigma - 1)$ as multiplication by $e$,
so $\ker(\sigma - 1) = \im \sum_{i = 0}^{e - 1} \sigma$;
that is,
$H^1(L/N, C) = H^1(\langle\sigma\rangle, C) = \sset 0$.
\end{proof}

\begin{lm}
\label{lem:quotient-H1-iso}
With the notation and hypotheses of Corollary
\ref{cor:gen-iwahori-factorization-iso},
suppose further that
$f_2(\alpha) < \infty$ whenever
$\alpha \in \Phi(\bG, \bT)$ and $f_1(\alpha) < \infty$.
Then the natural map
$$
H^1(N/F, \bT(N)_{f_1(0):f_2(0)})
\to H^1(L/F, \lsub\bT\bG(L)_{x, f_1:f_2})
$$
is an isomorphism.
\end{lm}

Recall that $N/F$ is the maximal unramified subextension of
$L/F$.

\begin{proof}
By Remark \ref{rem:p-filt} and Lemma \ref{lem:H1-of-pro-p},
$H^1(N/F, \bT(N)_{f_1(0):f_2(0)})
\cong H^1(L/F, \bT(L)_{f_1(0):f_2(0)})$;
so it suffices to show that the natural map
$H^1(L/F, \bT(L)_{f_1(0):f_2(0)})
\to H^1(L/F, \lsub\bT\bG(L)_{x, f_1:f_2})$ is an
isomorphism.
By Corollary \ref{cor:gen-iwahori-factorization-iso}
(with notation $a + c$ as in Lemma
\ref{lem:gen-iwahori-factorization}), the
multiplication map
$$
\bT(L)_{f_1(0):f_2(0)}
\times
\prod_{a \in \Phi(\bG, \bS^\sharp)}
	U_{(a + f_1(a)):(a + f_2(a))}
\to \lsub\bT\bG(L)_{x, f_1:f_2}
$$
is a $\Gal(L/F)$-equivariant isomorphism, so the natural map
\begin{multline*}
H^1(L/F, \bT(L)_{f_1(0):f_2(0)})
\times
H^1(L/F, \prod U_{(a + f_1(a)):(a + f_2(a))}) \\
\to H^1(L/F, \lsub\bT\bG(L)_{x, f_1:f_2})
\end{multline*}
is also an isomorphism.
By \S\ref{sec:unip-exp}, there is a $\Gal(L/F)$-equivariant
isomorphism
$\prod U_{a + f_1(a)} \cong \bigoplus \mf u_{a + f_1(a)}$
which restricts to a $\Gal(L/F)$-equivariant isomorphism
$\prod U_{a + f_2(a)} \cong \bigoplus \mf u_{a + f_2(a)}$,
hence induces a $\Gal(L/F)$-equivariant isomorphism
$\prod U_{(a + f_1(a)):(a + f_2(a))}
\cong \bigoplus \mf u_{(a + f_1(a)):(a + f_2(a))}$.
Put
$$
\mc L_j(L) := \bigoplus \mf u_{a + f_j(a)}
\subseteq \Lie(\bG)(L)
$$
for $j = 1, 2$.

We need only show that
$H^1(L/F, \mc L_1(L)/\mc L_2(L)) = \sset 0$.
Put
$$
\pmb{\mf u}
:= \bigoplus_{\substack{
	a \in \Phi(\bG, \bS^\sharp) \\
	f_1(a) < \infty
}}
	\Lie(\bU_a).
$$
Remember that $F$ contains a complete
subfield $F'$ such that $L/F'$ is Galois.
By Lemma \ref{lem:complete-subfield}, we may, and hence do,
assume in addition that
\begin{itemize}
\item $F/F'$ is unramified,
\item \bG, \bT, $\bS^\sharp$, and $\pmb{\mf u}$ are defined over $F'$,
and
\item $x \in \BB(\bG, F')$.
\end{itemize}
By another application of Lemma \ref{lem:complete-subfield},
there is a finite subextension $N'/F'$ of $N/F'$ such that
$N/N'$ is unramified.
By a third application of Lemma \ref{lem:complete-subfield},
\begin{equation}
\tag{$*$}
H^1(L/F, \mc L_1(L)/\mc L_2(L))
= \varinjlim H^1(L/F, \mc L_1(L')/\mc L_2(L')),
\end{equation}
where $\mc L_j(L') := \mc L_j(L) \cap \Lie(\bG)(L')$ for
$j = 1, 2$,
and $L'$ runs over the collection of finite extensions of
$N'$ such that $L'/F'$ is Galois and $\bS^\sharp$ is
$L'$-split.
(Note that, in this setting, $\Gal(L/F)$ actually acts on
$L'$, hence on $\mc L_1(L')$ and $\mc L_2(L')$; so the
cohomology above makes sense.)

\def\funfld#1{{\ensuremath{\wtilde #1'}}\xspace}
\def\funF{\funfld F}
\def\funN{\funfld N}
Fix $L'$ as above.
Put $\funN := L' \cap N$ and $\funF := L' \cap F$.
Then it is easy to check that
\funN and \funF satisfy the hypotheses on $N'$ and $F'$
above,
that $L'/\funF$ is tame,
and
that $\funN/\funF$ is its maximal unramified subextension.
By Lemma \ref{lem:H1-of-pro-p},
\begin{multline}
\tag{$\dag$}
H^1(L/F, \mc L_1(L')/\mc L_2(L')) \\
\cong H^1(N/F, \mc L_1(L')^{\Gal(L/N)}/\mc L_2(L')^{\Gal(L/N)}) \\
= H^1(N/F, \mc L_1(\funN)/\mc L_2(\funN)).
\end{multline}
By \cite{serre:galois}*{\S I.5.8(a)}, we have an exact
sequence
\begin{multline*}
1
\to H^1(N/\funN F, \mc L_1(\funN)/\mc L_2(\funN))
\to H^1(N/F, \mc L_1(\funN)/\mc L_2(\funN)) \\
\to H^1(\funN F/F, \mc L_1(\funN)/\mc L_2(\funN)).
\end{multline*}
Since $\Gal(N/\funN F)$ acts trivially on
$\mc L_1(\funN)/\mc L_2(\funN)$, we have
$$
H^1(N/\funN F, \mc L_1(\funN)/\mc L_2(\funN))
\cong
\Hom\bigl(\Gal(N/\funN F), \mc L_1(\funN)/\mc L_2(\funN)\bigr).
$$
Note that $\Gal(\funN F/F) \cong \Gal(\funN/\funF)$,
and that $\mc L_j(\funN)$ is a $\Gal(\funN/\funF)$-stable lattice in
$\pmb{\mf u}(\funN) \cong \pmb{\mf u}(\funF) \otimes_\funF \funN$
for $j = 1, 2$.
By Lemma 2.9 of \cite{yu:supercuspidal} (the proof of which
uses only completeness of \funN, not local compactness),
we have
\begin{equation}
\tag{$\ddag$}
H^1(\funN F/F, \mc L_1(\funN)/\mc L_2(\funN)) = \sset 0.
\end{equation}
Thus
\begin{equation}
\tag{$**$}
H^1(L/F, \mc L_1(L')/\mc L_2(L')
\cong
\Hom\bigl(\Gal(N/\funN F), \mc L_1(\funN)/\mc L_2(\funN)\bigr).
\end{equation}

By ($*$) and ($**$), since
$\bigcup_{L'} (L' \cap N)F = N$,
we have
\begin{multline*}
H^1(L/F, \mc L_1(L)/\mc L_2(L)) \\
\cong
\varinjlim \Hom\bigl(
	\Gal(N/(L' \cap N)F), \mc L_1(L' \cap N)/\mc L_2(L' \cap N)
\bigr) = \sset 0.\qedhere
\end{multline*}
\end{proof}

\begin{cor}
\label{cor:quotient-H1-triv}
With the notation and hypotheses of Lemma
\ref{lem:quotient-H1-iso},
$$
H^1(L/F, \lsub\bT\bG(L)_{x, f_1:f_2}) = \sset 0.
$$
\end{cor}

\begin{proof}
By Lemma \ref{lem:quotient-H1-iso} and
Proposition I.5.5.38 of \cite{serre:galois},
it suffices to show that
$H^1(L/F, \bT(L)_{x, c:c{+}}) = \sset 0$ for all $c \in \R_{> 0}$.
Since Remark \ref{rem:p-filt} and Lemma \ref{lem:H1-of-pro-p} show that
$H^1(N/F, \bT(N)_{x, c:c{+}})
\cong H^1(L/F, \bT(L)_{x, c:c{+}})$,
the desired equality is Hypothesis \eqref{hyp:torus-H1-triv}.
\end{proof}

The next result
is the analogue of Proposition 13.4 of
\cite{yu:supercuspidal}.

\begin{pn}
\label{prop:H1-iso}
Suppose that $f(0) > 0$
and $(f \vee f)(\alpha) > f(\alpha)$
whenever $f(\alpha) < \infty$.
Then $H^1(L/F, \lsub\bT\bG(L)_{x, f}) = \sset 0$.
\end{pn}

\begin{proof}
As in the proof of Lemma \ref{lem:quotient-H1-iso}, it
suffices to show that
\begin{equation}
\tag{$\dag$}
H^1(L/F, \lsub\bT\bG(L')_{x, f})
\cong H^1(N/F, \lsub\bT\bG(N')_{x, f})
\end{equation}
and
\begin{equation}
\tag{$\ddag$}
H^1(N'/F', \lsub\bT\bG(N')_{x,f}) = \sset 0
\end{equation}
for any tower of complete subfields $L'/N'/F'$ of $L$ such that
\begin{itemize}
\item $L/F'$ and $L'/F'$ are Galois,
\item $N' = L' \cap N$ and $F' = L' \cap F$,
\item $F/F'$, $N/N'$, and $L/L'$ are unramified,
\item \bG, \bT, $\bS^\sharp$, and $\pmb{\mf u}$ are defined
over $F'$,
where
$\pmb{\mf u} := \bigoplus_{\smash{\substack{
	a \in \Phi(\bG, \bS^\sharp) \\
	f(a) < \infty
}}} \Lie(\bU_a)$,
\item $\bS^\sharp$ is split over $L'$,
and
\item $x \in \BB(\bG, F')$.
\end{itemize}

Equation ($\dag$) follows from Remark \ref{rem:p-filt}
and Lemma \ref{lem:H1-of-pro-p}.
As in the proof of Corollary
\ref{cor:gen-iwahori-factorization-bij}, we may, and hence
do, assume that there is some $\varepsilon > 0$ such that
$f \vee f \ge f + \varepsilon$.
Then equation ($\ddag$) follows (as in the proof of Proposition
13.4 of \cite{yu:supercuspidal}) from
Proposition I.5.5.38 of \cite{serre:galois},
Lemma 2.8 of \cite{yu:supercuspidal},
and
Corollary \ref{cor:quotient-H1-triv}.
\end{proof}

The next result is the analogue of Lemma 13.3 of
\cite{yu:supercuspidal}.

\begin{pn}
\label{prop:heres-a-gp}
Suppose that, for $j = 1, 2$,
$f_j$ is a $\Gal(E/F)$-invariant, concave function
on $\wtilde\Phi(\bG, \bT)$
such that $f_j(0) > 0$,
and $(f_j \vee f_j)(\alpha) > f_j(\alpha)$
whenever $f_j(\alpha) < \infty$.
Then
$$
\bigl( \lsub\bT\bG(L)_{x, f_1}\cdot\lsub\bT\bG(L)_{x, f_2} \bigr)
	\cap G
= \lsub\bT G_{x, f_1}\cdot\lsub\bT G_{x, f_2}.
$$
If moreover $\min\sset{f_1, f_2}$ is concave
and
$\lsub\bT\bG(L)_{x, f_1}\cdot\lsub\bT\bG(L)_{x, f_2}$ is a
group,
then
$$
\lsub\bT G_{x, \min\sset{f_1, f_2}}
= \lsub\bT G_{x, f_1}\cdot\lsub\bT G_{x, f_2}\,.
$$
\end{pn}

\begin{proof}
As in the proof of Corollary
\ref{cor:gen-iwahori-factorization-bij}, we may, and hence
do, assume that there is some $\varepsilon \in \R_{> 0}$ such
that $f_j \vee f_j \ge f_j + \varepsilon$ for $j = 1, 2$.
Then
$\max\sset{f_1, f_2} \vee \max\sset{f_1, f_2}
\ge \max\sset{f_1 + \varepsilon, f_2 + \varepsilon}
= \max\sset{f_1, f_2} + \varepsilon$.

Consider the groups
$A = \lsub\bT\bG(L)_{x, \max\sset{f_1, f_2}}$
and
$B = \lsub\bT\bG(L)_{x, f_1} \times \lsub\bT\bG(L)_{x, f_2}$\,.
By identifying $A$ with its image under the diagonal map,
we may regard $A$ as a subgroup of $B$.
By Proposition I.5.4.36 of \cite{serre:galois}, we have an
exact sequence (of pointed sets)
$$
B^{\Gal(L/F)} \to (B/A)^{\Gal(L/F)}
\to H^1(L/F, A).
$$
By Proposition \ref{prop:H1-iso}, $H^1(L/F, A) = \sset 1$,
so the natural map $B^{\Gal(L/F)} \to (B/A)^{\Gal(L/F)}$
is surjective.
By Lemma \ref{lem:more-vGvr-facts},
$\lsub\bT\bG(L)_{x, f_1} \cap \lsub\bT\bG(L)_{x, f_2}
= \lsub\bT\bG(L)_{x, \max\sset{f_1, f_2}}$,
so the map $(b_1, b_2) \mapsto b_1 b_2\inv$ identifies
$(B/A)^{\Gal(L/F)}$
with
$(\lsub\bT\bG(L)_{x, f_1}\dotm\lsub\bT\bG(L)_{x, f_2})
	^{\Gal(L/F)}
= (\lsub\bT\bG(L)_{x, f_1}\dotm\lsub\bT\bG(L)_{x, f_2})
	\cap G$.
By Lemma \ref{lem:tame-descent},
$B^{\Gal(L/F)}
= \lsub\bT G_{x, f_1} \times \lsub\bT G_{x, f_2}$\,.
This proves the first statement of the lemma.

Now suppose that $\min\{f_1,f_2\}$ is concave.
By Lemma \ref{lem:gen-iwahori-factorization},
$\lsub\bT\bG(L)_{x, \min\sset{f_1, f_2}}$ is generated by
$\lsub\bT\bG(L)_{x, f_1}$ and $\lsub\bT\bG(L)_{x, f_2}$.
If $\lsub\bT\bG(L)_{x, f_1}\dotm\lsub\bT\bG(L)_{x, f_2}$
is a group, then in fact
$\lsub\bT\bG(L)_{x, \min\sset{f_1, f_2}}
= \lsub\bT\bG(L)_{x, f_1}\dotm\lsub\bT\bG(L)_{x, f_2}$.
Since
$\lsub\bT\bG(L)_{x, \min\sset{f_1, f_2}} \cap G
= \lsub\bT G_{x, \min\sset{f_1, f_2}}$ by Lemma
\ref{lem:tame-descent}, the second statement follows.
\end{proof}

\section{Normal approximations: basic definitions}
\label{sec:normal}

\begin{defn}
\label{defn:good}
\indexmem{\GG^G_d}
Let $\GG^G_0$ denote the set of elements $\gamma \in G$
such that $\ogamma \in \tG$ is absolutely semisimple (in the
sense of Definition \ref{defn:top-F-ss}).

For $d > 0$,
$\gamma \in G$ is
\emph{good of depth $d$}
\indexme{good element!of depth $d > 0$}%
if there is a tame-modulo-center torus \bS
such that
$\gamma \in S_d \smallsetminus S_{d+}$
and
$\alpha(\gamma) = 1$ or $\ord(\alpha(\gamma) - 1) = d$
for all $\alpha \in \Phi(\bG, \bS)$.
(This is analogous to the definition of a good element of a Lie
algebra in Definition 2.2.4 of \cite{adler:thesis}.)
Let $\GG^G_d$ denote the set of such elements.
\end{defn}

\begin{rem}
Note that an element of $\GG^G_d$ with $d > 0$ has depth
precisely $d$, but an element of $\GG^G_0$ may have positive
depth (in which case it belongs to $Z(G)$), or may not even
belong to $G_0$\,.
For every discretely valued
separable extension $E/F$
and $d \in \R_{\ge 0}$, we have
$\GG^G_d \subseteq \GG^{\bG(E)}_d$.
\end{rem}

\begin{dn}
\label{defn:S-is-good}
\indexmem{\text{\bGd{\bG(E)}}}
Let $E/F$ be a discretely valued algebraic extension.
We say that an $E$-torus $\bS \subseteq \bG$ has property \bGd{\bG(E)}
if, for all $d > 0$,
every coset in $\bS(E)_{d:d+}$ intersects
$\GG^{\bG(E)}_d \cup \sset 1$.
\end{dn}

\begin{defn}
\label{defn:funny-centralizer}
\indexme{good sequence}%
A collection $\ugamma = (\gamma_i)_{0\leq i < r'}$
of elements of $G$, where
$r'\in\tR_{\geq0}$ and
$\gamma_i\in\GG^G_i\cup\sset 1$ for $0 \le i < r'$,
	%%Built in the requirement that everybody belongs
	%%to a tame torus here, rather than in the definition
	%%of a normal $r$-approximation, so that the remark
	%%below about \CC\bG r(\ugamma) being
	%%compatibly filtered is correct.
is called a \emph{good sequence} (in $G$) if there is a tame
$F$-torus \bS in \bG such that $\gamma_i \in S$ for all
$0 \le i < r'$.
(We will often omit those terms $\gamma_i$ that are
equal to $1$ from the notation.)
For $r \in \tR$ with $r \le r'$, put
\indexmem{\CC{\bG}{r}(\ugamma)}%
\indexmem{\ZZ{\bG}{r}(\ugamma)}%
\indexmem{\CC G  r(\ugamma)}%
\indexmem{\ZZ G r(\ugamma)}%
\begin{align*}
\CC\bG{r}(\ugamma) & := \bigl(\bigcap_{0\leq i < r} C_\bG(\gamma_i)\bigr)\conn,
\\
\ZZ\bG{r}(\ugamma) & := Z(\CC\bG{r}(\ugamma)), \\
\CC G r(\ugamma) & := \CC\bG r(\ugamma)(F), \\
\ZZ G r(\ugamma) & := \ZZ\bG r(\ugamma)(F).
\end{align*}
By convention, $\CC\bG r(\ugamma) = \bG$,
so $\ZZ\bG r(\ugamma) = Z(\bG)$,
for $r \le 0$.
\end{defn}

Since the intersection defining $\CC\bG r(\ugamma)$ may be
taken over a finite set,
we see by repeated applications of
Propositions 9.1(1) and 13.19 of \cite{borel:linear}
that
$\CC\bG r(\ugamma)$ is reductive and defined over $F$.
By Theorem 18.2(ii) of \cite{borel:linear},
$\ZZ\bG r(\ugamma)$ is defined over $F$.
Thus, the definitions of $\CC G r(\ugamma)$
and $\ZZ G r(\ugamma)$ make sense.

Note that $\CC\bG r(\ugamma)$ is a compatibly filtered
$F$-subgroup of \bG, by Proposition
\ref{prop:compatibly-filtered-tame-rank} and the definition
of a good sequence.

For later use, it will be convenient to know that the groups
$\CC\bG r(\gamma)$ descend well to Levi subgroups of \bG.

\begin{lm}
\label{lem:funny-centralizer-descends}
Suppose that
\begin{itemize}
\item $E/F$ is a discretely valued separable field extension,
\item \bM is an $E$-Levi $F$-subgroup of \bG,
and
\item $\ugamma = (\gamma_i)_{0 \le i < r'}$
is a good sequence in $M$ and $G$.
\end{itemize}
Then $\CC\bM r(\ugamma)$ equals
$\CC\bG r(\ugamma) \cap \bM$,
and is an $E$-Levi $F$-subgroup of $\CC\bG r(\ugamma)$
for $0 \le r \le r'$.
\end{lm}

\begin{proof}
Put $\bH = \CC\bG r(\ugamma)$.
Clearly,
$(\bH \cap \bM)\conn = \CC\bM r(\ugamma)$.
Let \bS be the maximal $E$-split torus in the center of \bM,
so that $\bM = C_\bG(\bS)$.
Since \ugamma is a good sequence in $M$, we have that
$\bS \subseteq \bH$.
Then $\bH \cap \bM = C_\bH(\bS)$ is an $E$-Levi $F$-subgroup of
\bH.
In particular, it is connected, so that
$\bH \cap \bM = \CC\bM r(\ugamma)$, as desired.
\end{proof}

\begin{defn}
\label{defn:bracket}
Let \ugamma and $r \in \tR$ be as in Definition
\ref{defn:funny-centralizer}.
Fix $x \in \BB(\CC\bG r(\ugamma), F)$ and $j \in \tR$.
For this definition only, write
\begin{equation*}
\vbG(j) = (\CC\bG{r - i}(\ugamma))_{0 < i < j}\,;
\quad
\vec r(j) = (i)_{0 < i < j}\,;
\quad
\text{and}
\quad
\vec s(j) = (i/2)_{0 < i < j}
\end{equation*}
(where $i$ runs over the indicated elements in \R,
not \tR).
For $x \in \BB(\CC\bG r(\ugamma), F)$, put
\indexmem{[\ugamma; x, r]}%
\indexmem{\udc}%
\indexmem{[\ugamma; x, r]^{(j)}}%
\indexmem{\udc^{(j)}}%
\begin{equation*}
[\ugamma; x, r]^{(j)} := \vG(j)_{x, \vec r(j)}
\quad\text{and}\quad
\udc^{(j)} := \vG(2j)_{x, \vec s(2j)}.
\end{equation*}
(Here, we are using Remark \ref{rem:infinite-vec} to handle
``vectors'' $\vec r(j)$ and $\vec s(j)$ with infinitely many
entries.
Note also that $[\ugamma; x, r]^{(j)}$ and $\udc^{(j)}$ are just open
subgroups of $G$, not sequences of such subgroups.)
Put also $[\ugamma; x, r] = [\ugamma; x, r]^{(\infty)}$ and
$\udc = \udc^{(\infty)}$.
\end{defn}

\begin{rk}
\label{rem:bracket-facts}
If necessary, we will indicate the dependence of
$[\ugamma; x, r]$ on $G$ by denoting it by
$[\ugamma; x, r]_G$\,,
and similarly for $\udc$, $[\ugamma; x, r]^{(j)}$, and
$\udc^{(j)}$.
By Proposition \ref{prop:heres-a-gp},
each of the groups defined above has a more concrete description.
For example,
$[\ugamma; x, r]
= \prod_{0 < i \le r} \CC G{r - i}(\ugamma)_{x, i}$\,.
These concrete descriptions make it easy to see the
following (although they can also be verified directly).
\begin{inc_enumerate}
\item
\label{rem:bracket-facts-containment}
If $h + j \le r+$, and $r = r+$ or $h \in \R$ or $j \in \R$, then
\begin{gather*}
[\ugamma; x, r]^{(j)} \subseteq \CC G h(\gamma)_{x, 0+}\,, \\
\udc^{(j)} \subseteq [\ugamma; x, h]^{(j)}, \\
G_{x,r} \subseteq [\ugamma;x,r] \subseteq G_{x,0+} \, , \\
\intertext{and}
G_{x,r/2} \subseteq \llbracket\ugamma;x,r\rrbracket \subseteq G_{x,0+} \, .
\end{gather*}
\item
\label{rem:bracket-facts-field-descent}
By Lemma \ref{lem:sloppy-vGvr}
and Remark \ref{rem:little-split},
$[\ugamma; x, r]_{\bG(E)} \cap G_{x, 0}
= [\ugamma; x, r]_{\bG(E)} \cap G_{x, 0+}
= [\ugamma; x, r]_G$
for any discretely valued separable extension $E/F$.
\item
\label{rem:bracket-facts-torus-descent}
Suppose that \bT is an $F$-torus such that
$(\bT, \bG', \bG)$ is a tame reductive $F$-sequence,
$\gamma \in T$, and $x \in \BB(\bT, F)$.
Let $E/F$ be a splitting field for \bT,
$f_1$ the concave function appearing in the
definition of $[\ugamma; x, r]_{\bG(E)}$, and
$f_2$ the concave function appearing in the definition of
$(\bG', \bG)(E)_{x, (0{+}, \infty)}$.
(See Definition \ref{defn:vGvr-split}.)
By Remark \ref{rem:vGvr-facts},
$(G', G)_{x, (0{+}, \infty)} = G'_{x, 0+}$\,.
Since \bG and $\bG'$ are compatibly filtered over $F$, we
have that $G_{x, 0+} \cap G' = G'_{x, 0+}$\,.
Since, as above, $[\ugamma; x, r]_G \subseteq G_{x, 0+}$\,,
we have
$[\ugamma; x, r]_G \cap G'
= [\ugamma; x, r]_G \cap G'_{x, 0+}$\,.
By Lemma \ref{lem:more-vGvr-facts},
$[\ugamma; x, r]_G \cap G'_{x, 0{+}}
= \lsub\bT G_{x, \max\sset{f_1, f_2}}$\,.
It is easy to verify that
$\lsub\bT G_{x, \max\sset{f_1, f_2}} = [\ugamma; x, r]_{G'}$.
Thus,
$$
[\ugamma; x, r]_G \cap G'
= [\ugamma; x, r]_{G'}.
$$
The analogous facts for
$\udc$, $[\ugamma; x, r]^{(j)}$, and $\udc^{(j)}$
hold, with similar proofs.
\item
\label{rem:bracket-facts-decomp}
By an argument as in Remark
\ref{rem:bracket-facts-torus-descent}, we see that,
for $j > 0$ and $h \in \tR$,
$[\ugamma; x, r] \cap \CC G h(\ugamma)_{x, j}$ is the
group associated to the vector of groups
$(\CC\bG{\max\sset{r - i, h}}(\ugamma))_{j \le i}$ and the vector of
depths $(i)_{j \le i}$; and similarly for $\udc^{(j)}$.
Thus, by Proposition \ref{prop:heres-a-gp}, if $h + j \le r$
and $h^\sharp + 2j \le r$, then
\begin{gather*}
[\ugamma; x, r]^{(j)}([\ugamma; x, r] \cap \CC G h(\ugamma)_{x, j})
= [\ugamma; x, r]\\
\intertext{and}
\udc^{(j)}(\udc \cap \CC G{h^\sharp}(\ugamma)_{x, j}) = \udc. \\
\end{gather*}
\end{inc_enumerate}
\end{rk}

\begin{defn}
\label{defn:r-approx}
Suppose that
\begin{itemize}
\item $\gamma \in G$,
\item $\ugamma = (\gamma_i)_{0\leq i < r'}$
is a good sequence,
and
\item $r \in \tR$ with $0 \le r \le r'$.
\end{itemize}
Then \ugamma is an \emph{$r$-approximation to $\gamma$ (in $G$)}
\indexme{approximation@$r$-approximation}%
if there is a point $x \in \BB(\CC\bG r(\ugamma), F)$
such that
$\gamma\in \bigl(\prod_{0 \le i < r} \gamma_i\bigr)
			G_{x,r}$
(or $\gamma \in \stab_G(\ox)$, if $r = 0$).
For emphasis, we will sometimes write that the pair $(\ugamma,x)$
is an $r$-approximation to $\gamma$.
\indexme{approximation@$r$-approximation!normal}%
We say that it is a \emph{normal $r$-approximation} to
$\gamma$ if $\gamma \in \CC Gr(\ugamma)$.
\end{defn}

Note that the notion of a $0$-approximation to $\gamma$
is nearly trivial.
This is intentional, and allows us to state Lemma
\ref{lem:simult-approx} in a uniform fashion (i.e., without
separating the cases $d = 0$ and $d > 0$).

Note also that a (normal) $r$-approximation may have terms with
indices greater than $r$.  The point of this is that a
normal (say) $7$-approximation to $\gamma$ is also a normal
$3$-approximation to $\gamma$.

\begin{rk}
\label{rem:0+-approx-is-tJd}
If $(\gamma_0)$ is a normal $(0{+})$-approximation to
$\gamma$, and we put $\gamma_{> 0} := \gamma_0\inv\gamma$, then,
by Proposition \xref{J-prop:top-F-Jordan-defn} of
\cite{spice:jordan},
we have that $(\gamma_0, \gamma_{> 0})$ is a topological
Jordan decomposition of $\gamma$ modulo $Z(\bG)\conn$ (in
the sense of Definition \xref{J-defn:top-F-Jordan} of
loc.\ cit.).
Conversely, if $(\gamma\tsemi, \gamma\tunip)$ is a
topological Jordan decomposition of $\gamma$ modulo
$Z(\bG)\conn$, and $\gamma\tunip \in G_{0{+}}$\,,
then Lemma \xref{J-lem:bounded-and-top-F-unip} of loc.\ cit.\
and Proposition \ref{prop:abs-ss-good-descent}
show that there is a point
$x \in \BB(C_\bG(\gamma\tsemi), F)$
such that $((\gamma\tsemi), x)$ is a normal
$(0{+})$-approximation to $\gamma$.
\end{rk}

If $\ugamma = (\gamma_i)_{0 \le i < r'}$ is a normal $r$-approximation
to $\gamma$, then we will often write
$\gamma_{< r}$ and $\gamma_{\ge r}$
for
$\prod_{0 \le i < r} \gamma_i$ and $\gamma_{< r}\inv\gamma$,
respectively.  Of course, these depend on the choice of $\ugamma$,
not just on $\gamma$.  If we use these notations without
explicit mention of \ugamma, then the choice of normal
approximation will be irrelevant, or clear from the context.

\begin{rem}
\label{rem:approx-facts}
Suppose that $\ugamma = (\gamma_i)_{0 \le i < r}$ is an
$r$-approximation to an element $\gamma \in G$.
\begin{inc_enumerate}
\item
\label{rem:approx-facts-in-center}
For $0\leq i < r$, we have that $\gamma_i\in \ZZ Gr(\ugamma)$;
so also $\gamma_{< r} \in \ZZ G r(\ugamma)$.
\item
\label{rem:approx-facts-in-stab}
If $y \in \BB(\CC\bG r(\ugamma), F)$ and $0 \le i < r$, then
$\gamma_i \in \ZZ G r(\ugamma)$ fixes the image of $y$ in
$\rBB(\CC\bG r(\ugamma), F)$,
hence acts by a translation on $y$.
Since $\gamma_i$ is bounded modulo $Z(G)$, it actually fixes
the image of $y$ in $\rBB(\bG, F)$.
That is, $\gamma_i \in \stab_G(\ol y)$.
Since $0 \le i < r$ was arbitrary, also
$\gamma_{< r} \in \stab_G(\ol y)$.
\item
\label{rem:approx-facts-commutator}
Suppose that $y\in\BB(\CC\bG r(\ugamma),F)$,
$\gamma_{\ge r} \in G_{y, r}$\,,
and $k \in [\gamma; y, r]$.
Then
$$
[k, \gamma] \in G_{y, r}\,.
$$
Indeed, $[k, \gamma]$ is the product of the
$\gamma_{< i}$-conjugates of $[k, \gamma_i]$
(for $0 \le i < r$)
with the $\gamma_{< r}$-conjugate of $[k, \gamma_{\ge r}]$.
By Remark \ref{rem:approx-facts-in-stab},
$\gamma_{< r}$
and all $\gamma_{< i}$ lie in
$\stab_G(\ol y)$,
so it suffices to show that
$[k, \gamma_i], [k, \gamma_{\ge r}] \in G_{y, r}$
for all $0 \le i < r$.
Since $k \in G_{y, 0{+}}$\,, certainly
$[k, \gamma_{\ge r}]$ lies in $G_{y, r}$
(even $G_{y, r{+}}$).
For $0 \le i < r$, we have
$k \in (\CC G{i+}(\ugamma), G)_{y, (0{+}, r - i)}$.
Thus Lemma \ref{lem:shallow-comm}
(or Lemma \ref{lem:stab-norm}, if $i = 0$)
gives
$[\gamma_i, k] \in (\CC G{i+}(\gamma), G)_{y, (i{+}, r)}$,
and then Lemma \ref{lem:center-comm} gives
$[\gamma_i, k] \in G_{y, r}$\,.
\end{inc_enumerate}
\end{rem}

\begin{lm}
\label{lem:connected}
The $r$-approximation
$\ugamma = (\gamma_i)_{0 \le i < r'}$ to $\gamma \in G$
is normal if and only if $\gamma$ commutes with
$\gamma_i$ for $0 \le i < r$.
\end{lm}

\begin{proof}
The `only if' direction is obvious.

The `if' direction is vacuous if $r = 0$, so suppose that
$r > 0$.
Put $\bH := \bigcap_{0 \le i < r} C_\bG(\gamma_i)$.
By Definition \ref{defn:r-approx},
there is
$x \in \BB(\CC\bG r(\gamma), F) = \BB(\bH, F)$ such that
$\gamma_{\ge r} \in G_{x, r}$\,.
On the other hand, since $\gamma$ commutes with $\gamma_i$
for each $0 \le i < r$, we have
$\gamma_{\ge r} \in H$.
By Corollary \ref{cor:compatibly-filtered-tame-rank},
\bH is a compatibly filtered $F$-subgroup of \bG,
so $H \cap G_{x, r} = H_{x, r}$\,.
\end{proof}

\begin{rem}
We will see later that every $r$-approximation to $\gamma$
is conjugate to a normal $r$-approximation
(Lemma \ref{lem:normal-conj}).
Moreover, we don't have much freedom when choosing
a normal $r$-approximation (Proposition \ref{prop:unique-approx}).
\end{rem}

\begin{lemma}
\label{lem:compare-centralizers}
Suppose $\bT\subseteq \bG$ is an $F$-torus,
$\gamma\in T$,
and
$\ugamma = (\gamma_i)_{0\leq i < r'}$
is an $r$-approximation to $\gamma$
such that $\gamma_i \in T$ for $0 \le i < r$.
Then
\begin{multline*}
\set{\alpha\in \Phi(\bG,\bT)}{\ord(\alpha(\gamma)-1) \geq r} \\
=
\bigcap_{0\leq i < r} \set{\alpha\in\Phi(\bG,\bT)}{\alpha(\gamma_i)=1} \, .
\end{multline*}
\end{lemma}

\begin{proof}
It is clear that the right-hand side is contained in the left.
To prove the opposite containment,
suppose there is some root $\alpha$ such that
$\ord(\alpha(\gamma)-1) \geq r$,
but $\alpha(\gamma_i) \neq 1$ for some $0\leq i <r$.
Let $i_0$ be the minimal such $i$.
Then
$\alpha(\gamma_{i_0})\alpha(\gamma)\inv
= \bigl(\prod_{i_0 < i < r} \alpha(\gamma_i)\bigr)\inv
	\alpha(\gamma_{\ge r})\inv$,
so
$\ord(\alpha(\gamma_{i_0})\alpha(\gamma)\inv - 1)
> i_0$.
Note that, if $E/F$ is the splitting field of \bT, then
$\gamma_{\ge r}
\in G_r \cap T
	\subseteq \bG(E)_r \cap \bT(E)
	= \bT(E)_r$
by Lemmata \ref{lem:domain-field-ascent} and
\ref{lem:levi-descent},
so also
$\ord(\alpha(\gamma_{\ge r}) - 1) \ge r > i_0$.
Since
$$
\bigl(\prod_{0\leq j< i_0} \alpha(\gamma_j)\inv\bigr)
\bigl(\prod_{i_0< j< r} \alpha(\gamma_j)\inv\bigr)
\alpha(\gamma_{\geq r})\inv
=
\alpha(\gamma_{i_0})\alpha(\gamma)\inv,
$$
we have that
$\ord(\alpha(\gamma_{i_0})-1) > i_0$.
That $\gamma_{i_0} \in \GG_{i_0}^G \cup \sset 1$
then implies that $\alpha(\gamma_{i_0})=1$, which is a contradiction.
\end{proof}

\begin{cor}
\label{cor:compare-centralizers}
If $\ugamma = (\gamma_i)_{0 \le i < r'}$ is a normal
$r$-approximation to $\gamma$, then
$\CC\bG r(\ugamma) = C_\bG(\gamma_{< r})\conn$.
If, further, $\gamma$ is semisimple, then
$C_\bG(\gamma)\conn \subseteq \CC\bG r(\ugamma)$.
\end{cor}

\begin{proof}
We prove the second statement first.
By Definition \ref{defn:r-approx},
$\gamma \in \CC G r(\ugamma)$.
In particular, there is a maximal $F$-torus
$\bT \subseteq \CC\bG r(\ugamma)$ such that
$\gamma \in T$.
Then $\gamma_i \in \ZZ G r(\ugamma) \subseteq T$ for
$0 \le i < r$.
Now the result follows from Lemma \ref{lem:compare-centralizers},
together with
\cite{springer-steinberg:conj}*{\S II.4.1(b)}.

For the first statement, notice that the containment
$\CC\bG r(\ugamma) \subseteq C_\bG(\gamma_{< r})\conn$ is
clear.  The reverse containment follows from the second
statement, and the fact that \ugamma is a normal
$r$-approximation to the semisimple element $\gamma_{< r}$.
\end{proof}

\section{Good elements and commutators}
\label{sec:good-facts}

Analogues of most of the results in
this section are already known
for good elements in Lie algebras.  See
\cite{adler:thesis}*{\S 2}
and \cite{jkim-murnaghan:charexp}*{\S 2}.

In this section,
let
\begin{itemize}
\item $d\in \R_{\geq 0}$,
\item $\gamma_d \in \GG^G_d$,
\item $\bH = C_\bG(\gamma_d)$,
\item $\gamma_{>d} \in H_{d+}$\,,
and
\item $\gamma = \gamma_d\gamma_{>d}$.
\end{itemize}
Note that, by Corollary
\ref{cor:compatibly-filtered-tame-rank},
\bH is a compatibly filtered $F$-subgroup of \bG.

\begin{lemma}
\label{lem:GE2}
Suppose that $g \in \bG(F\sep)$ and
$\lsup g(\gamma_d H_{d+}) \cap \gamma_d H_{d+} \ne \emptyset$.
Then $g \in \bH(F\sep)$.
\end{lemma}

The proof
is based on a communication from Stephen DeBacker.

\begin{proof}
Suppose that $h_1, h_2 \in H_{d+}$ are such that
$\lsup g(\gamma_dh_1) = \gamma_dh_2$.
We claim that we may assume that $h_1$ and $h_2$
are semisimple.

By Lemma \ref{lem:complete-subfield}, there exists a
complete subfield $F'$ of $F$ such that
\begin{itemize}
\item $F/F'$ is unramified,
\item \bG is defined over $F'$, and
\item $\gamma_d, h_1, h_2 \in \bG(F')$.
\end{itemize}
Notice that it follows that $\bH = C_\bG(\gamma_d)$
is defined over $F'$.
By Lemma \ref{lem:domain-field-ascent}, we have
$\gamma_d \in \GG^{\bG(F')}_d$
and
$h_1, h_2 \in \bH(F')_{d+}$\,.
Since ${F'}\sep = F\sep$, we have that
$g \in \bG({F'}\sep)$.
Thus we may, and hence do, assume that $F$ is complete.
By another application of Lemma
\ref{lem:domain-field-ascent}, we see that it is harmless to
replace $F$ by finite separable extensions, so we do so as
necessary.
In particular, we may, and hence do, assume that $g \in G$.

If $\chr F = p > 0$, then let $n \in \Z_{\ge 0}$ be so large
that $h_i^{p^n}$ is semisimple for $i = 1, 2$.
Then
$$
\lsup g(\gamma_d^{p^n}h_1^{p^n})
= \bigl(\lsup g(\gamma_dh_1)\bigr)^{p^n}
= (\gamma_dh_2)^{p^n}
= \gamma_d^{p^n}h_1^{p^n}
\in \lsup g(\gamma_d^{p^n}H_{d+})
	\cap \gamma_d^{p^n}H_{d+}\,.
$$
Certainly, $\gamma_d^{p^n} \in \GG^G_d$.
An easy $\GL_n$ calculation shows that
$\bH = C_\bG(\gamma_d^{p^n})$,
so we may, and hence do, replace
$\gamma_d$ by $\gamma_d^{p^n}$
and $h_i$ by $h_i^{p^n}$ (which is semisimple)
for $i = 1, 2$.

If $\chr F = 0$, then Lemma 3.7.18 of
\cite{adler-debacker:bt-lie} shows that
$\depth_H (h_i)\semi = \depth_H(h_i) > d$
for $i = 1, 2$.
In particular, we may, and hence do, replace $h_i$ by
$(h_i)\semi$ for $i = 1, 2$.

In either case, for $i = 1, 2$,
let $\bT_i \subseteq \bH$ be a maximal $F$-torus such
that $h_i \in T_i$.
Since $\gamma_d \in Z(H\conn) \subseteq T_i$,
also $\gamma_d h_i \in T_i$ for $i = 1, 2$.
Upon enlarging $F$ if necessary, we may, and hence do,
assume that $\bT_1$ and $\bT_2$ are $F$-split.
By Lemma \ref{lem:levi-descent}, since $\bT_2$ is an $F$-Levi
subgroup of \bH, we have that
$H_{d+} \cap T_2 = (T_2)_{d+}$\,,
so $h_2 \in (T_2)_{d+}$\,.
Since
$(\lsup g\gamma_d)$ is a normal $(d{+})$-approximation
to the semisimple element $\lsup g(\gamma_dh_1)$,
we have by Corollary \ref{cor:compare-centralizers} that
$C_\bG(\lsup g(\gamma_dh_1))\conn
\subseteq C_\bG(\lsup g\gamma_d)\conn$.
In particular,
$\bT_2 \subseteq C_\bG(\lsup g\gamma_d)\conn$, so
$\lsup g\gamma_d \in Z(C_G(\lsup g\gamma_d)\conn) \subseteq T_2$.
Thus $\lsup gh_1 \in T_2$.
As above, we have that $\lsup gh_1 \in (T_2)_{d+}$\,.
Thus
$\lsup g\gamma_d
= \gamma_d(h_2(\lsup gh_1)\inv)
\in \gamma_d(T_2)_{d+}$\,.
We now discard $h_1$, $h_2$, and $\bT_1$,
put $\bT = \bT_2$,
and simply remember for later use that
$\lsup g\gamma_d \equiv \gamma_d \pmod{T_{d+}}$.

Since $\gamma_d, \lsup g\gamma_d \in T$, we have that
$\bT, \bT^g \subseteq \bH$,
hence $\bT, \bT^g \subseteq \bH\conn$.
Since both tori are maximal $F$-split in $\bH\conn$,
there is $h \in H\conn$ such that $\bT = \bT^{g h}$,
i.e., there is some element $w \in W(\bG, \bT)$
that represents $(gh)\inv$.
Then $\lsup g\gamma_d = w\inv(\gamma_d)$.

Let \bS be the subtorus of \bT generated by the images of
the $w$-fixed cocharacters of \bT.  Then any representative
for $w$ commutes with \bS, so
$w \in W(\bM, \bT)$, where $\bM = C_\bG(\bS)$.
Let $V = \bX_*(\bT) \otimes_\Z \Q$, and denote by $V^w$ and
$V_w$ the spaces of invariants and coinvariants,
respectively, of $V$ under $w$.  Since the action of $w$ on
$V$ is semisimple and $\bX_*(\bT)$ is torsion free,
$$
V = V^w \oplus V_w
= (\bX_*(\bT)^w \otimes_\Z \Q) \oplus V_w
= (\bX_*(\bS) \otimes_\Z \Q) \oplus V_w\,.
$$
Thus, if $\chi \in \bX^*(\bT/\bS) \otimes_\Z \Q$ is
$w$-fixed, then $\langle\chi, \lambda\rangle = 0$ for all
$\lambda \in V$; i.e., $\chi = 0$.

Thus, for $\alpha \in \Phi(\bM, \bT)$,
$(w - 1)\inv\alpha$ makes sense as an element of
$\bX^*(\bT/\bS) \otimes_\Z \Q$.
Now let \mo J be as in Hypothesis
\eqref{hyp:good-weight-lattice}.
By Proposition 14.2 of \cite{borel:linear}, the
multiplication map
$\mo J\semi \times Z(\mo J)\conn \to \mo J$ is an
$F$-isogeny, where $\mo J\semi$ is the derived group of \mo J.
The preimage of $\bT \subseteq \mo J$ under this map is
contained in
$(\mo J\semi \cap \bT) \times Z(\mo J)\conn$.
On the other hand, by Proposition 11.14(1) of
\cite{borel:linear}, there is a maximal torus $\bT'$ in
$\mo J\semi$ such that
$\bT' \times Z(\mo J)\conn$ maps onto \bT.
Thus $\bT' = (\mo J\semi \cap \bT)\conn$;
	%%By dimension counting:
	%%\dim \bT' + \dim Z(\bG)\conn
	%%= \dim \bT
	%%= \rk \mo J
	%%= \rk \mo J\semi + \dim Z(\bG)\conn
	%%\ge \dim (\mo J\semi \cap \bT)\conn + \dim Z(\bG)\conn.
	%%The containment
	%%$\bT' \subseteq (\mo J\semi \cap \bT)\conn$,
	%%together with connectedness,
	%%gives the desired equality.
in particular, $\bT'$ is defined over $F$, and $w$ preserves $\bT'$.
Further, there are $z \in Z(\mo J)\conn(\ol F)$ and
$\gamma_d' \in \bT'(\ol F)$ such that
$\gamma_d = \gamma_d'z$.
By restriction, we may regard $(w - 1)\inv\alpha$ as an
element of $\bX^*(\bT') \otimes_\Z \Q$.
Then, for all
$\beta\spcheck \in \Phi\spcheck(\mo J\semi, \bT')$,
$$
\langle(w - 1)\inv\alpha, \beta\spcheck\rangle
= \langle\alpha, (w - 1)\beta\spcheck\rangle
= \langle\alpha, w\beta\spcheck\rangle
	- \langle\alpha, \beta\spcheck\rangle
\in \Z;
$$
so $(w - 1)\inv\alpha \in P(\mo J\semi, \bT')$.
Then $\chi := n(w - 1)\inv\alpha$ belongs to $\bX^*(\bT')$,
where $n$ is the order of
$P(\mo J\semi, \bT')/\bX^*(\bT')$.
Denote again by $\chi$ any extension of $\chi$ to \bT.
Now recall that
$$
w\inv(\gamma_d) = \lsup g\gamma_d \equiv \gamma_d
\pmod{T_{d+}};
$$
in particular, by the definition of the filtration on $T$,
$\chi(w\inv(\gamma_d)\dotm\gamma_d\inv) \in F\cross_{d+}$\,.
Thus
$$
\alpha(\gamma_d)^n
= \alpha(\gamma_d')^n
= (w-1)\chi(\gamma_d')
= \chi(w\inv (\gamma_d') \blankdot \gamma_d'^{-1})
= \chi(w\inv(\gamma_d)\blankdot\gamma_d\inv)
\in F\cross_{d{+}}\,.
$$
By Hypothesis \eqref{hyp:good-weight-lattice},
$n$ is not divisible by $\chr\ff$,
so we have that
$\alpha(\gamma_d) \in F\cross_{d{+}}$\,.
Since $\gamma_d \in \GG^G_d \cup \sset 1$,
in fact $\alpha(\gamma_d) = 1$.
Since $\alpha \in \Phi(\bM, \bT)$ was arbitrary, $\gamma_d \in Z(M)$.
That is, $\lsup g\gamma_d = w\inv(\gamma_d) = \gamma_d$,
so $g \in \bH(F\sep)$.
\end{proof}

From now on, we assume Hypothesis \eqref{hyp:reduced}, so that
the results of \S\ref{sec:concave_tame} are
available.

The next result is the analogue of Lemma 2.3.4 of
\cite{jkim-murnaghan:charexp}.
The statement involves a concave function $f$ satisfying
some complicated conditions.  Because we need it
(or, rather, its consequence Lemma \ref{lem:good-comm})
in the proof of Lemma \ref{lem:x-depth}, we must state the
result in this generality; but, for most applications, we
are only interested in the special cases described in
Corollaries \ref{cor:iso-quotients} and \ref{cor:good-comm}
below.

\begin{lemma}
\label{lem:iso-quotients}
Let \bT be a maximal $F$-torus in \bH, and $f$
a positive, $\Gal(F\sep/F)$-invariant, concave function on
$\wtilde\Phi(\bG, \bT)$.
Put
$$
f'(\alpha)
= \begin{cases}
f(\alpha),    & \alpha \in \wtilde\Phi(\bH, \bT),     \\
f(\alpha){+}, & \alpha \not\in \wtilde\Phi(\bH, \bT).
\end{cases}
$$
Suppose that
\begin{itemize}
\item $(\bT, \bH, \bG)$ is a tame reductive $F$-sequence;
\item $x, y \in \BB(\bT, F)$;
\item $\gamma_{> d} \in H_{x, d} \cap H_{y, d{+}}$
(or $\stab_G(\ox) \cap H_{y, d{+}}$\,, if $d = 0$);
\item $f \vee f \ge f'$;
\item $f \vee f' \ge f' + \varepsilon$
for some $\varepsilon \in \R_{> 0}$;
and
\item either
	\begin{itemize}
	\item $d = 0$ and, for all discretely valued tame
extensions $L/F$, $\bH(L) \cap \stab_{\bG(L)}(\ox)$
normalizes $\lsub\bT\bG(L)_{x, f}$ and
$\lsub\bT\bG(L)_{x, f'}$\,;
	or
	\item $d > 0$,
$f \vee \ol d \ge f + d$,
and
$f' \vee \ol d \ge f' + d$,
where \ol d is the function on $\wtilde\Phi(\bG, \bT)$ which
is constant with value $d$ on $\wtilde\Phi(\bH, \bT)$,
and $\infty$ elsewhere.
	\end{itemize}
\end{itemize}
(Here, the operator $\vee$ is as in
Definition \ref{defn:concave-vee}.)
Then the map $g \mapsto [\gamma, g]$ induces
an isomorphism
\begin{equation*}
[\gamma, \blankdot ]
\colon
\lsub\bT G_{x, f:f'}
\longrightarrow
\lsub\bT G_{x, (f + d):(f' + d)}.
\end{equation*}
\end{lemma}

\begin{proof}
Let $\bS^\sharp$ be the maximal tame-modulo-$Z(\bG)\conn$
$F$-torus in \bT.
Then $\bS^\sharp$ is also a maximal
tame-modulo-$Z(\bG)\conn$ $F$-torus in \bH (and \bG).
Since $\gamma_d \in Z(H\conn)$ is tame modulo $Z(\bG)\conn$,
we have by Lemma \ref{lem:split-in-center} that
$\gamma_d \in S^\sharp$.
If $a(\gamma_d) = 1$ for all $a \in \Phi(\bG, \bS^\sharp)$,
then $\alpha(\gamma_d) = 1$ for all $\alpha \in \Phi(\bG, \bT)$,
so $\bH = \bG$ and the result is trivial.
Thus we may, and hence do, assume that there is some
$a \in \Phi(\bG, \bS^\sharp)$ such that
$a(\gamma_d) \ne 1$;
but then there is $\alpha \in \Phi(\bG, \bT)$ such that
$\alpha(\gamma_d) = a(\gamma_d) \ne 1$.
Since $\gamma_d \in \mc G^G_d$, we have
$\ord(\alpha(\gamma_d) - 1) = d$.
Let $L/F$ be a discretely valued tame extension such that
\bH and \bG are $L$-quasisplit, and the image of $\bS^\sharp$
in \tbG is $L$-split.  Then $a(\gamma_d) \in L\cross$, so
$d \in \ord(L\cross)$.

By hypothesis and Lemma \ref{lem:master-comm}, the conditions
of Corollary \ref{cor:gen-iwahori-factorization-iso} are
satisfied for $f$.
It is easy to verify that they are also satisfied for
$f + d$.
Therefore, by \S\ref{sec:unip-exp} and Corollary
\ref{cor:gen-iwahori-factorization-iso},
for $h = f$ and $h = f + d$,
\begin{equation}
\tag{$*$}
\begin{aligned}
\lsub\bT\bG(L)_{x, h:h'}
&\cong
\prod_{a \in \Phi(\bG, \bS^\sharp) \smallsetminus \Phi(\bH, \bS^\sharp)}
	U_{(a + h(a)):(a + h(a)){+}} \\
&\cong
\prod_{a \in \Phi(\bG, \bS^\sharp) \smallsetminus \Phi(\bH, \bS^\sharp)}
	(L_{\alpha(a)})_{h(a):h(a){+}}\,.
\end{aligned}
\end{equation}
Here, we have chosen (arbitrarily), for each $L$-root $a$,
an $F\sep$-root $\alpha(a)$ restricting to it;
denoted by $L_\alpha$ the fixed field in $F\sep$ of
$\stab_{\Gal(F\sep/L)} \alpha$, for $\alpha \in \Phi(\bG, \bT)$;
written $h' = f'$ if $h = f$ and $h' = f' + d$ if $h = f + d$;
and, as usual,
used, for any $L$-root $a$, the shorthand $a + k$
to denote the unique affine $L$-root $\varphi$ with gradient
$a$ such that $\varphi(x) = k$.
(Really, the product should have been over
$\Phi(\bG, \bS^{\prime\,\sharp})$, where
$\bS^{\prime\sharp}$ is the maximal
$L$-split torus in \bT; but, by Lemma
\ref{lem:torus-quotient-product}, we have that
$\bS^\sharp = Z(\bG)\conn\dotm\bS^{\prime\,\sharp}$,
so the restriction map
$\Phi(\bG, \bS^\sharp)
\to \Phi(\bG, \bS^{\prime\,\sharp})$ is a bijection.)

We will denote by $[\gamma_d, \blankdot]_L$ the map
$$
\lsub\bT\bG(L)_{x, f:f'}
\to \lsub\bT\bG(L)_{x, (f + d):(f' + d)}
$$
induced by taking commutators with $\gamma_d$.
Note that this map is well defined, and a homomorphism.
Similar notation (such as $[\gamma, \blankdot]_L$) will be
used, without further explanation, as necessary.

By \S\ref{sec:unip-exp}, since
$\gamma_d \in \bT(L)$, the map
$$
\prod_{a \in \Phi(\bG, \bS^\sharp) \smallsetminus \Phi(\bH, \bS^\sharp)}
	(L_{\alpha(a)})_{f(a):f(a){+}}
\to
\prod_{a \in \Phi(\bG, \bS^\sharp) \smallsetminus \Phi(\bH, \bS^\sharp)}
	(L_{\alpha(a)})_{(f(a) + d):(f(a) + d){+}}
$$
induced (via the isomorphisms in ($*$) for $h = f$
and $h = f + d$)
by $[\gamma_d, \blankdot]_L$ is
$$
(t_a)_a \mapsto ((a(\gamma_d) - 1)t_a)_a.
$$
Since $\gamma_d \in \GG^G_d$, we have that
$\ord(a(\gamma_d) - 1) = d$ for
$a \in \Phi(\bG, \bS^\sharp) \smallsetminus \Phi(\bH, \bS^\sharp)$,
so $[\gamma_d, \blankdot]_L$ is an isomorphism.

Our next step is to show that 
$[\gamma, \blankdot]_L$ is an isomorphism.
For the remainder of this proof, it will be convenient to do
many of our calculations in the $\ff$-algebra $\mc E$ of
endomorphisms of the finite-dimensional $\ff$-vector space
$V := \lsub\bT\bG(L)_{x, f:f'}$
or the $\ff$-vector space $\mc E'$ of homomorphisms
$V \to \lsub\bT\bG(L)_{x, (f + d):(f' + d)}$.

First, suppose $d=0$.
Then
$[\gamma,\blankdot]_L\in \mc E$.
As elements of $\mc E$,
\begin{multline*}
[\gamma,\blankdot]_L
= \Int(\gamma)-1 \\
= (\Int(\gamma_0)-1) + \Int(\gamma_0)(\Int(\gamma_{> 0}) - 1) \\
= [\gamma_0, \blankdot]_L
	+ \Int(\gamma_0)(\Int(\gamma_{> 0}) - 1).
\end{multline*}
Note that the summands commute, and the first summand is an
isomorphism.
Recall that $\gamma_{> 0}$ acts simplicially on
$\BB(\bG, F) = \rBB(\bG, F) \times V_F(Z(\bG))$.
Since
$\gamma_{> 0} \in H_{y, 0{+}} \subseteq \stab_G(y)$,
the action on the second factor is trivial.
Since $\gamma_{> 0}$ fixes \ox by assumption, it fixes
$x \in \sset\ox \times V_F(Z(\bG))$.
By Lemma \ref{lem:stab-deep}, $\gamma_{> 0} \in H_{x, 0}$\,.
By Lemma \ref{lem:unipotent}, the image of $\gamma_{> 0}$
in $(\ms H_x^L)\conn(\ff_L)$ is unipotent.
Since $\Int$ affords an algebraic representation
% Why the representation is algebraic
% (which we need in order to know that it takes unipotent
% elements to unipotent elements):
% By going up to a splitting field and using a Chevalley basis,
% we can write down the action explicitly.
of the $\ol\ff_L$-group
$(\ms H_x^L)\conn$ in the $\ol\ff_L$-vector space
in the $V \otimes_{\ff_L} \ol\ff_L$,
we have that $\Int(\gamma_{> 0}) - 1$ is a nilpotent operator
on $V$.
Therefore,
$[\gamma,\blankdot]_L$ is an isomorphism.

Now suppose $d>0$.
Recall that \ol d is the function on $\wtilde\Phi(\bG, \bT)$
which is constant with value $d$ on $\wtilde\Phi(\bH, \bT)$,
and $\infty$ elsewhere.
Since
$$
(f + d) \vee \ol d
\ge (f \vee \ol d) + d
\ge f + 2d > f' + d,
$$
we have by Lemma \ref{lem:master-comm} that
$$
[\gamma, \blankdot]_L
= [\gamma_d, \blankdot]_L
	+ \Int(\gamma_d)[\gamma_{> d}, \blankdot]_L
= [\gamma_d, \blankdot]_L
	+ [\gamma_{> d}, \blankdot]_L
$$
in $\mc E'$.
Further,
$[\gamma_d, \blankdot]_L = \Int(\gamma_d) - 1$
and
$[\gamma_{> d}, \blankdot]_L = \Int(\gamma_{> d}) - 1$
in $\mc E'$.

By Lemmata \ref{lem:domain-field-ascent} and \ref{lem:unipotent},
there exist a parabolic $L$-subgroup \bP of \bH
containing $\bS^\sharp$
and an element $u \in R\unip(\bP)(L)$
such that $\gamma_{> d} \in u\bH(L)_{x, d+}$\,.
By
Proposition 6.4.9 of
\cite{bruhat-tits:reductive-groups-1},
since $u \in \bH(L)_{x, d}$\,, there are elements
$u_a \in U_{a + d}$
for $a \in \Phi(R\unip(\bP), \bS^\sharp)$
such that $u = \prod_a u_a$.
We have that
$[\gamma_{> d}, \blankdot]_L = \sum_a [u_a, \blankdot]_L$
in $\mc E'$.
We claim that the element
$\mc T := \sum_a
	[\gamma_d, \blankdot]_L\inv \circ [u_a, \blankdot]_L$
of $\mc E$ is nilpotent.
Choose a basis of $\Phi(\bG, \bS^\sharp)$ with respect to
which the elements of $\Phi(R\unip(\bP), \bS^\sharp)$ are
positive,
and denote by $\hgt$ the associated height function.
For $i \in \Z$, denote by $V_i$
the subgroup of $\lsub\bT\bG(L)_{x, f:f'}$
generated by the images there of the affine root subgroups
$U_{b + f(b)}$ with
$b \in \Phi(\bG, \bS^\sharp) \smallsetminus \Phi(\bH, \bS^\sharp)$
and $\hgt(b) \ge i$.
Note that, for $i$ sufficiently large,
$V_i = \sset 0$ and $V_{-i} = V$.

Fix $i \in \Z$ and $a \in \Phi(R\unip(\bP), \bS^\sharp)$,
and suppose that
$b \in \Phi(\bG, \bS^\sharp) \smallsetminus \Phi(\bH, \bS^\sharp)$
satisfies $\hgt(b) \ge i$.

Notice that $a \in \Phi(\bH, \bS^\sharp)$, so $b \ne -a$.
Since the collection of restrictions of $\depth_x$
to the root groups is
a \emph{valuation de la donn\'e radicielle},
in the language of \cite{bruhat-tits:reductive-groups-1}*{\S 6.2},
D\'efinition 6.2.1(V3) of
loc.\ cit.\ implies that
$$
[u_a, U_{b + f(b)}]
\subseteq \prod U_{(m a + n b) + (m d + n f(b))}\,,
$$
the product taken over all $m, n \in \Z_{> 0}$.
If $m > 1$ or $n > 1$, then
$$
U_{(m a + n b) + (m d + n f(b))}
\subseteq \lsub\bT\bG(L)_{x, (f + d){+}}
\subseteq \lsub\bT\bG(L)_{x, f' + d}\,,
$$
so the image of $[u_a, U_{b + f(b)}]$ in
$\lsub\bT\bG(L)_{x, (f + d):(f' + d)}$ is
contained in the image there of $U_{(a + b) + (d + f(b))}$\,.
One sees as above that
the preimage under $[\gamma_d, \blankdot]_L$ of
$U_{(a + b) + (d + f(b))}$ in the latter case is precisely
$U_{(a + b) + f(b)}$\,.
Since $\hgt(a + b) > i$, we have
(whether or not $a + b \in \Phi(\bG, \bS^\sharp)$) that
$[\gamma_d, \blankdot]_L\inv \circ [u_a, \blankdot]_L$
carries the image in $V$ of $U_{b + f(b)}$ into $V_{i + 1}$.

Since
$b \in \Phi(\bG, \bS^\sharp)
	\smallsetminus \Phi(\bH, \bS^\sharp)$
was arbitrary, we have that \mc T carries
$V_i$ into $V_{i + 1}$ for all $i \in \Z$,
hence is nilpotent.
Thus $1 + \mc T$ is an isomorphism; so
$$
[\gamma, \blankdot]_L = [\gamma_d, \blankdot]_L \circ (1 + \mc T)
$$
is also.

Finally, we perform tame descent by reduction to the
complete case.
Suppose that $F'$ and $L'$ are complete subfields of $F$ and $L$,
respectively, such that
\begin{itemize}
\item $F/F'$ and $L/L'$ are unramified,
\item \bH and \bG are defined over $F'$ and quasisplit over
$L'$,
\item $\bS^\sharp$ is $L'$-split,
and
\item $x \in \BB(\bT, F')$.
\end{itemize}
We replace $L'$ by its Galois closure over $F'$,
then $F'$ by the intersection of $L'$ and $F$.
Then we have the following additional conditions.
\begin{itemize}
\item $L'/F'$ is tame and Galois,
\item $F' = L' \cap F$,
and
\item $f$ is $\Gal(L'/F')$-invariant.
% Here is why the last condition is true.
% Since $\bS^\sharp$ is $L'$-split, the action of
% $\Gal(L'\sep/L')$ on $\Phi(\bG, \bS^\sharp)$ is trivial.
% Since $f$ is constant on collections of roots restricting to
% a fixed root in $\Phi(\bG, \bS^\sharp)$, we have that $f$ is
% $\Gal(L'\sep/L')$-invariant.  Since $f$ is also
% $\Gal(F\sep/F)$-invariant, the fixed field of
% $\stab_{\Gal(F'\sep/F')} f$ includes $L' \cap F = F'$.}
\end{itemize}
By Lemma \ref{lem:tame-descent} and
Proposition \ref{prop:H1-iso},
for $h = f$ and $h = f + d$,
$$
\lsub\bT\bG(L')_{x, h:h'}^{\Gal(L'/F')}
= \lsub\bT\bG(L')_{x, h}^{\Gal(L'/F')}
	/
\lsub\bT\bG(L')_{x, h'}^{\Gal(L'/F')}
= \lsub\bT\bG(F')_{x, h:h'}\,,
$$
which is contained in $\lsub\bT G_{x, h:h'}$\,.
By Lemma \ref{lem:complete-subfield},
$\lsub\bT\bG(L)_{x, h:h'}^{\Gal(L/F)}
= \varinjlim \lsub\bT\bG(L')_{x, h:h'}^{\Gal(L'/F')}$,
the limit taken over pairs $(L', F')$ satisfying all the
above conditions; so
$\lsub\bT\bG(L)_{x, h:h'}^{\Gal(L/F)}
\subseteq \lsub\bT G_{x, h:h'}$\,.
The reverse containment, hence equality, is easy.
Since $[\gamma, \blankdot]_L$ is $\Gal(L/F)$-equivariant,
$[\gamma, \blankdot]_F$ is an isomorphism, as desired.
\end{proof}

\begin{cor}
\label{cor:iso-quotients}
Suppose that $\gamma_{> d} \in H_{x, d} \cap H_{d{+}}$\,.
For $r \in \tR_{> 0}$,
the map $g \mapsto [\gamma, g]$ induces an isomorphism
$$
[\gamma, \blankdot]
: (H, G)_{x, (r, r):(r, r{+})}
\to (H, G)_{x, (r + d, r + d):(r + d, (r + d){+})}.
$$
\end{cor}

\begin{proof}
Note that the domain and codomain of $[\gamma, \blankdot]$
are both trivial if $r = s{+}$ for some
$s \in \R$, or if $r = \infty$; so we assume that
$r \in \R$.

Let $y$ be any point of $\BB(\bH, F)$ such that
$\gamma_{> d} \in H_{y, d{+}}$\,,
and let \bT be a maximal $F$-torus in \bH such that
$(\bT, \bH)$ is a tame reductive $F$-sequence and
$x, y \in \BB(\bT, F)$.
(To see that such a torus \bT exists, let \bS be any maximal
$F$-split torus such that $x, y \in \AA(\bS, F)$.
By Lemma \ref{lem:ratl-maxl-torus}, there is a maximal
$F\tame$-split torus $\bS^\sharp$, defined over $F$,
containing \bS.
Then we may take \bT to be $C_\bG(\bS^\sharp)$.)
Let $f$ be the constant function on $\wtilde\Phi(\bG, \bT)$
with value $r$.
Note that $\lsub\bT G_{x, f} = G_{x, r} = (H, G)_{x, (r, r)}$
and $\lsub\bT G_{x, f'} = (H, G)_{x, (r, r{+})}$\,,
where $f'$ is as in Lemma \ref{lem:iso-quotients},
and similarly with $f + d$ and $f' + d$ in place of $f$ and
$f'$, respectively.
Thus it suffices to show that the hypotheses of Lemma
\ref{lem:iso-quotients} are satisfied.

We have that $f \vee f$ is the constant function with value
$2 r \ge r{+}$,
so $f \vee f \ge f'$;
and
$f \vee f'$ is constant with value $2 r \ge r$ on
$\wtilde\Phi(\bH, \bT)$, and takes the value
$2 r{+}$
elsewhere,
so $f \vee f' \ge f' + r$.
If $d > 0$, and we define \ol d as in Lemma
\ref{lem:iso-quotients}, then a similar calculation shows
that
$f \vee \ol d = f + \ol d \ge f + d$
and
$f' \vee \ol d = f' + \ol d \ge f' + d$.
If $d = 0$, then Lemma \ref{lem:stab-norm} gives that
$\bH(L) \cap \stab_{\bG(L)}(\ox)$ normalizes
both
$\lsub\bT\bG(L)_{x, f} = \bG(L)_{x, r}$
and
$\lsub\bT\bG(L)_{x, f'} = (\bH, \bG)(L)_{x, (r, r{+})}$.
\end{proof}

The next result is the analogue of Corollary 2.3.5 of
\cite{jkim-murnaghan:charexp}.

\begin{lemma}
\label{lem:good-comm}
With the notation and bulleted hypotheses of Lemma \ref{lem:iso-quotients},
suppose further that
\begin{itemize}
\item $f \vee f \ge f + (0{+})$;
and
\item either
	\begin{itemize}
	\item $d = 0$ and, for all discretely valued tame
extensions $L/F$ and all $t \in \tR_{\ge 0}$,
$\bH(L) \cap \stab_{\bG(L)}(\ox)$ normalizes
$\lsub\bT\bG(L)_{x, f + t}$ and $\lsub\bT\bG(L)_{x, f' + t}$;
or
	\item $d > 0$.
	\end{itemize}
\end{itemize}
Then the $\lsub\bT G_{x, f}$-orbit of
$\gamma\dotm\lsub\bT H_{x, f + d}$ is
$\gamma\dotm\lsub\bT G_{x, f + d}$\,.
\end{lemma}

\begin{proof}
We have
$\gamma\dotm\lsub\bT H_{x, f + d} \subseteq H_{x, d}$
(or $H \cap \stab_G(\ox)$, if $d = 0$).
Since the commutator of $H_{x, d}$ (or $H \cap \stab_G(\ox)$,
if $d = 0$) and
$\lsub\bT G_{x, f}$ lies in $\lsub\bT G_{x, f + d}$\,,
we have that the $\lsub\bT G_{x, f}$-orbit of
$\gamma\dotm\lsub\bT H_{x, f + d}$ is contained in
$\gamma\dotm\lsub\bT G_{x, f + d}$\,.
Thus
we need only show that every element of
$\gamma\dotm\lsub\bT G_{x, f + d}$ is
$\lsub\bT G_{x, f}$-conjugate to an element of
$\gamma\dotm\lsub\bT H_{x, f + d}$\,.
By Lemma \ref{lem:complete-subfield},
given any such element, we may find a complete subfield
of $F$ over which that element is defined, and of which $F$
is an unramified extension.
Upon replacing $F$ by this complete subfield,
we may, and hence do, assume that $F$ is complete.

For $t \in \tR_{\ge 0}$ and $\alpha \in \wtilde\Phi(\bG, \bT)$, put
$$
f'_t(\alpha)
= \begin{cases}
f(\alpha),     & \alpha \in \wtilde\Phi(\bH, \bT),     \\
f(\alpha) + t, & \alpha \not\in \wtilde\Phi(\bH, \bT).
\end{cases}
$$
(Thus, in the notation of Lemma \ref{lem:iso-quotients},
$f' = f'_{0{+}}$.)
By Lemma \ref{lem:more-vGvr-facts},
$$
\bigcap_{t \ge 0} \lsub\bT G_{x, f'_t + d}
	= \lsub\bT H_{x, f + d}\,.
$$
Thus, it suffices to prove that, for any $t \in \R_{\ge 0}$,
any element of $\gamma\dotm\lsub\bT G_{x, f'_t + d}$ is conjugate
by $\lsub\bT G_{x, f + t}$ to an element of
$\gamma\dotm\lsub\bT G_{x, f'_{t{+}} + d}$\,.

Choose $\delta \in \lsub\bT G_{x, f'_t + d}$\,.
By Proposition \ref{prop:heres-a-gp},
% Here's why we can use Proposition \ref{prop:heres-a-gp}:
% $$
% (f + t + d) \vee (f'_{t{+}} + d) \ge (f \vee f) + (t + d) \ge f + t + d,
% $$
we have
$$
\lsub\bT G_{x, f'_t + d}
= \lsub\bT G_{x, f + t + d}
	\dotm\lsub\bT G_{x, f'_{t{+}} + d}\,.
$$
In particular, $\lsub\bT G_{x, f'_{t{+}} + d}\dotm\delta\inv$ contains an
element of $\lsub\bT G_{x, f + t + d}$\,, say $\delta'$.
Note that $f + t$ satisfies the conditions of Lemma
\ref{lem:iso-quotients},
so there is an element
$g \in \lsub\bT G_{x, f + t}$ such that
$$
[\gamma\inv, g] \in \lsub\bT G_{x, f' + (t + d)}\dotm\delta'
	\subseteq \lsub\bT G_{x, f'_{t{+}} + d}\dotm\delta'.
$$
In particular,
$[\gamma\inv, g] \in \lsub\bT G_{x, f'_{t{+}} + d}\dotm\delta\inv$.
Since
$$
(f'_t + d) \vee (f + t)
\ge (f \vee f) + (t + d)
\ge f + (t + d){+}
\ge f'_{t{+}} + d
$$
by our hypothesis on $f$,
Lemma \ref{lem:master-comm} gives
$[\delta\inv, g] \in \lsub\bT G_{x, f'_{t{+}} + d}$\,.
Thus
$$
\lsup g(\gamma\delta)
= \gamma\dotm{[\gamma\inv, g]}\delta\dotm{[\delta\inv, g]}
\in \gamma\dotm\lsub\bT G_{x, f'_{t{+}} + d}\,.\qedhere
$$
\end{proof}

\begin{cor}
\label{cor:good-comm}
Suppose that $\gamma_{> d} \in H_{x, d} \cap H_{d{+}}$
(or $\stab_G(\ox) \cap H_{d{+}}$\,, if $d = 0$).
For any $r \in \tR_{> 0}$, the $G_{x, r}$-orbit of
$\gamma H_{x, r + d}$ is $\gamma G_{x, r + d}$\,.
\end{cor}

\begin{proof}
This corollary follows from Corollary \ref{cor:iso-quotients} in
exactly the same way as Lemma \ref{lem:good-comm}
follows from Lemma \ref{lem:iso-quotients}.
\end{proof}

The next result is the analogue of Lemma 2.3.3 of
\cite{jkim-murnaghan:charexp}.  Corollary 4.4.3 of
\cite{debacker:nilp} is a similar result.

\begin{lemma}
\label{lem:x-depth}
Let $x\in\BB(\bG,F)$.
If $d = 0$ and $\gamma \in \stab_G(\ox)$, or
$d>0$ and $\gamma\in G_{x, d}$\,, then
$x\in \BB(\bH,F)$.
\end{lemma}

\begin{proof}
If $d = 0$, then, by Proposition
\ref{prop:abs-ss-good-descent},
we have
$\ox \in \BB(C_\tbG(\ogamma_0), F)
	= \BB(C_\tbG(\ogamma_0)\conn, F)$.
By Proposition 9.6 of \cite{borel:linear},
$$
C_\tbG(\ogamma_0)\conn
= C_\bG(\gamma_0)\conn/Z(\bG)\conn
= \bH\conn/Z(\bG)\conn.
$$
Thus,
$\ox \in \BB(\bH\conn/Z(\bG)\conn, F)
	\subseteq \BB\red(\bG, F)$.
Since the preimage in $\BB(\bG, F)$ of
$\BB(\bH\conn/Z(\bG)\conn, F)$ is
$\BB(\bH\conn, F)$, we have that
$x \in \BB(\bH\conn, F) = \BB(\bH, F)$.

Now suppose that $d > 0$.
By Hypothesis \eqref{hyp:conn-cent}, $\bH\conn$ is a Levi
subgroup of \bG.
Let $L/F$ be a discretely valued separable extension such that
\bH, hence also \bG, is $L$-split.
Then $\bH\conn$ is an $L$-Levi subgroup of \bG.
Moreover, since $\gamma_d$ lies in any $L$-split maximal torus
of \bH, we have that $d = \depth(\gamma_d) \in \ord(L\cross)$.
Let \bP be a parabolic $L$-subgroup of \bG with Levi
component $\bH\conn$, and $\bP'$ its opposite parabolic
(with respect to any maximal torus in $\bH\conn$).
Let \bU and $\bU'$ be the unipotent radicals of \bP and $\bP'$,
respectively.
Put
$$
\BB(\gamma) = \set{y \in \BB(\bG, L)}{\gamma \in \bG(L)_{y, d}}.
$$
If $\BB(\gamma) \subseteq \BB(\bH, L)$, then
Lemma \ref{lem:levi-building-descent} gives
$x \in \BB(\gamma) \cap \BB(\bG, F) \subseteq \BB(\bH, F)$,
and we are done.
Otherwise, we claim that that $\BB(\gamma)$ is convex chamber closed
in the sense of \cite{moy:displacement}*{\S 3.1}.

Indeed, it is clear that $\BB(\gamma)$ is closed, convex,
and the union of closures of facets.
By Lemma \ref{lem:domain-field-ascent},
$H_{d+} \subseteq \bH(L)_{d+}$\,, so
$\gamma \in \gamma_d\bH(L)_{d+}$\,.
Thus there is some
chamber $C_0 \subseteq \BB(\bH, L)$
so that $\gamma \in \gamma_d\bH(L)_{x, d+}$ for all
$x \in C_0$.
By Lemma \ref{lem:depth-in-center},
$\gamma_d \in \bH(L)_{x, d}$ for all
$x \in C_0$,
so $C_0 \subseteq \BB(\gamma)$.
If $J$ is a facet in $\BB(\gamma)$, then there is an
apartment \fAA containing both $C_0$ and $J$, hence the convex
hull of $C_0$ and $J$.
Since \fAA is a Euclidean space, we know that the convex hull
of $C_0$ and $J$, hence $\BB(\gamma)$, contains a chamber $C'$
of $\BB(\bG,L)$
such that
$J \subseteq {\ol C}\,'$.
That is, $\BB(\gamma)$ is a union of closures of chambers.

Thus
there is a chamber $C \subseteq \BB(\gamma)$ such that
$C \cap \BB(\bH, L) = \emptyset$ and
$\ol C \cap \BB(\bH, L)$ contains a facet
$J$ of $C$ of codimension $1$.  Fix $y \in J$.
By Lemma 2.4.1 of
\cite{adler-debacker:bt-lie}, there is
$u \in \bU(L) \cap \bG(L)_y$ such that
$u C \subseteq \BB(\bH, L)$.
Note that the closure of $u C$ includes $u y = y$.
Fix $z \in C$.
By Theorem 4.2 of \cite{moy-prasad:jacquet}, 
since $u z \in u C \subseteq \BB(\bH, L)$,
the preimage of $\bG(L)_{u z, d}$ under the multiplication map
$$
\bU'(L) \times \bH\conn(L) \times \bU(L) \to \bG(L)
$$
is
$$
(\bU'(L) \cap \bG(L)_{u z, d})
\times
\bH(L)_{u z, d}
\times
(\bU(L) \cap \bG(L)_{u z, d}).
$$
Since $z \in C \subseteq \BB(\gamma)$,
we have that $\gamma \in \bG(L)_{z, d}$\,, so
$\gamma\cdot[\gamma\inv, u] = \lsup u\gamma
\in \bG(L)_{u z, d}$\,.
Since $\gamma \in \bH\conn(L)$ and
$[\gamma\inv, u] \in \bU(L)$,
in fact $\gamma \in \bH(L)_{u z, d}$
and
$[\gamma\inv, u] \in \bU(L) \cap \bG(L)_{u z, d}$\,.
By Lemma \ref{lem:unipotent},
$[\gamma\inv, u] \in \bU(L) \cap \bG(L)_{w, d{+}}$
for some point $w$ near $u z$.
In particular, we may choose $w$ to lie in $uC$.
Since $d \in \ord(L\cross)$, we have that
$\bG(L)_{w, d{+}} = \bG(L)_{u z, d{+}}$, so
$[\gamma\inv, u] \in \bU(L) \cap \bG(L)_{u z, d{+}}$\,.

That is,
$\lsup u\gamma \in \gamma(\bU(L) \cap \bG(L)_{u z, d{+}})
	= \gamma\dotm\lsub\bT\bG(L)_{u z, f + d}$\,,
where \bT is an $L$-torus in \bH of maximal $F\tame$-split rank,
and
$f$ is the function on $\wtilde\Phi(\bG, \bT)$ which
is $\varepsilon$ on $\Phi(\bU, \bT)$ and $\infty$ elsewhere
(for some suitably small $\varepsilon \in \R_{> 0}$).
By Lemma \ref{lem:good-comm}, there is
$k \in \lsub\bT\bG(L)_{u z, f}$ such that
$\lsup{k u}\gamma \in \gamma\dotm\lsub\bT\bH(L)_{u z, f + d}
	= \sset{\gamma}
	\subseteq \gamma_d\bH(L)_{d{+}}$\,.
Since $k u \in \bU(L)$, there is a
one-parameter subgroup $\lambda \in \bX_*(Z(\bH\conn))$ such that
$\lim_{t \to \infty} \lsup{\lambda(t)}(k u)
= 1$.
On the other hand,
by Lemma \ref{lem:GE2}, $k u \in \bH(L)$; so
$\lim_{t \to \infty} \lsup{\lambda(t)}(k u)
\equiv k u \pmod{\bH\conn(L)}$.
Thus $k u \in \bH\conn(L) \cap \bU(L) = \sset 1$,
so $u = k\inv \in \lsub\bT \bG(L)_{u z, f} \subseteq \stab_G(u z)$.
In particular,
$z = u\inv(u z) = u z \in C \cap u C
\subseteq C \cap \BB(\bH, L) = \emptyset$, which is
a contradiction.
\end{proof}

\section{Existence and uniqueness of normal approximations}
\label{sec:approx-facts}

We now discard the notation $d$, $\gamma$, $\gamma_d$, and
\bH of the preceding section.
Fix $\gamma \in G$ and $r \in \tR_{\ge 0}$.

In this section, we show that normal approximations
often exist.  When they exist, they are essentially unique.

\begin{lm}
\label{lem:simult-approx}
Suppose that
\begin{itemize}
\item $(\bG', \bG)$ is a tame reductive $F$-sequence;
\item $d \in \tR_{\ge 0}$;
\item $d > 0$ or
$\tG_{0+}$ is contained in the image of $G_{0+}$\,;
\item $\gamma$ is bounded modulo $Z(G)$;
and
\item $\gamma$ has a normal $d$-approximation
$(\gamma_i)_{0 \le i < d}$ with
$\gamma_i \in G'$ for $0 \le i < d$.
\end{itemize}
If there is a maximal tame-modulo-center (of \bG)
$F$-torus $\bS \subseteq \CC{\bG'}d(\gamma)$,
satisfying property \bGd G, such that $\gamma \in S$,
then
$\gamma_i \in S$ for $0 \le i < d$,
and there are elements $\gamma_i \in S$ for $d \le i < r$
such that $(\gamma_i)_{0 \le i < r}$ is a normal
$r$-approximation to $\gamma$.
\end{lm}

The hypothesis that $\gamma$ is bounded modulo
$Z(G)$ is redundant if $d > 0$.
If $d = 0$, the existence of a normal $d$-approximation is
vacuous.

The point of this result is to show that normal
approximations exist in abundance.  Since it is not much
additional work, we show that we may find normal
approximations in $G$ and $G'$ simultaneously.

\begin{proof}
If the result holds (with $d$ fixed)
for all $r'$ sufficiently close
to $r$, then taking the union of successively larger
$r'$-approximations gives the result for $r$.
Thus the set of all $r$ for which the result holds is
order-closed.
Further, if the result holds for the pair $(d, r)$, then it
holds for the pair $(d, r{+})$ if and only if it holds for
the pair $(r, r{+})$.
Finally, since
a normal $(r{+})$-approximation is precisely
a normal $(r + \varepsilon)$-approximation for
$\varepsilon \in \R_{> 0}$ sufficiently small,
the result holds for the pair $(d, r{+})$ if and only if it
holds for all pairs $(d, r + \varepsilon)$ with
$\varepsilon \in \R_{> 0}$ sufficiently small.
Thus, if we can show the result for all pairs
$(d, d{+})$, then we will have shown that (with $d$ fixed)
the set of all $r$ for which the result holds is also
order-open, hence consists of all $r \in \tR$.
That is, we may, and hence do, suppose throughout that
$r = d{+}$.

By Corollary \ref{cor:torus-quotient-product},
$\bT := C_\bG(\bS)$ is an $E$-Levi $F$-subgroup of \bG,
where $E/F$ is the splitting field of $\bS/Z(\bG)\conn$.
It is also a maximal torus in \bG.

Suppose $d = 0$.
By Proposition \xref{J-prop:tJd=cpct} of \cite{spice:jordan},
$\ogamma$ has a topological Jordan decomposition
$\ogamma = \ogamma\tsemi\ogamma\tunip$ (in the sense of
Definition \xref{J-defn:top-F-Jordan} of loc.\ cit.).
In particular, $\ogamma\tsemi$ belongs to a torus containing
$\ogamma = \ogamma\semi$,
hence, by Lemma \xref{J-lem:top-F-Jordan-central} of
loc.\ cit.,
to every maximal such torus
--- in particular, to the image in \tbG of
\bT.
By Proposition \xref{J-prop:tame-tunip} of loc.\ cit.,
we have that $\ogamma\tunip \in \tG_{0+}$\,;
so, by hypothesis, $\ogamma\tunip$ has a preimage
$\gamma_{> 0} \in G_{0+}$\,.
Put $\gamma_0 = \gamma\gamma_{> 0}\inv$.
Then $(\gamma_0)\spdash = \ogamma\tsemi$ belongs to the
image in \tG of $T$,
so $\gamma_0 \in Z(G)T = T$,
hence also $\gamma_{> 0} \in T$.
By Corollary \xref{J-cor:abs-F-ss-tame} of
loc.\ cit., $(\gamma_0)\spdash$ is tame, i.e.,
$\gamma_0$ is tame modulo center in \bG.
By Lemma \ref{lem:still-split-in-levi},
$\gamma_0$ remains tame modulo $Z(\bG)$ in \bT,
hence belongs to $S$, the
(set of $F$-rational points of the)
maximal tame-modulo-$Z(\bG)$ torus in \bT.
By Lemmata \ref{lem:domain-field-ascent},
\ref{lem:levi-descent},
and \ref{lem:field-descent},
we have
$\gamma_{> 0} \in T_{0{+}} \subseteq G_{x, 0{+}}$
for any $x \in \BB(\bT, F) \subseteq \BB(C_\bG(\gamma_0), F)$.
Thus, $(\gamma_0)$ is a normal $(0{+})$-approximation to
$\gamma$.

Now suppose $d > 0$.
By Definitions \ref{defn:good} and \ref{defn:r-approx},
$\gamma_i$ belongs to a tame-modulo-center
$F$-torus in $\CC\bG d(\gamma)$ for $0 \le i < d$.
By Remark \ref{rem:approx-facts-in-center},
$\gamma_i \in \ZZ G d(\gamma)$ for $0 \le i < d$.
Since \bS is a maximal tame-modulo-center $F$-torus in
$\CC\bG d(\gamma)$,
we have by Lemma \ref{lem:split-in-center} that
$\gamma_i \in S$ for $0 \le i < d$.
Thus $\gamma_{\ge d} \in S$,
and we may, and hence do, replace $\bG'$ by
$\CC{\bG'}d(\gamma)$.
Since $\gamma_{\ge d} \in G_d$\,, we have by Lemma
\ref{lem:domain-field-ascent} that
$\gamma_{\ge d} \in \bG(E)_d$\,,
hence by Lemma \ref{lem:levi-descent} that
$\gamma_{\ge d} \in \bT(E)_d$\,.
By Lemma 2.4 of \cite{rapoport:T1-is-T0}, we have
$\bS(E)_0 \subseteq \bT(E)_0$\,.
By the definition of the filtration on $\bS(E)$, we have
that
$\bS(E) \cap \bT(E)_d = \bS(E)_0 \cap \bT(E)_d
= \bS(E)_d$\,;
in particular, $\gamma_{\ge d} \in \bS(E)_d$\,.
By Lemma \ref{lem:field-descent},
$\gamma_{\ge d} \in S_d$\,.
By property \bGd G, there is an element
$\gamma_d
\in \gamma_{\ge d}S_{d+}
	\cap (\GG^G_d \cup \sset 1)$.
Put $\ugamma = (\gamma_i)_{0 \le i \le d}$
and choose
$x \in \BB(C_{\bG'}(\bS), F)
	\subseteq \BB(\CC{\bG'}{d+}(\ugamma), F)$.
Then
$\bigl(\prod_{0 \le i \le d} \gamma_i\bigr)\inv\gamma \in S_{d+}$\,,
and reasoning as above shows that
$S_{d+} \subseteq G'_{x, d+} \subseteq G_{x, d+}$\,.
That is,
\ugamma is a $(d{+})$-approximation to $\gamma$.
Since $\gamma_i \in S$ for $0 \le i \le d$,
we have that $\gamma \in S \subseteq \CC G{d+}(\ugamma)$,
so \ugamma is a
\emph{normal} $(d{+})$-approximation to $\gamma$.
\end{proof}

In Lemma \ref{lem:simult-approx}, we showed that, if
$(\bG', \bG)$ is a tame reductive $F$-sequence, then
some normal approximation in $G$ is also a normal
approximation in $G'$.  The next result shows that,
if $(\bG', \bG)$ is a tame \emph{Levi} $F$-sequence, then
any normal approximation in $G$ is already also a
normal approximation in $G'$.

Recall that $r \in \tR_{\ge 0}$.

\begin{lm}
\label{lem:G-approx-is-G'-approx}
Suppose that
\begin{itemize}
\item $(\bG', \bG)$ is a tame Levi $F$-sequence,
\item $\gamma \in G'$,
\item $\ugamma = (\gamma_i)_{0 \le i < r}$ is a normal
$r$-approximation to $\gamma$ in $G$,
and
\item $r \le 0{+}$, $\gamma$ is semisimple,
or $\gamma_{< r} \in G'$.
\end{itemize}
Then \ugamma is a normal $r$-approximation to $\gamma$
in $G'$, and
$$
\set{x \in \BB(\CC\bG r(\gamma), F)}
	{\gamma_{\ge r} \in G_{x, r}}
\cap \BB(\bG', F)
= \set{x \in \BB(\CC{\bG'}r(\gamma), F)}
	{\gamma_{\ge r} \in G'_{x, r}}.
$$
\end{lm}

\begin{proof}
We may, and hence do, assume that $r > 0$
(since otherwise the statement is vacuous).
To show that \ugamma is a normal $r$-approximation to
$\gamma$ in $G'$, we need to show two things:
\begin{enumerate}
\item that it is a good sequence in $G'$
(i.e., that there exists a tame-modulo-center $F$-torus
in $\bG'$ containing all the $\gamma_i$);
and then
\item that $\gamma_{\ge r} \in \CC{\bG'}r(\ugamma)_r$\,.
\end{enumerate}

Suppose that $\gamma_{< r} \in G'$.
Then also $\gamma_{\ge r} = \gamma_{< r}\inv\gamma \in G'$.
By Corollary \ref{cor:compare-centralizers},
$\ZZ G r(\gamma) = Z(C_G(\gamma_{< r})\conn) \subseteq G'$.
By Remark \ref{rem:approx-facts-in-center},
$\gamma_i \in \ZZ G r(\gamma) \subseteq G'$
for all $0 \le i < r$.
Since there is a tame $F$-torus in \bG containing
all the $\gamma_i$,
we have by Lemma \ref{lem:still-split-in-levi}
that there is a tame $F$-torus in $\bG'$ containing all the
$\gamma_i$.
Thus, \ugamma is a good sequence in $G'$.
Put $\bH = \CC\bG r(\ugamma)$
and $\bH' = \CC{\bG'}r(\ugamma)$.
By Lemma \ref{lem:funny-centralizer-descends},
$\bH' = \bH \cap \bG'$.
In particular, $\bH'$ is an $F\tame$-Levi $F$-subgroup of
\bH; say $L/F$ is a tame finite extension such that
$\bH'$ is an $L$-Levi $F$-subgroup of \bH.
Since \ugamma is a normal $r$-approximation to $\gamma$ in
$G$, we have that
$\gamma_{\ge r} \in H_r \subseteq \bH(L)_r$\,;
so, by Lemmata \ref{lem:levi-descent} and
\ref{lem:domain-field-ascent},
we have that
$\gamma_{\ge r} \in H' \cap \bH(L)_r
	= H' \cap \bH'(L)_r
	= H'_r$\,.
Thus, \ugamma is a normal $r$-approximation to $\gamma$ in
$G'$.
Now another application of
Lemmata \ref{lem:levi-descent} and
\ref{lem:domain-field-ascent},
together with
Lemma \ref{lem:levi-building-descent},
gives the desired equality of subsets of $\BB(\bG',F)$.

Thus, it remains only to show that
$\gamma_{< r} \in G'$.
Remember that we have assumed that
$r \le 0{+}$, $\gamma$ is semisimple, \emph{or}
$\gamma_{< r} \in G'$.
Thus, it suffices to show that either of the first two conditions
implies the third.

If $\gamma$ is semisimple, then, by Corollary
\ref{cor:compare-centralizers}, we have
$\ZZ G r(\gamma) \subseteq Z(C_G(\gamma)\conn) \subseteq G'$.
By Remark \ref{rem:approx-facts-in-center},
$\gamma_{< r} \in \ZZ G r(\gamma)$.

If $r = 0{+}$, then, by Remark \ref{rem:0+-approx-is-tJd},
$(\gamma_0, \gamma_{> 0})$ is a topological $F$-Jordan
decomposition of $\gamma$ in $G$, modulo $Z(\bG)$.
Denote by \ogamma and $\ogamma_0$ the images of $\gamma$ and
$\gamma_0$ in \tG.
By Lemma \xref{J-lem:top-F-Jordan-central} of
\cite{spice:jordan},
$\ogamma_0 \in Z(C_\tG(\ogamma))$.
Since $\gamma \in \CC G r(\gamma) = C_G(\gamma_0)\conn$,
we have that $\gamma_0 \in C_G(\gamma)\conn$,
so actually $\ogamma_0 \in Z(C_\tG(\ogamma)\conn)$.
Finally, by Proposition 9.6 of \cite{borel:linear},
we have that $C_\bG(\gamma)\conn$ surjects onto
$C_\tbG(\ogamma)\conn$, hence that
$Z(C_\bG(\gamma)\conn)$ surjects onto
$Z(C_\tbG(\ogamma)\conn)$.
Indeed, it is easy to see that $Z(C_\bG(\gamma)\conn)$
is the full preimage in \bG of $Z(C_\tbG(\ogamma)\conn)$.
Thus,
$\gamma_{< r} = \gamma_0 \in Z(C_G(\gamma)\conn) \subseteq G'$.
\end{proof}

The next two results concern the extent to which a normal
approximation to $\gamma$ is determined by $\gamma$.
If one is concerned only with semisimple elements (for
example, for computing the values of supercuspidal
characters of reductive $p$-adic groups, as in
\cite{adler-spice:explicit-chars}), the argument below can be
replaced with a root-value calculation in the spirit of
Lemma \ref{lem:compare-centralizers} and
Corollary \ref{cor:compare-centralizers}.
We prove the results for general elements $\gamma$ since the
argument is not much longer.

\begin{lemma}
\label{lem:normal-approx-init-1s}
If $\ugamma = (\gamma_i)_{0\leq i < r}$
is a normal $r$-approximation to $\gamma$ and $\gamma_0 = 1$, then
$\gamma_i = 1$ for all $0 \le i < \depth(\gamma)$.
\end{lemma}

\begin{proof}
Suppose that $j < \depth(\gamma)$, and that
$\gamma_i = 1$ for $0 \le i < j$.
Then, by Definition \ref{defn:r-approx}, there is
$x \in \BB(\CC\bG{j+}(\gamma), F)$ such that
$\gamma \in \gamma_j G_{x, j+} \cap \CC G{j+}(\ugamma)
	= \gamma_j\CC G{j+}(\ugamma)_{x, j+}$\,.
By Lemma \ref{lem:unipotent},
$\depth(\gamma_j) > j$.
Since
$\gamma_j \in \GG^G_j \cup \sset 1$, this means that
$\gamma_j = 1$.
That is, $\gamma_i = 1$
for all $0 \le i < \depth(\gamma)$.
\end{proof}

\begin{pn}
\label{prop:unique-approx}
If $\ugamma = (\gamma_i)_{0\leq i < r}$
and $\ugamma' = (\gamma_i')_{0\leq i < r}$
are normal $r$-approximations to $\gamma$, then,
for all $0\leq i \leq r$,
\begin{equation}
\tag{$\textbf{Cent}_i$}
\CC\bG{i}(\ugamma) = \CC\bG{i}(\ugamma');
\end{equation}
and, for all $0 \le i < r$,
\begin{equation}
\tag{$\textbf{Cong}_i$}
\gamma_i' \equiv \gamma_i \bmod
	\bigl(\ZZ G{i}(\ugamma)\cap \CC G{i}(\ugamma)_i \bigr)
	\bigl(\ZZ G{i+}(\ugamma)\cap \CC G{i+}(\ugamma)_{i+} \bigr).
\end{equation}
Further, 
\begin{equation}
\tag{$\textbf{Cong}_j'$}
\gamma_j' \equiv \gamma_j \bmod
	\bigl(\ZZ G{j+}(\ugamma) \cap \CC G{j+}(\ugamma)_{j+}\bigr)
\end{equation}
whenever $j < r$ is such that
$\gamma_i = \gamma_i'$ for $0 \le i < j$.
\end{pn}

Note that, for any index $i$,
($\textbf{Cong}_i'$) implies
($\textbf{Cong}_i$).

\begin{proof}
Note that ($\textbf{Cent}_0$) is true.
For $0 \le i \le j < r$, we have that
$$
\ZZ G i(\ugamma) \subseteq \ZZ G j(\ugamma)
\subseteq
\CC G j(\ugamma) \subseteq \CC G i(\ugamma).
$$
Suppose that $j \in \R_{\ge 0}$
is such that
$j < r$,
($\textbf{Cent}_i$) holds
for all $0 \le i \le j$,
and
($\textbf{Cong}_i$) holds
for all $0 \le i < j$.
In particular, there are elements
$z_i \in \ZZ G{i+}(\ugamma) \subseteq \ZZ G j(\ugamma)$
such that
$\gamma_i' = z_i\gamma_i$ for $0 \le i < j$.
If $\gamma_i = \gamma_i'$ for $0 \le i < j$,
then $z_i = 1$ for $0 \le i < j$.
Put
\begin{gather*}
\gamma_{\ge j}
= \bigl(\prod_{0 \le i < j} \gamma_i\bigr)\inv\gamma,
\quad
\gamma_{> j}
= \gamma_j\inv\gamma_{\ge j}, \\
\gamma_{\ge j}'
= \bigl(\prod_{0 \le i < j} \gamma_i'\bigr)\inv\gamma,
\quad
\gamma_{> j}'
= \gamma_j'^{-1}\gamma_{\ge j}'.
\end{gather*}
Then
$\gamma_{\ge j} = z\gamma'_{\ge j}$,
where
$z = \prod_{0 \le i < j} z_i \in \ZZ G j(\ugamma)$.
If $\gamma_i = \gamma'_i$ for $0 \le i < j$, then $z = 1$.

Since $\CC\bG{j+}(\ugamma) = \CC{\CC\bG j(\ugamma)}{j{+}}(\ugamma)$,
and similarly for $\ugamma'$, we may, and hence do, replace
\bG by $\CC\bG j(\ugamma) = \CC\bG j(\ugamma')$.

Since a normal $r$-approximation to $\gamma$ is \emph{a fortiori} a
normal $(j{+})$-approximation to it,
there are
$x \in \BB(\CC\bG{j+}(\ugamma), F)$
and
$x' \in \BB(\CC\bG{j+}(\ugamma'), F)$
such that $\depth_x(\gamma_{> j}) > j$ and
$\depth_{x'}(\gamma'_{> j}) > j$.
Then
$$
\gamma_{\ge j} \in \Upsilon \cap z\Upsilon',
$$
where
$$
\Upsilon := \gamma_j\CC G{j+}(\ugamma)_{x, j+}
\quad\text{and}\quad
\Upsilon' := \gamma'_j\CC G{j+}(\ugamma')_{x', j+}\,.
$$
Since
$$
\Upsilon \subseteq G_{x, j}
\quad\text{and}\quad
\Upsilon' \subseteq G_{x', j}\,,
$$
Lemma \ref{lem:depth-mod-center} gives
\begin{equation}
\tag{$*$}
z \in G_{x', j}\,.
\end{equation}
Thus
$\gamma_j\CC G{j+}(\ugamma)_{j+}
\cap G_{x', j}
\ne \emptyset$.
By Lemma \ref{lem:x-depth},
$x' \in \BB(\CC\bG{j+}(\ugamma), F)$.

Since $\Upsilon$ and $\Upsilon'$ are open, there is a semisimple element
$\delta$ of $\Upsilon \cap z\Upsilon'$.
Since
$\gamma_j\inv\delta \in \gamma_j\inv\Upsilon
	\subseteq \CC G{j{+}}(\ugamma)_{j{+}}$\,,
we have that $\gamma_j$ is a normal $(j{+})$-approximation
to $\delta$.  Thus, by Corollary
\ref{cor:compare-centralizers},
$C_G(\delta)\conn \subseteq \CC G{j{+}}(\ugamma)$.
In particular, $\gamma_j' \in C_G(\delta)\conn$ belongs to
$\CC G{j{+}}(\ugamma)$.
Since $x' \in \BB(\CC\bG{j{+}}(\ugamma), F)$ and
$\gamma_j'^{-1}z\inv\delta \in G_{x', j{+}}$\,,
we have that
$\gamma_j'^{-1}z\inv\delta \in \CC G{j{+}}(\ugamma)_{x', j{+}}$\,.
On the other hand, 
$z\inv\delta \in z\inv\Upsilon
	= z\inv\gamma_j\CC G{j{+}}(\ugamma)_{x, j{+}}$\,,
which (by Remark \ref{rem:approx-facts-in-center})
is contained in
$\ZZ G{j{+}}(\ugamma)\CC G{j{+}}(\ugamma)_{x, j{+}}$\,.
That is,
$$
\gamma_j'\CC G{j{+}}(\ugamma)_{x', j{+}}
	\cap \ZZ G{j{+}}(\ugamma)\CC G{j{+}}(\ugamma)_{x, j{+}}
\ne \emptyset.
$$
By Lemma \ref{lem:degen-mod-center},
$\gamma_j' \in \ZZ G{j{+}}(\ugamma)\CC G{j{+}}(\ugamma)_{j{+}}$\,.

Put $\bH = \CC\bG{j{+}}(\ugamma)$, and
let $\bT'$ be a maximal torus in \bH containing
$\gamma_j'$.
Denote by $E/F$ the splitting field of $\bT'$.
By Lemma \ref{lem:domain-field-ascent},
$\gamma_j' \in Z(\bH)(E)\bH(E)_{j{+}}$\,.
By Lemma \ref{lem:levi-descent},
$\bT'(E) \cap Z(\bH)(E)\bH(E)_{j{+}} = Z(\bH)(E)\bT'(E)_{j{+}}$\,.
Thus the root values of $\gamma'_j$ for \bH all
lie in $E\cross_{j+}$\,.
(The term ``root value'' is defined in Definition
\ref{defn:root-value}.)
Since
$\gamma'_j \in \GG_j^H \cup \sset 1$, the set of root
values of $\gamma'_j$ for \bH is therefore $\sset 1$; i.e.,
$\gamma'_j \in Z(H)$.  Thus
$\bH \subseteq C_\bG(\gamma'_j)$.
Since \bH is connected, in fact
$\bH \subseteq C_\bG(\gamma'_j)\conn
= \CC\bG{j+}(\ugamma')$.
The same argument (with the roles of \ugamma and
$\ugamma'$ reversed) gives the reverse containment; so
$\CC\bG{j+}(\ugamma) = \bH = \CC\bG{j+}(\ugamma')$.
This is ($\textbf{Cent}_{j+}$).

Now recall that
$\Upsilon = \gamma_j\CC G{j+}(\ugamma)_{x, j{+}} = \gamma_j H_{x, j+}$
and, similarly,
$\Upsilon' = \gamma_j'H_{x', j+}$\,,
so
$$
\gamma_j H_{x, j+}
\cap z\gamma'_j H_{x', j+}
=
\Upsilon \cap z\Upsilon'
\ne \emptyset;
$$
and
$\gamma_j$, $\gamma'_j$, and (by ($*$)) $z$ belong to
$Z(H) \cap H_j$\,.
By Lemma \ref{lem:depth-in-center},
$Z(H) \cap H_j = Z(H) \cap H_{x', j}$\,.
Therefore, the coset
$\gamma_j\inv z\gamma'_j H_{x', j+}$ is in
$H_{x', j:j+}$ and intersects $H_{x, j+}$\,.
By Lemma \ref{lem:unipotent},
$\gamma_j\inv z\gamma'_j \in H_{j+}$\,;
in fact, $\gamma_j\inv z\gamma'_j \in Z(H) \cap H_{j+}$\,.
If $\gamma_i = \gamma_i'$ for $0 \le i < j$, then $z = 1$,
so we have ($\textbf{Cong}_j'$).
In general, $z \in \ZZ G j(\ugamma)$
and (by ($*$)) $z \in \CC G j(\ugamma)_j$\,;
so we have ($\textbf{Cong}_j$).
\end{proof}

\section{Normal approximations and conjugation}
\label{sec:normal-conj}

We preserve the notation $\gamma$ and $r$ of the previous
section, so $\gamma \in G$ and $r \in \tR_{\ge 0}$.

We begin with a result which will play a technical role in
%Proposition \ref{prop:step1-formula1} of
\cite{adler-spice:explicit-chars}.

\begin{lm}
\label{lem:perp-commute}
Suppose that $i_0, j_0, t_0 \in \tR_{\ge 0}$
satisfy $i_0 + j_0 = t_0$.
Suppose further that
\begin{itemize}
\item $\gamma$ has a normal
$i_0$-approximation \ugamma,
\item $x \in \BB(\CC\bG{i_0}(\ugamma), F)$,
and
\item $\gamma_{\ge i_0} \in G_{x, i_0}$
(or $\stab_G(\ox)$, if $i_0 = 0$).
\end{itemize}
If $k \in G_{x, j_0+}$ and
$[\gamma, k] \in G_{x, t_0}$\,, then
$[\gamma, k] \in (\CC G{i_0}(\gamma), G)_{x, (t_0{+}, t_0)}$.
\end{lm}

\begin{proof}
We have that
$[\gamma, k] = [\gamma_{< i_0}, k]
		\dotm\lsup{\gamma_{< i_0}}[\gamma_{\ge i_0}, k]$.
Note that
$[\gamma_{\ge i_0}, k] \in G_{x, t_0+}$\,.
By Remark \ref{rem:approx-facts-in-stab},
$\gamma_{< i_0} \in \stab_G(\ox)$; in particular,
$\gamma_{< i_0}$ normalizes $G_{x, t_0{+}}$\,.
Therefore, $[\gamma_{< i_0}, k] \in G_{x, t_0}$\,.
By Lemma \ref{lem:center-comm}, in fact
$[\gamma_{< i_0}, k] \in (\CC G{i_0}(\gamma), G)_{x, (t_0{+}, t_0)}$.
The result follows.
\end{proof}

The next result shows that every $r$-approximation to
$\gamma$ is conjugate to a normal one.  After this result,
we will assume that $\gamma$ has a normal $r$-approximation.

\begin{lemma}
\label{lem:normal-conj}
If $(\ugamma,x)$ is an
$r$-approximation to $\gamma$,
then there is some $k\in [\ugamma; x, r]$
such that $\lsup k \ugamma$ is a normal $r$-approximation
to $\gamma$.
\end{lemma}

\begin{proof}
Trivially, every $0$-approximation to $\gamma$ is normal.
If $\gamma$ has an $\infty$-approximation, then $\gamma$ is
semisimple.
Then, by Lemma \ref{lem:connected} and Corollary
\ref{cor:compare-centralizers}, such an approximation
to $\gamma$ is normal.
Thus, we may, and hence do, assume that $0 < r < \infty$.
Write $\ugamma = (\gamma_i)_{0 \le i < r'}$\,.

Since all but finitely
many $\gamma_i$ are equal to $1$,
by Lemma \ref{lem:complete-subfield},
there is a
complete subfield $F'$ of $F$ such that $F/F'$ is separable,
and $\gamma$ and each $\gamma_i$ are defined over $F'$.
Since $[\ugamma; x, r]_{\bG(F')} \subseteq [\ugamma; x, r]$
by Remark \ref{rem:bracket-facts-field-descent},
it suffices to assume
that $F$ is complete.

In this setting, we need only show that,
if $d \in \R_{\ge 0}$ is such that $d < r$
and $(\ugamma, x)$ is a normal $d$-approximation
(as well as an $r$-approximation), then there is some $k\in [\ugamma;x,r]$
such that $\lsup k\ugamma$ is a normal $(d{+})$-approximation.
Let $\bH = \CC\bG{d}(\ugamma)$.
By Definition \ref{defn:r-approx},
$x\in \BB(\bH,F)$ and
$$
\gamma
\in \bigl(\prod_{0 \le i < r} \gamma_i\bigr)G_{x,r} \cap H
	= \bigl(\prod_{0 \le i < r} \gamma_i\bigr)H_{x,r}\,.
$$
By Corollary \ref{cor:good-comm},
there is some $k\in H_{x,r-d} \subseteq [\ugamma; x, r]$ such that
$$
\bigl(\prod_{0 \le i < d} \gamma_i\bigr)\inv(\lsup{k\inv}\gamma)
=
\lsup{k\inv}\bigl[ \bigl(\prod_{0\leq i < d} \gamma_i\bigr)\inv\gamma \bigr]
\in
\bigl(\prod_{d\leq i < r} \gamma_i\bigr) C_H(\gamma_d)_{x,r}\,.
$$
That is, $\lsup{k\inv}\gamma$ commutes with $\gamma_i$ for
$0 \le i \le d$.
By Lemma \ref{lem:connected},
$\lsup k\ugamma = (\lsup k\gamma_i)_{0 \le i < r'}$ is a
normal $(d{+})$-approximation to $\gamma$.
\end{proof}

From now on, assume that
$\gamma$ has a normal $r$-approximation
$\ugamma = (\gamma_i)_{0 \le i < r}$\,.
By Proposition \ref{prop:unique-approx}, the choice of this
approximation will not affect any of the results
or definitions that follow.

\begin{dn}
\label{defn:fancy-centralizer-no-underline}
\indexmem{\CC{\bG}r(\gamma)}
\indexmem{\CC G r(\gamma)}
\indexmem{\ZZ{\bG} r(\gamma)}
\indexmem{\ZZ G r(\gamma)}
\indexmem{[\gamma; x, r]^{(j)}}
\indexmem{[\gamma; x, r]}
\indexmem{\dc^{(j)}}
\indexmem\dc
Put
$\CC\bG r(\gamma) = \CC\bG r(\ugamma)$,
$\CC G r(\gamma) = \CC G r(\ugamma)$,
$\ZZ\bG r(\gamma) = \ZZ\bG r(\ugamma)$,
$\ZZ G r(\gamma) = \ZZ G r(\ugamma)$,
$[\gamma; x, r]^{(j)} = [\ugamma; x, r]^{(j)}$,
$[\gamma; x, r] = [\gamma; x, r]^{(\infty)}$,
$\dc^{(j)} = \udc^{(j)}$,
and
$\dc = \dc^{(\infty)}$
for $(j, x) \in \tR \times \BB(\CC\bG r(\gamma), F)$.
\end{dn}

\begin{rk}
\label{rem:sloppy-bracket}
If $(\ugamma', x)$ is a (not necessarily normal)
$r$-approximation to $\gamma$, then, by Lemma
\ref{lem:normal-conj}, there is $k \in [\ugamma'; x, r]$
such that $\lsup k\ugamma'$ is a normal $r$-approximation to
$\gamma$.  Then 
$$
[\gamma; x, r]
= [\lsup k\ugamma'; x, r] = \lsup k[\ugamma'; x, r]
= [\ugamma'; x, r].
$$
\end{rk}

\begin{dn}
\label{defn:Brgamma}
Put
\indexmem{\BB_r(\gamma)}
\begin{align*}
\BB_r(\gamma)
& = \set{x \in \BB(\bG, F)}{\gamma\ox = \ox}
& \text{if } & r = 0, \\
\BB_r(\gamma)
& = \set{x \in \BB(\CC\bG r(\gamma), F)}
	{\ZZ G r(\gamma)\gamma \cap G_{x, r} \ne \emptyset}
& \text{if } & r \in \tR_{> 0}.
\end{align*}
\end{dn}

The next result shows that, for $r > 0$,
$\BB_r(\gamma)$ is precisely the set of points $x$ such
that $(\ugamma, x)$ is a normal $r$-approximation.
(In fact, that is also true for $r = 0$.)

\begin{lm}
\label{lem:Brgamma}
If $r > 0$, then
$$
\BB_r(\gamma)
= \set{x \in \BB(\CC\bG r(\gamma), F)}
      {\depth_x(\gamma_{\ge r}) \ge r}.
$$
\end{lm}

\begin{proof}
By Remark \ref{rem:approx-facts-in-center},
$\gamma_{\ge r} = \gamma_{< r}\inv\gamma
\in \ZZ G r(\gamma)\gamma$.
By Definition \ref{defn:r-approx},
$\gamma_{\ge r} \in G_{y, r} \cap \CC G r(\gamma)
	= \CC G r(\gamma)_{y, r} \subseteq \CC G r(\gamma)_r$
(for some $y \in \BB(\CC\bG r(\gamma), F)$).

Suppose that $x \in \BB_r(\gamma)$.  Then
$\ZZ G r(\gamma)\gamma \subseteq \ZZ G r(\gamma)G_{x, r}$\,,
so
$\gamma_{\ge r} \in \ZZ G r(\gamma)G_{x, r} \cap \CC G r(\gamma)
                = \ZZ G r(\gamma)\CC G r(\gamma)_{x, r}$\,.
By Corollary \ref{cor:depth-mod-center},
$\gamma_{\ge r} \in G_{x, r}$\,.

On the other hand, suppose that
$\depth_x(\gamma_{\ge r}) \ge r$.  Then
$\gamma_{\ge r} \in \ZZ G r(\gamma)\gamma \cap
G_{x, r}$\,, so the intersection is non-empty.
\end{proof}

\begin{cor}
\label{cor:x-depth}
If $\gamma$ has a normal $(r{+})$-approximation, then
$\BB_r(\gamma) \subseteq \BB(\CC\bG{r+}(\gamma), F)$.
\end{cor}

\begin{proof}
This follows immediately from Lemmata \ref{lem:x-depth} and
\ref{lem:Brgamma}.
\end{proof}

The next result is an analogue of Lemma 3.6 of
\cite{murnaghan:chars-gln}.

\begin{lm}
\label{lem:bracket}
Fix $d \in \tR_{\ge 0}$.
Suppose that
\begin{itemize}
\item $x \in \BB_r(\gamma)$,
\item $g \in G_{x, \max\sset{d, 0{+}}}$\,, and
\item $[\gamma, g] \in G_{x, d + r}$\,.
\end{itemize}
Then $g \in [\gamma; x, d + r]$.
\end{lm}

\begin{proof}
Since $[\gamma; x, 0] = G_{x, 0+}$\,, the result is trivial
if $d + r = 0$; so we assume that $d + r > 0$.

Let \bT be a maximal $F$-torus in $\CC\bG r(\gamma)$, and
$E/F$ the splitting field of \bT.
By Lemma \ref{lem:field-descent},
$$
g \in \bG(E)_{x, d}
\quad\text{and}\quad
[\gamma, g] \in \bG(E)_{x, d + r}\,.
$$
If $r = 0$, then, since $\gamma = \gamma_{\ge r}$
stabilizes \ox, hence $\bG(E)_{x, d}$\,, we have
$[\gamma_{\ge r}, g] \in \bG(E)_{x, d} = \bG(E)_{x, d + r}$\,.
If $r > 0$, then
$\gamma_{\ge r} \in G_{x, r} \subseteq \bG(E)_{x, r}$\,, so
again
$[\gamma_{\ge r}, g] \in G_{x, d + r} \subseteq \bG(E)_{x, d + r}$\,.
By Remark \ref{rem:approx-facts-in-stab},
$\gamma_{< r} \in \stab_{\bG(E)}(\ox)$, so
$\lsup{\gamma_{< r}}[\gamma_{\ge r}, g] \in \bG(E)_{x, d + r}$\,.
Since
$[\gamma, g]
= \lsup{\gamma_{< r}}[\gamma_{\ge r}, g]
	\cdot[\gamma_{< r}, g]$,
also $[\gamma_{< r}, g] \in \bG(E)_{x, d + r}$\,;
i.e.,
$\Int(\gamma_{< r})g \equiv g \pmod{\bG(E)_{x, d + r}}$.
By Remark \ref{rem:actually-product}, the multiplication
map
\begin{equation}
\tag{$*$}
\prod_{\alpha \in \wtilde\Phi(\bG, \bT)}
      (\bU_\alpha(E) \cap \bG(E)_{x, \max\sset{d, 0{+}}})
\to \bG(E)_{x, \max\sset{d, 0{+}}}
\end{equation}
is a bijection.
Denote by $(g_\alpha)_{\alpha \in \wtilde\Phi(\bG, \bT)}$ the
preimage of $g$ under ($*$).
Let
$(\mexp_\alpha : E \to \bU_\alpha(E))
	_{\alpha \in \Phi(\bG, \bT)}$
be the isomorphisms of \S\ref{sec:unip-exp},
and, for $\alpha \in \wtilde\Phi(\bG, \bT)$ and $t \in E$,
let $\mult_{\alpha, t}$ be the endomorphism of $\bU_\alpha(E)$
given by
$$
u \mapsto
\begin{cases}
\mexp_\alpha(t\dotm\mexp_\alpha\inv(u)), &
	\alpha \in \Phi(\bG, \bT) \\
u, &
	\alpha = 0.
\end{cases}
$$
Then the preimage of
$\Int(\gamma_{< r})g$ under ($*$) is
$(\mult_{\alpha, \alpha(\gamma_{< r})}(g_\alpha))
	_{\alpha \in \wtilde\Phi(\bG, \bT)}$.
Recall that $\Int(\gamma_{< r})g \equiv g \pmod{\bG(E)_{x, d + r}}$.
By Corollary \ref{cor:gen-iwahori-factorization-bij}
(applied to the constant functions $f_1$ and $f_2$ with
values $\max\sset{d, 0{+}}$ and $d + r$, respectively)
gives that
$\mult_{\alpha, \alpha(\gamma_{< r})}(g_\alpha)
\equiv g_\alpha \pmod{\bG(E)_{x, d + r}}$\,,
i.e.,
$$
\mult_{\alpha, \alpha(\gamma_{< r}) - 1}(g_\alpha)
= \mult_{\alpha, \alpha(\gamma_{< r})}(g_\alpha)
  \cdot g_\alpha\inv
\in \bU_\alpha(E) \cap \bG(E)_{x, d + r}\,,
$$
for $\alpha \in \Phi(\bG, \bT)$.

Fix $\alpha \in \Phi(\bG, \bT)$.
If $\alpha(\gamma_i) = 1$ for $0 \le i < r$, then
$\bU_\alpha \subseteq \CC\bG r(\gamma)$, so
$g_\alpha \in \CC\bG r(\gamma)(E) \cap \bG(E)_{x, d}
\subseteq [\gamma; x, d + r]$.
If $\alpha(\gamma_i) \ne 1$ for some $0 \le i < r$,
then, by Lemma \ref{lem:compare-centralizers},
$i_\alpha := \ord(\alpha(\gamma_{< r}) - 1)
= \min \set i{\alpha(\gamma_i) \ne 1}$.
By Remark \ref{rem:actually-product},
$\bU_\alpha(E) \cap \bG(E)_{x, d + r}
= U_\psi$\,,
where $\psi$ is the affine function on $\AA(\bT, E)$ with
gradient $\alpha$ satisfying $\psi(x) = d + r$.
Then
$\bU_\alpha \subseteq \CC\bG{i_\alpha}(\gamma)$,
so
$$
U_{\psi - i_\alpha}
\subseteq \CC\bG{i_\alpha}(\gamma)(E)_{x, (d + r) - i_\alpha}
\subseteq [\gamma; x, d + r]_{\bG(E)}.
$$
In particular,
$$
g_\alpha
\in \mult_{\alpha, \alpha(\gamma_{< r}) - 1}\inv
    (U_\psi)
= U_{\psi - i_\alpha}
\subseteq [\gamma; x, d + r]_{\bG(E)}.
$$
By Remark \ref{rem:bracket-facts-field-descent},
$[\gamma; x, d + r]_{\bG(E)} \cap G_{x, 0+}
= [\gamma; x, d + r]_G$\,,
so $g = \prod g_\alpha
\in [\gamma; x, d + r]_G$\,.
\end{proof}

\begin{rk}
We will often apply the above lemma with
$d = 0$ and $r = \depth_x([\gamma, g])$.
\end{rk}

The following ``rigidity'' result will help us
in \cite{adler-spice:explicit-chars}
to apply Harish-Chandra's character formula.
% Specifically, we will use this result in 
% Theorem \ref{thm:char-tau|pi-1} of \cite{adler-spice:explicit-chars}.

\begin{lm}
\label{lem:rigidity}
Suppose that
\begin{itemize}
\item $(\bG', \bG)$ is a tame reductive $F$-sequence;
\item $x \in \BB(\CC\bG r(\gamma), F) \cap \BB(\bG', F)$;
\item $k \in G_{x, 0+}$\,;
and
\item $\ZZ G r(\gamma)$ and $\lsup k\ZZ G r(\gamma)$ are
contained in $G'$.
\end{itemize}
Then
$k \in G'_{x, 0+}\CC G r(\gamma)_{x, 0+}$\,.
\end{lm}

Note that the statement can also be rewritten as
$\lsup{G_{x, 0+}}\gamma \cap G'
= \lsup{G'_{x, 0+}}\gamma$
(under the stated conditions on $\gamma$ and $x$).

\begin{proof}
Since the result does not change if we replace $\gamma$ by
$\gamma_{< r}$, we do so.
Then $\gamma_{\ge r} = 1$, so, by Lemma \ref{lem:Brgamma},
$\BB_r(\gamma) = \BB(\CC\bG r(\gamma), F)$.
In particular, $x \in \BB_r(\gamma)$.

Let $i_0$ be the greatest index $i$ such that
$k \in G'_{x, 0+}\CC G i(\gamma)_{x, 0{+}}$\,.
If $i_0 \ge r$, then we are done, so we may, and hence do,
assume that $i_0 < r$.
Since the hypotheses and conclusion do not change if we replace $k$ by a
left $G'_{x, 0+}$-translate, we may, and hence do, assume
further that
$k \in \CC G{i_0}(\gamma)_{x, 0{+}}$\,.

By Remark \ref{rem:approx-facts-in-center},
$\gamma_i \in \ZZ G r(\gamma) \subseteq G'$,
and similarly $\lsup k\gamma_i \in G'$,
for $0 \le i < r$.
By applying Lemma \ref{lem:x-depth} repeatedly, one sees
that $x$ belongs to $\BB(\CC{\bG'}{i_0}(\gamma), F)$
(in fact, to $\BB(\CC{\bG'}r(\gamma), F)$).
Thus we may, and hence do, replace \bG and $\bG'$ by
$\CC\bG{i_0}(\gamma)$ and $\CC{\bG'}{i_0}(\gamma)$,
respectively.
Now $\gamma_{< i_0}$ and $\gamma_{\le i_0}$, as well as
their $k$-conjugates, belong to $G'$;
so $\gamma_{\ge i_0}$ and $\gamma_{> i_0}$, as well as their
$k$-conjugates, belong to $G'$.

Since $x \in \BB_r(\gamma)$, \emph{a fortiori}
$x \in \BB_{i_0}(\gamma)$ and $x \in \BB_{i_0{+}}(\gamma)$;
so $\gamma_{\ge i_0} \in G_{x, i_0}$ (or $\stab_G(\ox)$, if
$i_0 = 0$) and $\gamma_{> i_0} \in G_{x, i_0{+}}$\,.
Thus
$$
\lsup k\gamma_{\ge i_0}
\in \gamma_{\ge i_0}G_{x, i_0{+}} \cap G'
= \gamma_{i_0}G_{x, i_0{+}} \cap G'
= \gamma_{i_0}G'_{x, i_0{+}}\,,
$$
so $((\gamma_{i_0}), x)$ is an $(i_0{+})$-approximation (in
$G'$) to $\lsup k\gamma_{\ge i_0} \in G'$.
By Lemma \ref{lem:normal-conj}, there exists
$h \in G'_{x, 0{+}}$ such that $((\lsup h\gamma_{i_0}), x)$ is a
normal $(i_0{+})$-approximation to $\lsup k\gamma_{\ge i_0}$;
that is,
$$
\lsup k\gamma_{\ge i_0}
\in (\lsup h\gamma_{i_0})\dotm C_G(\lsup h\gamma_{i_0})_{x, i_0{+}}
= \lsup h(\gamma_{i_0}C_G(\gamma_{i_0})_{x, i_0{+}}).
$$
By Lemma \ref{lem:GE2},
$h\inv k \in C_G(\gamma_{i_0})$.
In fact, since $h, k \in G_{x, 0{+}}$
and
$x \in \BB(C_\bG(\gamma_{i_0}), F)$,
we have by Corollary \ref{cor:compatibly-filtered-tame-rank}
that $h\inv k \in C_G(\gamma_{i_0})_{x, 0{+}}$\,.
Note that, since $\bG = \CC\bG{i_0}(\gamma)$, we have
$C_\bG(\gamma_{i_0})\conn = \CC\bG{i_0{+}}(\gamma)$, so that
$h\inv k \in \CC G{i_0{+}}(\gamma)_{x, 0{+}}$\,.
Since $h \in G'_{x, 0{+}}$\,, we have
$k \in G'_{x, 0{+}}\CC G{i_0{+}}(\gamma)_{x, 0{+}}$\,, which
is a contradiction of the definition of $i_0$.
\end{proof}

Now we prove a few technical lemmata which will come in
handy in the proof of Proposition \ref{prop:aniso-Levi}.

\begin{lm}
\label{lem:x-depth-aniso}
Suppose that
\begin{itemize}
\item $t \in \tR_{\ge 0}$,
\item $(\bG', \bG)$ is a tame Levi $F$-sequence,
\item $\bG'/Z(\bG')$ is $F$-anisotropic,
and
\item $x \in \BB(\bG', F)$.
\end{itemize}
Then
$(G' \cap \stab_G(\ox))G_{x, t} \cap G_0 \subseteq G_{x, 0}$\,,
and $\depth_x(g) \ge \min\sset{t, \depth(g)}$ for
$g \in (G' \cap \stab_G(\ox))G_{x, t} \cap G_0$\,.
\end{lm}

\begin{proof}
The containment follows from Lemma
\ref{lem:stab-deep}.

Suppose that
$g \in (G' \cap \stab_G(\ox))G_{x, t} \cap G_0$
and
$j := \depth_x(g) < \depth(g)$.
If $j < t$, then, by Lemma \ref{lem:unipotent},
$g G_{x, j+} \subseteq G_{j+}$\,.
By Lemmata \ref{lem:domain-field-ascent},
\ref{lem:levi-descent},
and \ref{lem:field-descent}, we have
$G' \cap G_{j+} = G'_{j+} = G'_{x, j+}$\,.
Thus, since $g G_{x, j+} \cap G' \ne \emptyset$,
in fact
$g G_{x, j+} \cap G'_{x, j+} \ne \emptyset$;
so $g \in G_{x, j+}$\,, which is a contradiction.
\end{proof}

\begin{lm}
\label{lem:dont-conj}
Suppose that
\begin{itemize}
\item $t \in \tR_{\ge 0}$,
\item $(\bG', \bG)$ is a tame Levi $F$-sequence,
\item $\ZZ G t(\gamma) \subseteq \lsup{[\gamma; x, t]}G'$,
\item $x \in \BB(\bG', F) \cap \BB_t(\gamma)$,
and
\item $\gamma \in \lsup{[\gamma; x, t]}(G'G_{x, t{+}})$.
\end{itemize}
Then $\gamma \in (G' \cap \stab_G(\ox))G_{x, t+}$\,.
\end{lm}

\begin{proof}
Note that $\gamma_{< t} \in \stab_G(\ox)$ by
Remark \ref{rem:approx-facts-in-stab},
and $\gamma_{\ge t} \in \stab_G(\ox)$ by
Lemma \ref{lem:Brgamma} (if $t > 0$) or by Definition
\ref{defn:Brgamma} (if $t = 0$).
Thus also $\gamma \in \stab_G(\ox)$.

In particular, if $t = 0$, then, since
$[\gamma; x, 0] \subseteq G_{x, 0+}$ (in fact, we have
equality),
any $[\gamma; x, t]$-conjugate of $\gamma$ lies in
$\gamma G_{x, 0+}$\,; so the result is easy.
Thus we may, and hence do, assume that $t > 0$.

By hypothesis, there is $k \in [\gamma; x, t]$ such that
$\lsup k\gamma_{< t} \in G'$.
By Corollary \ref{cor:compare-centralizers},
$\CC\bG t(\lsup k\gamma) = C_\bG(\lsup k\gamma_{< t})\conn$,
so
$\ZZ\bG t(\gamma) = Z(C_\bG(\lsup k\gamma_{< t})\conn)
\subseteq \bG'$.
Upon replacing $\gamma$ by $\lsup k\gamma$,
we may, and hence do, assume
that $\ZZ G t(\gamma) \subseteq	G'$.  In particular, by
Remark \ref{rem:approx-facts-in-center},
$\gamma_{< t} \in G'$; so,
by hypothesis, $\gamma \in G'G_{x, t}$\,.
Since $\bG'$ is an $F\tame$-Levi $F$-subgroup of \bG, there
is a tame $F$-torus $\bS' \subseteq \bG$ such that
$\bG' = C_\bG(\bS')$.
Since $\gamma_{< t} \in G'$, we have
$\bS' \subseteq C_\bG(\gamma_{< t})\conn$;
so, by Corollary \ref{cor:compare-centralizers},
$\bS' \subseteq \CC\bG t(\gamma)$.

Let \bS be a maximal tame-modulo-center $F$-torus
containing $\bS'$, hence contained in $\bG'$.
By Definitions \ref{defn:good} and \ref{defn:r-approx},
the terms $\gamma_i$ for $0 \le i < t$ belong to a
tame-modulo-center $F$-torus in $\CC\bG t(\gamma)$.
Therefore, by Lemma \ref{lem:split-in-center}, the terms
$\gamma_i$ for $0 \le i < t$ in fact belong to $S$.
Thus $(\gamma_i)_{0 \le i < t}$ is a good sequence in $G'$.
By abuse of language, we will write (for example)
$\CC{\bG'}t(\gamma)$ in place of
$\CC{\bG'}t((\gamma_i)_{0 \le i < t})$,
even though possibly $\gamma \not\in G'$.
Let \bT be a maximal $F$-torus in $\bG'$ containing
$\bS$.
Now recall that
$[\gamma; x, t]_G = \vG_{x, \vec t}$\,, where
$\vG := (\CC\bG i(\gamma))_{0 \le i < t}$ and
$\vec t := (t - i)_{0 \le i < t}$; and, similarly,
$[\gamma; x, t]_{G'} = \vG^{\,\prime}_{x, \vec t}$\,, where
$\vG^{\,\prime} := (\CC{\bG'}i(\gamma))_{0 \le i < t}$.
Let $f^\perp$ be the function on $\wtilde\Phi(\bG, \bT)$
such that
\begin{align*}
f^\perp(\alpha) &
= \min\set{(t - i){+}}{\alpha \in \Phi(\CC\bG i(\gamma), \bT)} \\
\intertext{for $\alpha \in \Phi(\bG', \bT)
	\smallsetminus \Phi(\CC\bG t(\gamma), \bT)$;}
f^\perp(\alpha) & = 0{+} \\
\intertext{for $\alpha \in \wtilde\Phi(\CC\bG t(\gamma), \bT)$;
and}
f^\perp(\alpha) &
= \min\set{t - i}{\alpha \in \Phi(\CC\bG i(\gamma), \bT)} \\
\end{align*}
otherwise.  Note that this function is concave.
By hypothesis, we may, and hence do,
choose $k \in [\gamma; x, t]_G$ such that
$\lsup k\gamma \in G'G_{x, t+}$\,.
Then, by Proposition \ref{prop:heres-a-gp},
there exist $k' \in [\gamma; x, t]_{G'}$ and
$k^\perp \in \lsub\bT G_{x, f^\perp}$ such that
$k = k'k^\perp$.
Modulo $G_{x,t+}$\,, we have
\begin{multline*}
\lsup k\gamma
= \lsup {k'} (\gamma [\gamma\inv, k^\perp])
\equiv \lsup {k'} \gamma \dotm {[\gamma\inv, k^\perp]} \\
= [k', \gamma] \dotm \gamma \dotm {[\gamma\inv, k^\perp]}
\equiv [k', \gamma_{<t}] \dotm \gamma_{<t} \gamma_{\geq t}
	\dotm {[\gamma\inv_{<t}, k^\perp]}.
\end{multline*}
Since $\gamma_{< t} \in G'$, also
$[k', \gamma_{< t}] \in G'$.
Thus
$$
\gamma_{\ge t}\dotm{[\gamma_{< t}\inv, k^\perp]}
	\in G'G_{x, t+}\,.
$$
By Remark \ref{rem:approx-facts-commutator},
$[\gamma_{< t}\inv, k^\perp] \in G_{x, t}$\,.
Thus in fact
$$
\gamma_{\ge t}\dotm{[\gamma_{< t}\inv, k^\perp]}
	\in G'G_{x, t+} \cap G_{x, t} = (G', G)_{x, (t, t{+})}.
$$
Recall that $\gamma_{\ge t} \in \CC G t(\gamma)_{x, t}$\,, so
\begin{equation}
\tag{$*$}
[\gamma_{< t}\inv, k^\perp]
\in \gamma_{\geq t}\inv(G', G)_{x, (t, t{+})}
	\subseteq \lsub\bT G_{x, f}\,,
\end{equation}
where $f$ is the function on $\wtilde\Phi(\bG, \bT)$ which
takes the value $t$ on
$\wtilde\Phi(\bG', \bT) \cup \wtilde\Phi(\CC\bG t(\gamma), \bT)$
and the value $t{+}$ elsewhere.
In particular, $[\gamma_{< t}\inv, k^\perp] \in G_{x, t}$\,.
Then we have
\begin{equation}
\tag{$**$}
[\gamma_{< t}\inv, k^\perp] \in (\CC G t(\gamma), G)_{x, (t{+}, t)}
\end{equation}
by Lemma \ref{lem:center-comm}.
Finally, by Proposition \ref{prop:heres-a-gp}, we may write
$k^\perp = \prod_{0 \le i < t} k^\perp_i$, where
$k^\perp_i \in (\CC{G'}i(\gamma), \CC G i(\gamma))
	_{x, ((t - i){+}, t - i)}$ for $0 \le i < t$.
Put $\gamma' := \gamma_{< t}$.
Then, with the obvious notation, we have
$$
[\gamma'^{-1}, k^\perp]
= \prod_{0 \le i < t} \lsup{k^\perp_{< i}}
	[\gamma'^{-1}, k^\perp_i]
= \prod_{0 \le i < t} \lsup{k^\perp_{< i}}
	[\gamma_{\ge i}'^{-1}, k^\perp_i].
$$
By Lemma \ref{lem:shallow-comm},
$[\gamma_{\ge i}'^{-1}, k^\perp_i]
\in (\CC{G'}i(\gamma), \CC G i(\gamma))_{x, (t{+}, t)}
\subseteq (G', G)_{x, (t{+}, t)}$.
In particular,
$[\gamma_{\ge i}'^{-1}, k^\perp_i] \in G_{x, t}$\,.
Since $k^\perp_{< i} \in G_{x, 0+}$\,, the commutator of
$k^\perp_{< i}$ with $[\gamma_{\ge i}'^{-1}, k^\perp_i]$
lies in $G_{x, t{+}} \subseteq (G', G)_{x, (t{+}, t)}$,
so
$\lsup{k^\perp_{< i}}[\gamma_{\ge i}'^{-1}, k^\perp_i]
\in (G', G)_{x, (t{+}, t)}$ also.
Thus
\begin{equation}
\tag{$*{*}*$}
[\gamma'^{-1}, k^\perp] \in (G', G)_{x, (t{+}, t)}.
\end{equation}
(It is also possible to show this directly, without
appealing to Proposition \ref{prop:heres-a-gp}, as in
Remark \ref{rem:approx-facts-commutator};
but the argument is more complicated.)
By Lemma \ref{lem:more-vGvr-facts} and
Remark \ref{rem:vGvr-facts},
it follows from ($*$), ($**$), and ($*{*}*$) that
$[\gamma_{< t}\inv, k^\perp] \in G_{x, t{+}}$\,,
hence also $[\gamma\inv, k^\perp] \in G_{x, t+}$\,.
Further,
$\lsup{k^\perp}\gamma = \gamma [\gamma\inv, k^\perp] \in \gamma G_{x, t+}$
and
$ \lsup{k^\perp}\gamma = \lsup{k'^{-1}k}\gamma
\in \lsup{k'^{-1}}(G'G_{x, t+}) = G'G_{x, t+}$\,,
so $\gamma \in G'G_{x, t+}$\,.
\end{proof}

Recall that $r \in \tR_{\ge 0}$.

\begin{lm}
\label{lem:aniso-Brgamma}
Suppose that
\begin{itemize}
\item $(\bG', \bG)$ is a tame Levi $F$-sequence,
\item $\bG'/Z(\bG')$ is $F$-anisotropic,
\item $x \in \BB(\bG', F)$,
and
\item $\gamma \in (G' \cap \stab_G(\ox))G_{x, r}$\,.
\end{itemize}
Then $x \in \BB_r(\gamma)$.
\end{lm}

\begin{proof}
If $r = 0$, then, since
$\BB_0(\gamma)$ is the set of points in $\BB(\bG, F)$ whose image
in $\rBB(\bG, F)$ is fixed by $\gamma$, in particular
$\sset\ox \times V_F(Z(\bG)) \subseteq \BB_0(\gamma)$;
so we are done.
If $r = \infty$, then
$\BB_r(\gamma) = \BB(\CC\bG\infty(\gamma), F)
	= \BB(\CC\bG{r'}(\gamma), F)$
for $r' \in \R$ sufficiently large.
Thus we may, and hence do, assume that $0 < r < \infty$.

Now it suffices to show that, if
$t \in \R_{\ge 0}$ with $t < r$,
and
$x \in \BB_t(\gamma)$,
then $x \in \BB_{t+}(\gamma)$.
Since $\BB_t(\gamma) \subseteq \BB(\CC\bG{t+}(\gamma), F)$
by Corollary \ref{cor:x-depth}, Lemma \ref{lem:Brgamma}
shows that it suffices to prove $\gamma_{> t} \in G_{x, t+}$\,.
By Lemma \ref{lem:x-depth-aniso},
since $\depth(\gamma_{> t}) > t$,
it suffices to prove that
$\gamma_{> t}
\in \lsup{\stab_G(\ox)}
	\bigl((G' \cap \stab_G(\ox))G_{x, t+}\bigr)$.

First suppose that $t = 0$.
Write $\gamma = \gamma'\gamma^\perp$,
with $\gamma' \in \stab_{G'}(\ox)$
and $\gamma^\perp \in G_{x, r}$\,.
By Proposition \xref{J-prop:tJd=cpct} of \cite{spice:jordan},
$\ogamma^{\,\prime}$ has a topological Jordan decomposition
$\ogamma^{\,\prime}
= \ogamma^{\,\prime}\tsemi\ogamma^{\,\prime}\tunip$
in $\tbG'$
(where $\tbG' = \bG'/Z(\bG)\conn$).
Then $((\ogamma^{\,\prime}\tsemi), \ox)$ is a normal $(0{+})$-approximation to
$\ogamma'$, so $((\ogamma^{\,\prime}\tsemi), \ox)$
is a $(0{+})$-approximation to \ogamma.
By Lemma \ref{lem:normal-conj},
there is $\ol k \in \tG_{\ox, 0{+}}$ such that
$(\lsup{\ol k}\ogamma^{\,\prime}\tsemi)$ is a normal
$(0{+})$-approximation to \ogamma;
in particular,
$\ogamma = \lsup{\ol k}\ogamma^{\,\prime}\tsemi\dotm
	(\lsup{\ol k}\ogamma^{\,\prime\,{-1}}\tsemi)\ogamma$
is a topological Jordan decomposition.
On the other hand, $\ogamma = \ogamma_0\ogamma_{> 0}$ is
also a topological Jordan decomposition, so, by
Proposition \xref{J-prop:unique-top-F-Jordan},
$\ogamma_0 = \lsup{\ol k}\ogamma^{\,\prime}\tsemi
	\in \ogamma^{\,\prime}\tsemi\tG_{\ox, 0{+}}
	= \ogamma\tG_{\ox, 0{+}}$\,.
Thus
$\ogamma_{> 0} = \ogamma_0\inv\ogamma \in \tG_{\ox, 0{+}}$\,,
so there is a small neighborhood $U$ of $x$ in
$\BB(\bG, F)$ so that
$\ogamma_{> 0} \in \tG_{\ox', 0{+}}$
for $x' \in U$.
Thus, by Lemma \ref{lem:stab-deep},
$\gamma_{> 0}
\in G_{0+} \cap \bigcap_{x' \in U} \stab_G(\ox')
= G_{0+} \cap \bigcap_{x' \in U} G_{x'}$\,.
By Lemma \ref{lem:unipotent}, the image of
$\gamma_{> 0}$ in $(\ms G^F_x)\conn(\ol\ff)$ is unipotent.
By Proposition 5.1.32(i) of
\cite{bruhat-tits:reductive-groups-2},
for every parabolic subgroup
$\ms P \subseteq \ms G_x$, there is a facet
$J = J_{\ms P}$ of $\BB(\bG, F\unram)$ containing
$x$ in its closure such that $\ms P(\ol\ff)$ is the image in
$\ms G_x(\ol\ff)$ of $\bG(F\unram)_y$ for any $y \in J$.
Now
$\tG_{\ox, 0{+}} \subseteq \tbG(F\unram)_{\ox, 0{+}}
\subseteq \tbG(F\unram)_{\ol y, 0{+}}$\,, so
$\gamma_{> 0} \in \stab_{\bG(F\unram)}(\ol y)$,
for $y \in \bigcup_{\ms P} J_{\ms P}$.
Since
$\gamma_{> 0} \in G_{0{+}} \subseteq \bG(F\unram)_{0{+}}$\,,
Lemma \ref{lem:stab-deep} gives
$\gamma_{> 0} \in \bG(F\unram)_y$
for $y \in \bigcup_{\ms P} J_{\ms P}$;
that is, the image of $\gamma_{> 0}$ in
$\ms G_x(\ol\ff)$ lies in
$\bigcap_{\ms P} \ms P(\ol\ff)$.
By Theorem 13.16 of \cite{borel:linear} and the fact that
$\ms G_x$ is reductive, we have that the image of
$\gamma_{> 0}$ lies in every maximal torus in $\ms G_x$;
in particular, is semisimple.
On the other hand, since $\gamma_{> 0} \in G_{0{+}}$\,,
by Lemma \ref{lem:unipotent} there is a facet $J$ containing
$x$ in its closure such that
$\gamma_{> 0} \in G_{y, 0{+}}$ for $y \in J$.
By another application of Proposition 5.1.32(i) of
\cite{bruhat-tits:reductive-groups-2},
the image of $\gamma_{> 0}$ in $\ms G_x(\ff)$ is unipotent.
Thus this image is trivial, so
$\gamma_{> 0} \in G_{x, 0{+}}$\,, as desired.

Now suppose that $t > 0$.
Since $\gamma_{\ge t} \in G_{x, t}$\,, we have that
$\gamma_{< t}G_{x, t} \cap (G' \cap \stab_G(\ox)) \ne \emptyset$.
Since the intersection is an open subset of $G'$, it
contains a semisimple element, say $g'$.
Then \ugamma is a $t$-approximation to $g'$, so, by Lemma
\ref{lem:normal-conj}, there is an element
$k \in [\gamma; x, t]$ such that $\lsup k\ugamma$
is a normal $t$-approximation to $g'$.
In particular, by Corollary \ref{cor:compare-centralizers},
$\ZZ\bG t(\lsup k\ugamma) \subseteq Z(C_\bG(g')\conn)
\subseteq \bG'$;
so, by Remark \ref{rem:approx-facts-in-center},
$\lsup k\gamma_i \in G'$ for $0 \le i < t$
(hence also $\lsup k\gamma_{< t} \in G'$).
Put $\gamma' = \lsup k\gamma$,
and note that $x \in \BB_t(\gamma')$
(hence, by Corollary \ref{cor:x-depth},
$x \in \BB(\CC\bG{t+}(\gamma'), F)$).
Then
$[\gamma'; x, t] = \lsup k[\gamma; x, t] = [\gamma; x, t]$,
so
$\gamma'
\in \lsup{[\gamma'; x, t]}
	\bigl((G' \cap \stab_G(\ox))G_{x, t+}\bigr)$.
By Lemma \ref{lem:dont-conj},
$\gamma' \in G'G_{x, t+}$\,.
Since $\gamma'_{< t} = \lsup k\gamma_{< t} \in G'$, also
$\gamma'_{\ge t} \in G'G_{x, t+}$\,.

As before, there is a semisimple element
$g'' \in \gamma'_{\ge t}G_{x, t+} \cap G'$.
Since $x \in \BB_t(\gamma)$,
and $\BB_t(\gamma) \subseteq \BB(\CC\bG{t+}(\gamma), F)$ by
Corollary \ref{cor:x-depth},
we have by Lemma \ref{lem:Brgamma} that
$\gamma'_{\ge t} \in G_{x, t}$
and by Lemma \ref{lem:depth-in-center} that
$\gamma'_t \in G_{x, t}$\,.
Thus $\gamma'_{> t} \in G_{x, t}$\,.
By Corollary \ref{cor:good-comm}, there is $h \in G_{x, 0+}$
such that
$\lsup h g'' \in \gamma'_{\ge t}C_G(\gamma'_t)_{x, t+}$\,.
By Lemma \ref{lem:unipotent},
$\gamma'_{\ge t}C_G(\gamma'_t)_{x, t+}
\subseteq \gamma'_t C_G(\gamma'_t)_{t+}$\,.
In particular, $(\gamma'_t)$ is a normal
$(t{+})$-approximation to $\lsup h g''$; so, by Corollary
\ref{cor:compare-centralizers},
$Z(C_\bG(\gamma'_t)\conn) \subseteq Z(C_\bG(\lsup h g'')\conn)
\subseteq \lsup h\bG'$.
Since $\gamma'_t \in Z(C_G(\gamma'_t)\conn)$, we have
$\gamma'_t \in \lsup h G'$.
Thus
$\gamma'_t \in \lsup h G' \cap G_{x, t}
	= \lsup h G'_{x, t}
	\subseteq (G', G)_{x, (t, t{+})}$.
Since we saw above that
$\gamma'_{\ge t} \in G'G_{x, t{+}} \cap G_{x, t}
	= (G', G)_{x, (t, t{+})}$,
we have that
$\gamma'_{> t}
\in (G', G)_{x, (t, t{+})}
\subseteq (G' \cap \stab_G(\ox))G_{x, t{+}}$\,,
as desired.
\end{proof}

\begin{pn}
\label{prop:aniso-Levi}
Suppose that
\begin{itemize}
\item $r \in \tR_{> 0}$;
\item $(\bG', \bG)$ is a tame reductive $F$-sequence;
\item $x \in \BB(\bG', F)$;
and
\item
$x \in \BB_r(\gamma)$
or $\bG'/Z(\bG')$ is $F$-anisotropic.
\end{itemize}
Then $\gamma \in (G' \cap \stab_G(\ox))G_{x, r}$ if and only if
\begin{enumerate}
\item
\label{item:aniso-Levi-intersection}
$\ZZ G r(\gamma)\gamma \cap G_{x, r} \ne \emptyset$, and
\item
\label{item:aniso-Levi-in-G'}
$\lsup k\ZZ G r(\gamma) \subseteq G'$
for some $k \in [\gamma; x, r]$.
\end{enumerate}
\end{pn}

\begin{rk}
As in the proof of Lemma \ref{lem:dont-conj}, we could, if
desired, replace the condition that
$\lsup k\ZZ G r(\gamma) \subseteq G'$
for some $k \in [\gamma; x, r]$
by the condition that
$\ZZ G r(\gamma) \subseteq \lsup{[\gamma; x, r]}G'$;
but the stated form is more convenient for our purposes.
An appropriate modification of Proposition
\ref{prop:aniso-Levi} also holds for $r = 0$.
We do not state it, since the modified result is vacuous.
\end{rk}

\begin{proof}
For the `only if' direction, suppose that
$\gamma \in (G' \cap \stab_G(\ox))G_{x, r}$\,.
If $\bG'/Z(\bG')$ is $F$-anisotropic, then
Lemma \ref{lem:aniso-Brgamma} implies that
$x \in \BB_r(\gamma)$.
Therefore, we may assume $x \in \BB_r(\gamma)$.
By Definition \ref{defn:Brgamma},
(\ref{item:aniso-Levi-intersection})
holds.

By Lemma \ref{lem:Brgamma},
we have $\gamma \in \gamma_{< r}G_{x, r}$\,.
Thus
$\gamma_{< r}G_{x, r}
\cap (G' \cap \stab_G(\ox))
   \ne \emptyset$.
If $r < \infty$, then
the intersection is open in $G'$, hence contains a semisimple
element, say $\gamma'$.
If $r = \infty$, then $\gamma$ belongs to the intersection
and is semisimple, so we may take $\gamma' = \gamma$.
Then $(\ugamma, x)$ is also
an $r$-approximation to $\gamma'$.  By Lemma
\ref{lem:normal-conj}, there is
$k \in [\ugamma; x, r]$
such that $(\lsup k\gamma_i)_{0 \le i < r}$ is a
normal $r$-approximation to $\gamma'$.
By Corollary \ref{cor:compare-centralizers},
$\lsup k\ZZ G r(\gamma) = \ZZ G r(\gamma')
\subseteq Z(C_G(\gamma')\conn) \subseteq G'$.

For the `if' direction, suppose that
$x \in \BB_r(\gamma)$ and
$\ZZ G r(\gamma) \subseteq \lsup{[\gamma; x, r]}G'$.
By Definition \ref{defn:Brgamma},
$\ZZ G r(\gamma)\gamma \cap G_{x, r} \ne \emptyset$.
Therefore,
$\gamma \in (\lsup{[\gamma; x, r]}G')G_{x, r}$\,, so
there is $k \in [\gamma; x, r]$ with
$\gamma \in (\lsup kG')G_{x, r} = \lsup k(G' G_{x, r})$.
By Remark \ref{rem:approx-facts-commutator},
$\lsup{k\inv}\gamma \in \gamma G_{x, r}$\,.
Thus $\gamma \in G'G_{x, r}$\,.
Since $x \in \BB_r(\gamma)$, we have
by Lemma \ref{lem:Brgamma}
that $\gamma_{\ge r} \in G_{x, r} \subseteq \stab_G(\ox)$
and by Remark \ref{rem:approx-facts-in-stab}
that
$\gamma_{< r} \in \stab_G(\ox)$;
so in fact
$\gamma \in G' G_{x, r} \cap \stab_G(\ox)
	= (G' \cap \stab_G(\ox))G_{x, r}$\,.

For the `if' direction,
it remains to prove that, if
$\ZZ G r(\gamma)\gamma \cap G_{x, r} \ne \emptyset$,
$\ZZ G r(\gamma) \subseteq \lsup{[\gamma; x, r]}G'$, and
$\bG'/Z(\bG')$ is $F$-anisotropic, then
$x \in \BB_r(\gamma)$.
From Remark \ref{rem:approx-facts-in-center},
for each $0 \le i < r$,
there is $k_i \in [\gamma; x, r]$ so that
$\lsup{k_i}\gamma_i \in G'$.
In particular, since $\bG'/Z(\bG')$ is $F$-anisotropic
and $x \in \BB(\bG', F)$,
we have that
$\lsup{k_i}\gamma_i$, hence also $\gamma_i$, acts by
a translation on $x$.
Since $\gamma_i$ is bounded modulo center, in fact
it fixes the image \ox of $x$ in $\BB\red(\bG, F)$.
If $i > 0$, then
$\depth_x(\lsup{k_i}\gamma_i) = \depth(\lsup{k_i}\gamma_i)$,
so also $\depth_x(\gamma_i) = \depth(\gamma_i) = i$;
i.e., $\gamma_i \in G_{x, i}$\,.
If $r < \infty$, then,
by applying Lemma \ref{lem:x-depth} repeatedly, one sees
that $x$ belongs to $\BB(\CC\bG r(\gamma), F)$.
If $r = \infty$, then
$\BB(\CC\bG r(\gamma), F) = \BB(\CC\bG{r'}(\gamma), F)$
for $r' < \infty$ sufficiently large,
so again $x \in \BB(\CC\bG r(\gamma), F)$.
By Definition \ref{defn:Brgamma},
since we knew already that
$\ZZ G r(\gamma)\gamma \cap G_{x, r} \ne \emptyset$,
this means that $x \in \BB_r(\gamma)$.
\end{proof}

\begin{cor}
\label{cor:aniso-Levi}
With the notation and hypotheses of Proposition
\ref{prop:aniso-Levi},
$\gamma \in \lsup{G_{x, 0+}}
	\bigl((G' \cap \stab_G(\ox))G_{x, r}\bigr)$
if and only if
$\gamma_{< r} \in \lsup{G_{x, 0+}}G'$.
\end{cor}

\begin{proof}
Suppose that
$\gamma \in \lsup{G_{x, 0+}}
	\bigl((G' \cap \stab_G(\ox))G_{x, r}\bigr)$.
Upon replacing $\gamma$ by a $G_{x, 0+}$-conjugate, we may,
and hence do, assume that
$\gamma \in (G' \cap \stab_G(\ox))G_{x, r}$\,.
By Proposition \ref{prop:aniso-Levi},
$\gamma_{< r} \in \ZZ G r(\gamma) \subseteq \lsup{k\inv}G'
	\subseteq \lsup{G_{x, 0+}}G'$
for some $k \in [\gamma; x, r]$.

On the other hand, suppose that
$\gamma_{< r} \in \lsup{G_{x, 0+}}G'$.
Upon replacing $\gamma$ by a $G_{x, 0+}$-conjugate, we may,
and hence do, assume that $\gamma_{< r} \in G'$.
Then $\ZZ G r(\gamma) = Z(C_G(\gamma_{< r})\conn) \subseteq G'$.
Since $x \in \BB_r(\gamma)$
(either directly from the hypothesis,
or from Lemma \ref{lem:aniso-Brgamma} if
$\bG'/Z(\bG')$ is $F$-anisotropic), certainly
$\ZZ G r(\gamma)\gamma \cap G_{x, r} \ne \emptyset$
(by Definition \ref{defn:Brgamma}).
By another application of Proposition \ref{prop:aniso-Levi},
$\gamma \in \stab_{G'}(\ox)G_{x, r}$\,.
\end{proof}

\appendix

\section{Generalities on reductive groups}
\label{sec:general}

In this section,
$F$ is an arbitrary field and
\bG is a connected reductive $F$-group.
The following results might already be known, but we
were unable to find references.

\begin{lemma}
\label{lem:torus-quotient-product}
Let \bT be an $F$-torus,
and \bS a closed (not necessarily connected) $F$-subgroup of \bT.
Then there is an $F$-subtorus $\bT'\subseteq\bT$, split over
the splitting field of $\bT/\bS$,
such that $\bT=\bT'\bS$.
\end{lemma}

\begin{proof}
Note that $\bS\conn$ is a torus.
From Proposition 1.8 of \cite{borel-tits:reductive-groups},
there is an $F$-torus $\bT'\subseteq \bT$ such that
$\bT'\bS\conn=\bT$ and $\bT'\cap \bS\conn$ is finite.
Thus,
$\bT'\bS=\bT$ and $\bT'\cap \bS$ is finite.
Therefore $\bT/\bS = \bT'\bS/\bS \cong \bT'/(\bT'\cap\bS)$
is isogenous to $\bT'$,
so Corollaire 1.9(a)
of \cite{borel-tits:reductive-groups}
implies that the splitting fields of $\bT'$ and $\bT/\bS$
are the same.
\end{proof}

\begin{cor}
\label{cor:torus-quotient-product}
Suppose that
\begin{itemize}
\item \bS is an $F$-torus in \bG,
\item \mo N is a closed $F$-subgroup of $Z(\bG)$,
and
\item $\bS/(\mo N \cap \bS)$ is $F$-split.
\end{itemize}
Then
$C_\bG(\bS)$ is an $F$-Levi subgroup of \bG.
\end{cor}

\begin{proof}
From Lemma \ref{lem:torus-quotient-product}, there
is an $F$-split torus $\bS'$ such that
$\bS = \bS'(\mo N \cap \bS)$.  Then
$C_\bG(\bS) = C_\bG(\bS')$, which is an $F$-Levi subgroup
of \bG.
\end{proof}

\begin{lm}
\label{lem:split-in-center}
If \mo N is a closed $F$-subgroup of $Z(\bG)$ and
$\gamma \in Z(\bG)(\ol F)$ belongs to a torus which is $F$-split
modulo \mo N,
then $\gamma$ belongs to every maximal $F$-split modulo \mo N
torus.
\end{lm}

\begin{proof}
By Corollaire 1.9 of \cite{borel-tits:reductive-groups}, a torus is
$F$-split modulo \mo N if and only if it is $F$-split modulo
$\mo N\conn$.  Thus we may, and hence do, assume that \mo N
is connected.
By %Corollary 3.18(ii), Proposition 6.7, and Theorem 6.8 of
Proposition 22.4 of
\cite{borel:linear}, the map $\bG \to \bG/\mo N$ is central.
Let \bT be a maximal $F$-split modulo \mo N torus such that
$\gamma \in \bT(\ol F)$,
and let $\bT'$ be any other maximal $F$-split modulo \mo N torus.
In particular, by maximality, \bT and $\bT'$ contain \mo N.
By Theorem 22.6 of \cite{borel:linear},
the images $\wtilde\bT$ and $\wtilde\bT'$
of \bT and $\bT'$, respectively,
in $\bG/\mo N$ are maximal $F$-split tori there.
Thus, by Theorem 20.9 of \cite{borel:linear},
there is $\ol g \in (\bG/\mo N)(\ol F)$
such that $\lsup{\ol g}\wtilde\bT = \wtilde\bT'$.
(In fact, \ol g may be chosen in $(\bG/\mo N)(F)$.)
This means that
$\ogamma = \lsup{\ol g}\ogamma \in \wtilde\bT'(\ol F)$,
where \ogamma is the image in $\bG/\mo N$ of $\gamma$.
Thus $\gamma$ belongs to $(\mo N\bT')(\ol F) = \bT'(\ol F)$,
as desired.
\end{proof}

\begin{defn}
\label{defn:root-value}
Let $\gamma\in G$ be semisimple.
Then the \emph{root values of $\gamma$ for $\bG$} are the elements of
\indexme{root values of a semisimple element}
$\set{\alpha(\gamma)}{\alpha\in\Phi(\bG,\bT)}$,
and the \emph{character values of $\gamma$} are the elements of
\indexme{character values of a semisimple element}
$\set{\chi(\gamma)}{\chi\in\bX^*(\bT)}$,
where \bT is any maximal torus in $C_\bG(\gamma)$.
\end{defn}

Note that, by \cite{borel:linear}*{Corollary 8.5},
the character values of $\gamma$ do not depend on the
choice of $\bT$ or $\bG$, and thus can be defined even
when $\bG$ is non-connected and $\gamma \notin G\conn$.
We will not pursue this here.

\begin{dn}
\label{defn:bad-prime}
Recall
from \cite{springer-steinberg:conj}*{\S I.4.1}
that a prime $p$ is \emph{bad} for a root system $\Phi$
if there is some
(integrally) closed subsystem $\Phi_1\subseteq\Phi$
such that $\Z\Phi/\Z\Phi_1$ has $p$-torsion;
and that $p$ is \emph{bad} for $\bG$
if it is bad for 
the absolute root system for $\bG$ with
respect to some (hence any) maximal torus.
\end{dn}

The sets of bad primes for the irreducible root systems
are given below.
$$
\begin{array}{l l}
\Phi  & \text{bad primes for $\Phi$}  \\
\cline{1-2}
\mathsf{A}_n			 &\emptyset	 \\
\mathsf{B}_n, \mathsf{C}_n, \mathsf{D}_n &	 \sset 2	 \\
\mathsf{E}_6,\mathsf{E}_7,\mathsf{F}_4,\mathsf{G}_2 &	 \sset{2,3} \\
\mathsf{E}_8			 &\sset{2,3,5}
\end{array}
$$
In general, a prime is bad for a root system $\Phi$ if and only if
it is bad for some irreducible factor of $\Phi$.

\begin{defn}
\label{defn:good-root-values}
\indexme{good root values}
We say that an element of $F\cross$ is \emph{bad} for \bG
if it is a non-trivial root of unity whose order has only
bad primes for \bG as prime divisors.
If $\gamma \in G$ is semisimple, and no root value of
$\gamma$ is bad for \bG, then we say that $\gamma$
\emph{has only good root values} for \bG.
\end{defn}

\begin{prop}
\label{prop:levi}
Suppose that $\gamma \in G$ is semisimple and has
only good root values for $\bG$.
Then $\bL := C_\bG(\gamma)\conn$ is a Levi subgroup of \bG.

Moreover, let $E$ denote the minimal Galois extension of $F$
such that \bL is an $E$-Levi subgroup of \bG,
and $E'$ the extension of $F$ generated by
the root values of $\gamma$ for $\bG$.
Then $E$ contains $E'$.
Suppose that no root value of $\gamma$ for $\bG$
is a non-trivial root of unity.
Then $E=E'$.
\end{prop}

The proof is adapted and generalized from part of the proof
of Theorem 4.14 of \cite{roche:thesis}.
An early draft was discussed with Jonathan Korman.

\begin{proof}
Let \bT be a maximal $F$-torus in \bG such that $\gamma\in T$.
To prove the first claim,
it is enough to construct an $F$-torus
$\bS\subseteq\bT$ such that $\bL= C_{\bG}(\bS)$.

For a subset $\bY\subseteq \bX^*(\bT)$,
let $\bY^\perp \subseteq \bT$ denote the intersection of the
kernels of the elements of $\bY$.
For a subgroup $\bT'\subseteq\bT$, let $\smash{\bT'}^\perp$
denote the lattice of characters in $\bX^*(\bT)$
that are trivial on $\bT'$.

Let $\Phi=\Phi(\bG,\bT)$,
$\Phi_\gamma = \set{\alpha\in\Phi}{\alpha(\gamma)=1}$,
$\wtilde\bS = \Phi_\gamma^\perp$,
and $\bS = \wtilde\bS^\circ$.
Since (by \cite{springer-steinberg:conj}*{\S II.4.1(b)})
$$
C_{\bG}(\bS) = \langle \bT,\bU_\beta :
\beta\in\bS^\perp\cap\Phi\rangle
$$
and
$$
\bL = \langle\bT, \bU_\alpha : \alpha\in\Phi_\gamma\rangle,
$$
it will be enough to show that $\bS^\perp \cap \Phi = \Phi_\gamma$.

We show that $\bS^\perp \cap \Phi \subseteq  \Phi_\gamma$,
since the opposite containment is trivial.
Put $m = \smcard{\wtilde\bS/\bS}$.
For $\alpha\in \bS^\perp\cap \Phi$,
we have
$m\alpha\in\wtilde\bS^\perp
= \Phi_\gamma^{\perp\perp}$,
which contains $\Z\Phi_\gamma$ with finite index,
so that some positive integer multiple of $\alpha$ lies in $\Z\Phi_\gamma$.
Let $n$ be the smallest positive integer such that $n\alpha \in \Z\Phi_\gamma$.
Since $\Phi_\gamma$ is closed in $\Phi$,
$n$ is by definition a product of bad primes.
Since $(n\alpha)(\gamma) = 1$, we have that
$\alpha(\gamma)$ is an $n$th root of unity.
But $\gamma$ has only good root values,
so this implies that $\alpha(\gamma)=1$,
and thus that $\alpha\in\Phi_\gamma$.

It is clear that $E$ must contain $E'$.

We now prove the final statement of the proposition.
Let $\bY_\gamma$ denote the lattice of all characters in $\Z\Phi$
that are trivial on $\gamma$.
Put $\wtilde\bS'= \bY_\gamma^\perp$,
and $\bS' = (\wtilde\bS')^\circ$.
Then $\gamma\in\wtilde\bS'$.
Pick a positive integer $n$ such that $\gamma^n\in\bS'$.
Our hypothesis on the root values of $\gamma$ implies that
for $\alpha\in\Phi$, $\alpha(\gamma^n)=1$ if and only if
$\alpha(\gamma)=1$.
Thus, by another application of
\cite{springer-steinberg:conj}*{\S II.4.1(b)},
$C_{\bG}(\gamma)\conn = C_{\bG}(\gamma^n)\conn$.
Since $\bS' \subseteq \bS$,
we have that
$\bL = C_{\bG}(\gamma^n)^\circ
\supseteq C_{\bG}(\bS') \supseteq C_{\bG}(\bS) = \bL$,
and thus that
$\bL =C_{\bG}(\bS')$.
Therefore, it is enough to show that $\bS'/(Z(\bG)\conn \cap \bS')$
splits over $E'$.

In order to do so,
we must show that for a Galois splitting field $L$ of \bT over $E'$,
the group $\Gamma=\Gal(L/E')$ acts trivially
on the lattice generated by the restrictions to $\bS'$ of roots
in $\Phi$.
That is, we must show that $\Gamma$ acts trivially on
$\Z\Phi/(\Z\Phi \cap \smash{\bS'}^\perp)$.
Let $\sigma\in\Gamma$ and $\alpha\in\Phi$.
By the construction of $E'$, $\sigma(\alpha) - \alpha$
belongs to the lattice $\bY_\gamma$.
Since
$\bY_\gamma \subseteq \smash{\wtilde\bS'}^\perp
	\subseteq \smash{\bS'}^\perp$,
we are finished.
\end{proof}

\begin{lemma}
\label{lem:power-no-roots-unity}
Suppose that $\gamma \in G$ is semisimple.
Then there is a positive integer $n$,
not divisible by $\chr F$, 
such that
no character value of $\gamma^n$ is a non-trivial root of unity.
\end{lemma}

\begin{proof}
Let \bT be an $F$-torus such that $\gamma \in T$,
and put
$X=\bX^*(\bT)$,
$Y=\set{\chi\in X}{\chi(\gamma) = 1}$,
and
$n= \card{(X/Y)\textsub{tor}}$.
Then no character value of $\gamma^n$ is a non-trivial root
of unity.

It only remains to show that $n$ is not divisible by $p := \chr F$.
Indeed, suppose that $\chi \in X$ and $p\chi \in Y$.
Then $\chi(\gamma)^p = 1$.
Since $\ol F$ has no non-trivial $p$th roots of unity,
$\chi(\gamma) = 1$, so
$\chi \in Y$.
That is, $X/Y$ has no $p$-torsion, so $p$ does not divide
$n$.
\end{proof}

\begin{lemma}
\label{lem:no-roots-unity-torus}
Suppose that \bT is a torus, and \mc S is any subset of \bT
such that no character value of an element of $\mc S$ is a
non-trivial root of unity.
Then the Zariski closure of the group
generated by $\mc S$ is a torus.
\end{lemma}

\begin{proof}
Let $X = \bX^*(\bT)$ and
$Y = \sett{\chi\in X}{$\chi(\gamma)=1$ for all $\gamma\in \mc S$}$.
By assumption, $X/Y$ is torsion free,
so there is some lattice $X_0\subseteq X$
such that $X = X_0 \oplus Y$.
Let $\bS$ denote the Zariski closure of the group generated
by $\mc S$.
By Corollary 8.5 and Proposition 8.7 of \cite{borel:linear},
there is an $F$-torus $\bS'\subseteq\bT$ such that
$\bT = \bS\conn \times \bS'$.
Fix $\gamma \in \mc S$, and write
$\gamma = \gamma_0\gamma'$, with
$\gamma_0\in \bS\conn$
and
$\gamma' \in \bS'$.
Every character $\chi\in \bX^*(\bS')$
has a unique extension (which we will also denote $\chi$)
to \bT that is trivial on $\bS\conn$.
Let $n = \card{\bS/\bS\conn}$,
so $\chi(\gamma)^n=1$.
By our assumption on $\mc S$,
$\chi(\gamma)=1$.
Therefore, $\chi(\gamma') = \chi(\gamma)\chi(\gamma_0\inv) = 1$.
Since this is true for all $\chi\in \bX^*(\bS')$,
we have that $\gamma'=1$;
that is, $\gamma = \gamma_0 \in \bS\conn$.
Thus $\bS = \bS\conn$; that is, \bS is connected.
By Proposition 8.5 of \cite{borel:linear}, it is a torus.
\end{proof}

\begin{cor}
\label{cor:no-roots-unity-levi}
Suppose that $\gamma \in G$ is semisimple, and that
no character value of $\gamma$ is a non-trivial root of unity.
Then $C_\bG(\gamma)$ is an $E$-Levi subgroup of $\bG$,
where $E$ is the extension of $F$ generated by the root values
of $\gamma$.
\end{cor}

\begin{proof}
Let $\bS$ denote the Zariski closure of the group generated
by $\gamma$.  Then \bS is a torus by Lemma
\ref{lem:no-roots-unity-torus}, so
$C_\bG(\gamma) = C_\bG(\bS)$ is a Levi subgroup of \bG.
By Proposition \ref{prop:levi},
it is actually an $E$-Levi subgroup.
\end{proof}

\begin{lm}
\label{lem:still-split-in-levi}
If \bM is an $F$-Levi subgroup of \bG,
\bT is an $F$-split torus in \bG,
and \mo D is a subvariety of $\bM \cap \bT$ defined over
$F$,
then \mo D is contained in an $F$-split torus in
\bM.
\end{lm}

Of course, we may take \mo D to be $\bM \cap \bT$ if that
variety is defined over $F$ (in particular, if $F$ is
perfect).

\begin{proof}
Put $\bH = C_\bG(\mo D)\conn$,
which is defined over $F$ by Corollary 9.2 (and Proposition
1.2) of \cite{borel:linear}
and reductive by \cite{springer-steinberg:conj}*{\S II.4.1(b)}.
Let \bS be the $F$-split part of the center of \bM, so
that $\bM = C_\bG(\bS)$.
Then \bS and \bT are $F$-split tori in \bH.
Upon enlarging \bT if necessary, we may, and hence do, assume that
it is a maximal $F$-split torus in \bH.
Thus, there is $h \in H$ such that
$\bS \subseteq \lsup h\bT$.
Then $\lsup h\bT$ is an $F$-split torus that is
contained in $C_\bG(\bS) = \bM$
and contains
$\lsup h\mo D = \mo D$.
\end{proof}

For the remainder of this appendix, suppose that
\bG is $F$-quasisplit, and
\bH is a reductive $F$-subgroup of \bG of the same
absolute and $F$-split ranks as \bG.

\begin{lm}
\label{lem:rank-ascent}
We have that \bH is $F$-quasisplit,
and that
\bH and \bG have the same $E$-split rank for every
extension $E/F$.
\end{lm}

\begin{proof}
Let \bS be an $F$-split torus maximal in \bH, hence in \bG.
Since \bG is $F$-quasisplit, the centralizer
$\bT := C_\bG(\bS)$ is a maximal torus in \bG.
Clearly $C_\bH(\bS)$ is contained in \bT and contains a
maximal torus in \bH.
Since \bH and \bG have the same absolute rank,
$\bT = C_\bH(\bS)$.
This implies that \bH is $F$-quasisplit.

Now let $E/F$ be an extension.  Then
\bS is contained in a maximal $E$-split torus in \bG.  This
torus is necessarily contained in $C_\bG(\bS) = \bT$, hence
in \bH.  Thus, \bG and \bH have the same $E$-split rank.
\end{proof}

\begin{lm}
\label{lem:no-root-ambience}
Let \bS be a maximal $F$-split torus in \bH, and fix
$a \in \Phi(\bH, \bS)$.
If $2a \in \Phi(\bG, \bS)$, then $2a \in \Phi(\bH, \bS)$.
The unique closed connected subgroup of \bH with Lie
algebra the sum of the $a$- and $2a$-weight spaces in
$\Lie(\bH)$ for the action of \bS
equals
the unique closed connected subgroup of \bG with Lie algebra
the sum of the $a$- and $2a$-weight spaces in $\Lie(\bG)$
for the action of \bS.
\end{lm}

\begin{proof}
Put $\bT = C_\bG(\bS)$, a maximal torus in \bG.  Since \bH
has the same absolute rank as \bG, also $\bT = C_\bH(\bS)$,
so \bH is also $F$-quasisplit.

We will use throughout this proof the following fact:
If $\alpha \in \Phi(\bH, \bT)$, then, since the
$\alpha$-weight space in $\Lie(\bG)$ is one-dimensional
by Theorem 13.18(4b) of \cite{borel:linear},
the $\alpha$-weight spaces in $\Lie(\bH)$ and $\Lie(\bG)$
are equal.

Choose a character $\alpha$ of \bT which restricts to $a$
such that the $\alpha$-weight space in $\Lie(\bH)$ is non-trivial.
Since \bG is $F$-quasisplit,
by Proposition 15.5.3(iii) of \cite{springer:lag},
any such character is of the
form $\sigma\alpha$ for some $\sigma \in \Gal(F\sep/F)$.
Since $\Lie(\bH)$ is $\Gal(F\sep/F)$-stable, for any
$\sigma \in \Gal(F\sep/F)$, the $\sigma\alpha$-weight space
in $\Lie(\bH)$ is also non-trivial, hence equal to the
$\sigma\alpha$-weight space in $\Lie(\bG)$.
Since the $a$-weight space in $\Lie(\bH)$, respectively
$\Lie(\bG)$, is equal to the sum of the
$\sigma\alpha$-weight spaces in $\Lie(\bH)$, respectively
$\Lie(\bG)$, as $\sigma$ runs over $\Gal(F\sep/F)$, in fact
the $a$-weight spaces in $\Lie(\bH)$ and $\Lie(\bG)$ are the
same.

If $2a \not\in \Phi(\bG, \bS)$, then the $2a$-weight spaces
in $\Lie(\bH)$ and $\Lie(\bG)$ are both trivial.  Otherwise,
by another application of Proposition 15.5.3(iii) of
\cite{springer:lag}, there is an element
$\sigma \in \Gal(F\sep/F)$ such that 
$\alpha + \sigma\alpha \in \Phi(\bG, \bT)$.
By \cite{springer:lag}*{\S 17}, the irreducible factor of
$\Phi(\bG, \bT)$ containing $\alpha$ and $\sigma\alpha$
is of type $\ms A_{2n}$.  By
\cite{borel-tits:homomorphismes}*{\S 4.3}, the $\alpha$- and
$\sigma\alpha$-weight spaces in $\Lie(\bG)$ do not commute.
Since both weight spaces lie in $\Lie(\bH)$, so does their
commutator.  Since their commutator lies in (in fact, equals) the
$(\alpha + \sigma\alpha)$-weight space of $\Lie(\bG)$, we have
that the $(\alpha + \sigma\alpha)$-weight space of $\Lie(\bH)$
is non-trivial.  Thus the $2a$-weight space of $\Lie(\bH)$ is
non-trivial; i.e., $2a \in \Phi(\bH, \bS)$.
This is the first part of the lemma.
Now an argument similar to the above shows
that, actually, the $2a$-weight spaces of $\Lie(\bH)$ and
$\Lie(\bG)$ are equal.

Thus the two Lie algebras mentioned in the statement of the
lemma are equal, so the groups mentioned there are also
equal.  This is the second part of the lemma.
\end{proof}

The following geometric result is well known in case \bH is
a maximal $F$-torus in \bG.

\begin{lm}
\label{lem:Zariski}
Suppose that
\begin{itemize}
\item \bH is connected;
\item \bS is a maximal $F$-split torus in \bH;
and
\item $h \in \bH(\ol F)$ and, for each
$a \in \Phi(\bG, \bS) \smallsetminus \Phi(\bH, \bS)$
with $a/2 \not\in \Phi(\bG, \bS)$,
$u_a \in \bU_a(\ol F)$.
\end{itemize}
Denote by $\mult_\bH$ the multiplication map
$$
\bH \times
\prod_{\substack{a \in \Phi(\bG, \bS) \smallsetminus \Phi(\bH, \bS) \\
		a/2 \not\in \Phi(\bG, \bS)}}
	\bU_a
\to \bG
$$
(the latter product taken in any order).
Then $\mult_\bH$ is injective, and
$\mult_\bH(h, (u_a)) \in G$ if and only if $h \in H$
and $u_a \in U_a$ for
$a \in \Phi(\bG, \bS) \smallsetminus \Phi(\bH, \bS)$.
\end{lm}

\begin{proof}
Put $\bT = C_\bG(\bS)$, a maximal torus in \bG;
and abbreviate $\mult_\bH$ to \mult.
By Proposition 14.4(2a) and Corollary 14.14 of
\cite{borel:linear},
the multiplication maps
$\bT \times \prod_{\beta \in \Phi(\mo J, \bT)} \bU_\beta
\to \mo J$
are open immersions for any connected reductive subgroup \mo J
of \bG containing \bT.
By Lemma \ref{lem:no-root-ambience},
if $\alpha \in \Phi(\bH, \bT)$,
then the root subgroups of \bG and \bH associated to
$\alpha$ are the same, so we may write $\bU_\alpha$ without
ambiguity.
By Proposition 14.5 and Remark 21.10(1) of
\cite{borel:linear},
for $a \in \Phi(\bG, \bS)$,
$\bU_a = \prod_{\alpha\bigr|_\bS \in \sset{a, 2a}} \bU_\alpha$
(as varieties).
Fix, temporarily, $\alpha \in \Phi(\bG, \bT)$, and put
$a := \alpha\bigr|_\bS$.
If $a \in \Phi(\bH, \bS)$, then $\bU_\alpha \subseteq \bH$ is
contained in the domain of \mult.
If $a \not\in \Phi(\bH, \bS)$ and $a/2 \not\in \Phi(\bG, \bS)$,
then $\bU_\alpha \subseteq \bU_a$ is obviously contained in
the domain of \mult.
If $a \not\in \Phi(\bH, \bS)$ but $a/2 \in \Phi(\bG, \bS)$,
then, by Lemma \ref{lem:no-root-ambience}, also
$a/2 \not\in \Phi(\bH, \bS)$.  Then
$\bU_\alpha \subseteq \bU_{a/2}$ is contained in the domain
of \mult.
Thus there is an open subset
$\mo V :=
\bT \times \prod_{\alpha \in \Phi(\bG, \bT)} \bU_\alpha$
of the domain of \mult on which \mult is an open immersion.
Moreover, the orbits of \mo V under the action
of \bH fill out the domain of \mult.
Since \mult is \bH-equivariant,
we have that the differential of \mult at every point of its
domain is an isomorphism.
By Zariski's main theorem (see
Corollaire 18.12.13 of \cite{grothendieck:EGA-IV}),
\mult is an open immersion; in particular, it is injective.

Denote by \mo C the image of \mult.
(For example, if
$\bH = \bT$, then
\mo C is a translate of the `big cell'.)
Then
\mult is an isomorphism onto \mo C,
so $\mo C(F)$ is the image
under $\mult$ of the set of $F$-rational points
in the domain of $\mult$, namely
$H \times
\prod_{\substack{a \in \Phi(\bG, \bS) \smallsetminus \Phi(\bH, \bS) \\
		a/2 \not\in \Phi(\bG, \bS)}}
	U_a$.
The remainder of the result follows from the injectivity of
\mult.
\end{proof}

\begin{lm}
\label{lem:rel-root-red}
Suppose that
$E/F$ is a finite Galois extension of odd degree.
If \bG splits over $E$, then $\lsub F\Phi(\bG)$ is reduced.
\end{lm}

\begin{proof}
Let \bS be a maximal $F$-split torus in \bG and
put $\bT=C_\bG(\bS)$, a maximal $F$-torus in \bG.
If $\Phi(\bG, \bS)$ is not reduced, then,
by \cite{springer:lag}*{\S 17},
$\Phi(\bG,\bT)$ has an irreducible
factor of type $\ms A_{2n}$ for some $n$
and $\Gal(E/F)$ acts non-trivially on the Dynkin diagram
of $\ms A_{2n}$.
This implies that $E/F$ has even degree, which is a contradiction.
\end{proof}

\section{Generalities on concave functions}
\label{sec:concave-fns}

Let $\Phi$ be a root system in some \Q-vector space $V$.
Put $\wtilde\Phi = \Phi \cup \sset 0$.

\begin{dn}
\label{defn:concave-vee}
\indexmem{f_1 \vee f_2}
Given $C \in \tR$ and functions
$f_j : \wtilde\Phi \to \tR$ for $j = 1, 2$,
we define the function
$(f_1 \vee f_2)_C : \wtilde\Phi \to \tR \cup \sset{-\infty}$
by
$$
(f_1 \vee f_2)_C(\alpha)
= \inf \sset{\sum_m f_1(a_m) + \sum_n f_2(b_n)}
\quad\text{for $\alpha \in \wtilde\Phi$,}
$$
the infimum (in \tR) taken over all pairs of non-empty finite
sequences $(a_m)_m$ and $(b_n)_n$ in $\wtilde\Phi$,
both of length at most $C$,
such that $\sum_m a_m + \sum_n b_n = \alpha$.
Put $f_1 \vee f_2 = (f_1 \vee f_2)_\infty$.
\end{dn}

It is straightforward to verify that, for any functions
$f_1$ and $f_2$ as in Definition \ref{defn:concave-vee},
$f_1 \vee f_2$ is concave (in the sense of Definition
\ref{defn:concave}) as long as
$(f_1 \vee f_2)(\alpha) \ne -\infty$ for all
$\alpha \in \wtilde\Phi$
(in particular, by the next result, as long as
$(f_1 \vee f_2)(0) \ne -\infty$%
	%%We have that
	%%$(f_1 \vee f_2)(0) \le (f_1 \vee f_1)(0) + f_2(0)$,
	%%with a similar bound in terms of $f_2 \vee f_2$.
	%%Thus, $(f_1 \vee f_2)(0) \ne -\infty$ means that
	%%the hypotheses of the next result are satisfied.
);
and that a function $f : \wtilde\Phi \to \tR$ is
concave if and only if $f \vee f \ge f$.

\begin{lm}
\label{lem:concave-vee-achieved}
There is a constant $C(\Phi)$, depending only on $\Phi$, with
the following property.
Fix functions
$f_j : \wtilde\Phi \to \tR$
for $j = 1, 2$.
If
$(f_1 \vee f_1)(0) \ne -\infty \ne (f_2 \vee f_2)(0)$
(in particular, if $f_1$ and $f_2$ are non-negative),
then $f_1 \vee f_2 = (f_1 \vee f_2)_{C(\Phi)}$.
\end{lm}

\begin{proof}
Note that, in fact,
$(f_j \vee f_j)(0) \ge 0$, so
$f_j(0) \ge \frac 1 2(f_j \vee f_j)(0) \ge 0$,
for $j = 1, 2$.
Put
$C(0) = 2\card{\smash{\wtilde\Phi}}$
and
$C(\alpha) = C(\Phi) = 3\card{\smash{\wtilde\Phi}}$
for $\alpha \in \Phi$.

Fix $\alpha \in \wtilde\Phi$, and put
$C = C(\alpha)$.
We claim that there does not exist a pair of non-empty
finite sequences
$(a_m)_m$ and $(b_n)_n$ in $\wtilde\Phi$
such that
$\sum_m a_m + \sum_n b_n = \alpha$
but
$$
\sum_m f_1(a_m) + \sum_n f_2(b_n)
< (f_1 \vee f_2)_C(\alpha).
$$
Indeed, suppose there does exist such a pair,
and let $((a_m)_m, (b_n)_n)$ be one of
minimal total length, say $\ell$.

Let $(c_p)_{0 \le p < \ell}$ be a sequence containing
precisely the terms of $(a_m)_m$ and $(b_n)_n$, counted with
multiplicities.
Recall that
$\alpha = \sum_{p = 0}^{\ell - 1} c_p \in \wtilde\Phi$.
We claim that there is a permutation $\pi$ of
$\sset{0, 1, \dotsc, \ell - 1}$ such that
$c^{(i)} := \sum_{p = 0}^{i - 1} c_{\pi(p)} \in \wtilde\Phi$
for all $0 < i \le \ell$.
By induction, it suffices to show that there is some
index $0 \le p_0 < \ell$ such that
$\alpha - c_{p_0} \in \wtilde\Phi$.
If $\alpha = 0$, then any index $p_0$ will do.
If $\alpha \ne 0$, then, since
$$
\Bigl\langle\alpha, \sum_{p = 0}^{\ell - 1} c_p\Bigr\rangle
= \langle\alpha, \alpha\rangle > 0,
$$
there is some index $p_0$ such that
$\langle\alpha, c_{p_0}\rangle > 0$.
Then, by the corollary to Theorem VI.1.3.1 of
\cite{bourbaki:lie-gp+lie-alg_4-6},
$\alpha - c_{p_0} \in \wtilde\Phi$.

Since, necessarily,
$\ell > \card{\smash{\wtilde\Phi}}$,
there are indices $0 < i < i' \le \ell$ such that
$c^{(i)} = c^{(i')}$, hence such that
$\sum_{i \le p < i'} c_{\pi(p)} = 0$.  That is, there is some
proper subsequence of $(a_m)_m \cup (b_n)_n$ the sum of
whose terms (with multiplicities) is $0$.
Let $(\gamma_q)_q$ be of maximal length among such
subsequences.
There are four cases.
\begin{enumerate}
\item\label{case:c-swallowed-in-a}
The terms of $(\gamma_q)_q$ form a proper
subsequence of the terms of $(a_m)_m$, or a proper
subsequence of the terms of $(b_n)_n$ (counted with
multiplicities).
\item\label{case:c-swallows-some-a-some-b}
The terms of $(\gamma_q)_q$ include some, but not
all, of the terms of $(a_m)_m$, and some, but not all, of
the terms of $(b_n)_n$ (counted with multiplicites).
\item\label{case:c-swallows-a}
The terms of $(\gamma_q)_q$ are precisely the terms
of $(a_m)_m$, or precisely the terms of $(b_n)_n$ (counted
with multiplicities).
\item\label{case:c-swallows-a-and-some-b}
The terms of $(\gamma_q)_q$ include all of the terms
of $(a_m)_m$ and some (but, necessarily, not all) of the
terms of $(b_n)_n$ (counted with multiplicites), or
\emph{vice versa}.
\end{enumerate}

In case \eqref{case:c-swallowed-in-a}, if the first
possibility obtains,
let $(a_{m'})_{m'}$ be the
subsequence of $(a_m)_m$ containing the terms of
$(\gamma_q)_q$, and $(a_{m''})_{m''}$ the
complementary subsequence.
Then
\begin{multline*}
(f_1 \vee f_2)_C(\alpha)
> \sum_{m'} f_1(a_{m'}) + \sum_{m''} f_1(a_{m''})
	+ \sum_n f_2(b_n) \\
\ge \min\sset{f_1(0), (f_1 \vee f_1)(0)}
	+ \sum_{m''} f_1(a_{m''}) + \sum_n f_2(b_n).
\end{multline*}
Since $\min\sset{f_1(0), (f_1 \vee f_1)(0)} \ge 0$, we have
$$
(f_1 \vee f_2)_C(\alpha)
> \sum_{m''} f_1(a_{m''}) + \sum_n f_2(b_n),
$$
contradicting the minimality of $((a_m)_m, (b_n)_n)$.
The second possibility is handled similarly.

In case \eqref{case:c-swallows-some-a-some-b},
let $(a_{m'})_{m'}$ be the subsequence of $(a_m)_m$
containing the terms it shares with $(\gamma_q)_q$,
and $(a_{m''})_{m''}$ the complementary subsequence.
Define $(b_{n'})_{n'}$ and $(b_{n''})_{n''}$ similarly.
Then
\begin{multline*}
(f_1 \vee f_2)_C(\alpha)
> \sum_{m'} f_1(a_{m'}) + \sum_{n'} f_2(b_{n'})
	+ \sum_{m''} f_1(a_{m''}) + \sum_{n''} f_2(b_{n''}) \\
\ge (f_1 \vee f_2)(0)
	+ \sum_{m''} f_1(a_{m''}) + \sum_{n''} f_2(b_{n''}).
\end{multline*}
Since $(f_1 \vee f_2)(0) \ge 0$, we obtain a contradiction
as before.

Suppose we are in case \eqref{case:c-swallows-a}.
Then $\sum_m a_m, \sum_n b_n \in \wtilde\Phi$.
Necessarily, $(a_m)_m$ or $(b_n)_n$ has length greater than
$\card{\smash{\wtilde\Phi}}$; say $(a_m)_m$ does.
(The other possibility is treated similarly.)
By reasoning as above, we see that there is some proper
subsequence $(a_{m'})_{m'}$ of $(a_m)_m$ such that
$\sum_{m'} a_{m'} = 0$.
Let $(a_{m''})_{m''}$ be the complementary subsequence.
Then
\begin{multline*}
(f_1 \vee f_2)_C(\alpha)
> \sum_{m'} f_1(a_{m'})
	+ \sum_{m''} f_1(a_{m''}) + \sum_n f_2(b_n) \\
\ge \min\sset{f_1(0), (f_1 \vee f_1)(0)}
	+ \sum_{m''} f_1(a_{m''}) + \sum_n f_2(b_n).
\end{multline*}
Since $\min\sset{f_1(0), (f_1 \vee f_1)(0)} \ge 0$,
we obtain a contradiction as before.

In case \eqref{case:c-swallows-a-and-some-b}, if the first
possibility obtains, let $(b_{n'})_{n'}$ be the
subsequence of $(b_n)_n$ containing the terms it shares
with $(\gamma_q)_q$, and $(b_{n''})_{n''}$ the
complementary subsequence.
If $\alpha = 0$, then
\begin{multline*}
(f_1 \vee f_2)_C(\alpha)
> \sum_m f_1(a_m) + \sum_{n'} f_2(b_{n'})
	+ \sum_{n''} f_2(b_{n''}) \\
\ge \sum_m f_1(a_m) + \sum_{n'} f_2(b_{n'})
	+ \min\sset{f_2(0), (f_2 \vee f_2)(0)}.
\end{multline*}
Since $\min\sset{f_2(0), (f_2 \vee f_2)(0)} \ge 0$,
we obtain a contradiction in the usual fashion.
The second possibility is handled similarly.
Thus we have shown that
$(f_1 \vee f_2)_{C(0)}(0) = (f_1 \vee f_2)(0)$.

Now suppose $\alpha \ne 0$.
If $(b_{n''})_{n''}$ had more than
$\card{\smash{\wtilde\Phi}}$ terms, then the usual
argument would show that $(\gamma_q)_q$ was not maximal,
which is a contradiction.
Thus
\begin{multline*}
(f_1 \vee f_2)_C(\alpha)
> \sum_m f_1(a_m) + \sum_{n'} f_2(b_{n'})
	+ \sum_{n''} f_2(b_{n''}) \\
\ge (f_1 \vee f_2)(0) + \sum_{n''} f_2(b_{n''})
= (f_1 \vee f_2)_{C(0)}(0) + \sum_{n''} f_2(b_{n''}) \\
\ge (f_1 \vee f_2)_{C(0) + \card{\smash{\wtilde\Phi}}}
	(\alpha)
= (f_1 \vee f_2)_C(\alpha),
\end{multline*}
a contradiction.  Again, the second possibility is handled
similarly.

We have shown that
$(f_1 \vee f_2)(\alpha)
\ge (f_1 \vee f_2)_C(\alpha)
\ge (f_1 \vee f_2)_{C(\Phi)}(\alpha)$.
The reverse inequality is obvious.
Since $\alpha \in \wtilde\Phi$ was arbitrary, we have the
desired equality of functions.
\end{proof}

\begin{lm}
\label{lem:concave-in-R}
For any function
$h : \wtilde\Phi \to \tR$ and any
$\delta \in \R_{> 0}$,
put
\indexmem{h_\delta}
$$
h_\delta(\alpha) =
\begin{cases}
h(\alpha),  & h(\alpha) \in \R \cup \sset\infty, \\
r + \delta, & h(\alpha) = r{+}, r \in \R
\end{cases}
$$
for $\alpha \in \wtilde\Phi$.
Fix functions
$f_j : \wtilde\Phi \to \tR$
for $j = 1, 2, 3$
such that $f_1 \vee f_2 \ge f_3$.
Then
$f_{1, \delta} \vee f_{2, \delta} \ge f_{3, \delta}$
for $\delta \in \R_{> 0}$ sufficiently small.
If also $(f_1 \vee f_2)(\alpha) > f_3(\alpha)$
whenever $f_3(\alpha) \ne \infty$, then
$f_{1, \delta} \vee f_{2, \delta}
\ge f_{3, \delta} + \varepsilon$
for $\delta, \varepsilon \in \R_{> 0}$ sufficiently small.
\end{lm}

\begin{proof}
Consider the set of all $r - s$,
where $r, s \in \R$ with $s < r$
are such that there exist
non-empty finite sequences
$(a_m)_m$ and $(b_n)_n$ in $\wtilde\Phi$,
each of total length no more than $C(\Phi)$,
such that
\begin{itemize}
\item
$\sum_m f_1(a_m) + \sum_n f_2(b_n) \in \sset{r, r{+}}$;
\item
$\alpha := \sum_m a_m + \sum_n b_n \in \wtilde\Phi$;
and
\item
$f_3(\alpha) \in \sset{s, s{+}}$.
\end{itemize}
(Here, $C(\Phi)$ is as in Lemma
\ref{lem:concave-vee-achieved}.)
This is a finite set of positive real numbers.
Let $c$ be its minimum.
We claim that the first inequality holds whenever
$\delta \le c$,
and the second holds
(if $(f_1 \vee f_2)(\alpha) > f_3(\alpha)$
whenever $f_3(\alpha) < \infty$)
whenever $\varepsilon \le \delta \le c/2$.
Indeed, fix $\delta \le c$ and any
$\alpha \in \wtilde\Phi$.

If $f_3(\alpha) = \infty$,
then $(f_1 \vee f_2)(\alpha) = \infty$, 
so
\begin{equation}
\tag{$*$}
(f_{1, \delta} \vee f_{2, \delta})(\alpha)
\ge (f_1 \vee f_2)(\alpha) = \infty
= f_{3, \delta}(\alpha) + (c - \delta).
\end{equation}

Now suppose that
$f_3(\alpha) \in \sset{s, s{+}}$ with $s \in \R$.
By Lemma \ref{lem:concave-vee-achieved}, there exists a pair
$(a_m)_m$ and $(b_n)_n$ of non-empty finite
sequences in $\wtilde\Phi$, each of length no more
than $C(\Phi)$, such that $\sum_m a_m + \sum_n b_n = \alpha$
and
$\sum_m f_{1, \delta}(a_m) + \sum_n f_{2, \delta}(b_n)
= (f_{1, \delta} \vee f_{2, \delta})(\alpha)$.

If some $f_1(a_m)$ or $f_2(b_n)$ is $\infty$, then
$(f_{1, \delta} \vee f_{2, \delta})(\alpha)
= \infty \ge f_{3, \delta}(\alpha)$.

If all $f_1(a_m)$ and $f_2(b_n)$ are in \R
and $\sum_m f_1(a_m) + \sum_n f_2(a_n) = s$, then
$f_3(\alpha) = s = (f_1 \vee f_2)(\alpha)$
and
\begin{equation}
\tag{$\dag$}
(f_{1, \delta} \vee f_{2, \delta})(\alpha)
= \sum_m f_1(a_m) + \sum_n f_2(b_n)
= s = f_{3, \delta}(\alpha).
\end{equation}
If all $f_1(a_m)$ and $f_2(b_n)$ are in \R
and $\sum_m f_1(a_m) + \sum_n f_2(a_n) > s$, then
\begin{equation}
\tag{$**$}
(f_{1, \delta} \vee f_{2, \delta})(\alpha)
= \sum_m f_1(a_m) + \sum_n f_2(b_n)
\ge s + c
\ge f_{3, \delta}(\alpha) + (c - \delta).
\end{equation}

Now suppose that some $f_1(a_m)$ or $f_2(b_n)$ is not in
$\R \cup \sset\infty$; say $f_1(a_{m_0}) = d{+}$
with $d \in \R$.
Since
$$
\Bigl(d + \sum_{m \ne m_0} f_1(a_m) + \sum_n f_2(b_n)\Bigr){+}
\ge (f_1 \vee f_2)(\alpha)
\ge f_3(\alpha) \ge s,
$$
we have
$d + \sum_{m \ne m_0} f_1(a_m) + \sum_n f_2(b_n)
\ge s$, hence
\begin{multline}
\tag{$\ddag$}
(f_{1, \delta} \vee f_{2, \delta})(\alpha)
% intermediate step:
% = f_{1, \delta}(a_{m_0}) + \sum_{m \ne m_0} f_{1, \delta}(a_m)
%	+ \sum_n f_{2, \delta}(b_n) \\
\ge (d + \delta) + \sum_{m \ne m_0} f_1(a_m)
	+ \sum_n f_2(b_n)
\ge s + \delta \ge f_{3, \delta}(\alpha).
\end{multline}
If further $(f_1 \vee f_2)(\alpha) > f_3(\alpha)$, then similar
reasoning gives
\begin{equation}
\tag{$*{*}*$}
(f_{1, \delta} \vee f_{2, \delta})(\alpha)
% intermediate step:
% = f_{1, \delta}(a_{m_0}) + \sum_{m \ne m_0} f_1(a_m)
%	+ \sum_n f_2(b_n) \\
\ge f_{3, \delta}(\alpha) + \min \{c,\delta\}.
\end{equation}

By ($*$), ($\dag$), ($**$), and ($\ddag$), we have
$(f_{1, \delta} \vee f_{2, \delta})(\alpha)
\ge f_{3, \delta}(\alpha)$.
By ($*$), ($**$), and ($*{*}*$), if
$(f_1 \vee f_2)(\alpha) > f_3(\alpha)$ and
$\varepsilon \le \delta \le c/2$,
then we have
$(f_{1, \delta} \vee f_{2, \delta})(\alpha)
\ge f_{3, \delta}(\alpha) + \min\{c - \delta, \delta\}
\ge f_{3, \delta}(\alpha) + \varepsilon$.
\end{proof}

%%%% Here are the indexes of notation and terminology
%%%% referred to in comments near the top of the file.
%%%% In the final version, several lines below are commented out.
%%%% They are for including indexes that were automatically generated.
%%%% In the final version, instead we explicitly include modified
%%%% versions of the already-generated indexes.

%\Printindex{exp-notation}{Index of notation}
%\Printindex{exp-terminology}{Index of terminology}

%%%% Here are the modified versions of the generated files
\newpage %%% explicit formatting
%%% This is a hand-edited version of an automatically-generated index.
%%% Edited how?
%%% 1. Items are in order of appearance in the document.
%%% 2. We refer to section numbers, not page numbers.
%%% 3. Section heading added.
%%% 4. "theindex" environment changed to "trivlist" inside "multicols"
%%% 5. Every \item X, Y changed to \item X, {Y}.
%%% 6. \item redefined locally.
%%% 7. When we submitted a revised version, I didn't regenerate the
%%%	index and repeat the process above.  Instead, I just looked for
%%%	\indexme and \indexmem tags that had been added since the previous
%%%	submission, and inserted appropriate entries below.

\setlength{\columnseprule}{.5pt}
\begin{multicols}{3}[\section*{Index of notation}]
\let\olditem\item
\def\item#1,#2{\olditem \makebox[.3\textwidth]{#1 \hfill #2}}

\begin{trivlist}

  \item \ensuremath {\tbG }, {\S\ref{sec:notation-generalities}}%7
  \item \ensuremath {F ^{\mathrm {sep}}}, {\S\ref{sec:buildings}}% 8
  \item \ensuremath {F ^{\mathrm {tame}}}, {\S\ref{sec:buildings}}%8
  \item \ensuremath {F ^{\mathrm {un}}}, {\S\ref{sec:buildings}}%8
  \item \ensuremath {\ff _F}, {\S\ref{sec:buildings}}%8
  \item \ensuremath{\AA(\bS, F)}, {\S\ref{sec:buildings}}
	%%Manually inserted.
  \item \ensuremath {U _\varphi }, {\S\ref{sec:buildings}}%9
  \item \ensuremath {\lsub F U_\varphi }, {\S\ref{sec:buildings}}%9
  \item \ensuremath {V _F(Z({\ensuremath {\mathbf G}}\xspace  ))}, {\S\ref{sec:buildings}}%9
  \item \ensuremath {\PPsiGTF {\wtilde }{{\ensuremath {\mathbf G}}\xspace  }{{\ensuremath {\mathbf S}}\xspace  }{E}}, {\S\ref{sec:buildings}}%9
  \item \ensuremath {\PPsiGTF {}{{\ensuremath {\mathbf G}}\xspace  }{{\ensuremath {\mathbf S}}\xspace  }{E}}, {\S\ref{sec:buildings}}%9
  \item \ensuremath {G _r}, {\S\ref{sec:filtrations-and-depth}}%13
  \item \ensuremath {G _{x,r}}, {\S\ref{sec:filtrations-and-depth}}%13
  \item \ensuremath {\Lie (G)_{x, r}}, {\S\ref{sec:filtrations-and-depth}}%13
  \item \ensuremath {\depth _x}, {\S\ref{sec:filtrations-and-depth}}%13
  \item \ensuremath {\depth }, {\S\ref{sec:filtrations-and-depth}}%13
  \item \ensuremath {\ms G_x}, {\S\ref{sec:filtrations-and-depth}}%14
  \item \ensuremath{\lsub\bT G_{x, f}}, {\S\ref{sec:concave_gen}}
	%%Manually edited.
  \item \ensuremath {\vG_{x, \vec r}}, {\S\ref{sec:concave_gen}}%28
	%%Manually edited.
  \item \ensuremath{\GG^G_d}, {\S\ref{sec:normal}}
	%%Manually inserted.
  \item \ensuremath {\CC {{\ensuremath {\mathbf G}}\xspace  }{r}({\ensuremath {\underline {\gamma }}}\xspace  )}, {\S\ref{sec:normal}}%43
  \item \ensuremath {\ZZ {{\ensuremath {\mathbf G}}\xspace  }{r}({\ensuremath {\underline {\gamma }}}\xspace  )}, {\S\ref{sec:normal}}%43
  \item \ensuremath{\CC G r(\ugamma)}, {\S\ref{sec:normal}}
  \item \ensuremath{\ZZ G r(\ugamma)}, {\S\ref{sec:normal}}
	%%Both manually inserted.
  \item \ensuremath {\text {(\textbf {Gd}$^{{\ensuremath {\mathbf G}}\xspace  (E)}$)}}, {\S\ref{sec:normal}}%43
  \item \ensuremath {[ {\ensuremath {\underline {\gamma }}}\xspace  ; x, r]^{(j)}}, {\S\ref{sec:normal}}%44
  \item \ensuremath {[ {\ensuremath {\underline {\gamma }}}\xspace  ; x, r]}, {\S\ref{sec:normal}}%44
  \item \ensuremath {\CC G r({\ensuremath {\underline {\gamma }}}\xspace  )}, {\S\ref{sec:normal}}%44
  \item \ensuremath {\ZZ G r({\ensuremath {\underline {\gamma }}}\xspace  )}, {\S\ref{sec:normal}}%44
  \item \ensuremath {\udc ^{(j)}}, {\S\ref{sec:normal}}%44
  \item \ensuremath {\udc }, {\S\ref{sec:normal}}%44
  \item \ensuremath {\BB _r(\gamma )}, {\S\ref{sec:normal-conj}}%60
  \item \ensuremath {\CC G r(\gamma )}, {\S\ref{sec:normal-conj}}%60
  \item \ensuremath {\CC {{\ensuremath {\mathbf G}}\xspace  }r(\gamma )}, {\S\ref{sec:normal-conj}}%60
  \item \ensuremath {\ZZ G r(\gamma )}, {\S\ref{sec:normal-conj}}%60
  \item \ensuremath {\ZZ {{\ensuremath {\mathbf G}}\xspace  } r(\gamma )}, {\S\ref{sec:normal-conj}}%60
  \item \ensuremath {[ \gamma ; x, r]^{(j)}}, {\S\ref{sec:normal-conj}}%60
  \item \ensuremath {[ \gamma ; x, r]}, {\S\ref{sec:normal-conj}}%60
  \item \ensuremath {\dc ^{(j)}}, {\S\ref{sec:normal-conj}}
  \item \ensuremath {\dc }, {\S\ref{sec:normal-conj}}
  \item \ensuremath {f _1 \vee f_2}, {\S\ref{sec:concave-fns}}%72
  \item \ensuremath{h_\delta}, {\S\ref{sec:concave-fns}}
	%%Manually inserted.

\end{trivlist}
\end{multicols}

%%% This is a hand-edited version of an automatically-generated index.
%%% Edited how?
%%% 1. "theindex" environment changed to "trivlist" inside a section
%%% 2. page numbers changed to definition numbers.

\section*{Index of terminology}

\begin{trivlist}

  \item absolutely semisimple, Definition \ref{defn:top-F-ss}
  \item admissible sequence (of numbers), Definition \ref{defn:admissible}
  \item affine root (unusual usage), \S\ref{sec:buildings}
  \item $r$-approximation, Definition \ref{defn:r-approx}

  \indexspace

  \item character values of a semisimple element, Definition \ref{defn:root-value}
  \item compatibly filtered $F$-subgroup, Definition \ref{defn:compatibly-filtered}

  \indexspace

  \item filtration-preserving map, Definition \ref{defn:compatibly-filtered}

  \indexspace

  \item good element, Definition \ref{defn:good}
%    \subitem of depth $d > 0$, Definition \ref{defn:good}
  \item good root values, Definition \ref{defn:good-root-values}
  \item good sequence, Definition \ref{defn:bracket}

  \indexspace

  \item Levi subgroup
     \subitem $F$-Levi subgroup, p.~\pageref{defn:levi}
     \subitem $E$-Levi $F$-subgroup, p.~\pageref{defn:levi}
	%%Above item (and subitems) manually added.  Page
	%%number to be filled in.

  \indexspace

  %\item normal $r$-approximation, Definition \ref{defn:r-approx}
	%%Manually deleted (subitem of $r$-approximation).

  %\indexspace

  \item $(x,f)$-positive, Definition \ref{defn:f-positive}
%    \subitem pair, Definition \ref{defn:f-positive}
%    \subitem point, Definition \ref{defn:f-positive}
  \item preserves filtrations, Definition \ref{defn:compatibly-filtered}

  \indexspace

  \item root values of a semisimple element, Definition \ref{defn:root-value}

  \indexspace

  \item sufficiently large field, Definition \ref{defn:sufficiently-large}

  \indexspace

  \item $F$-tame, Definition \ref{defn:tame}
%    \subitem element, Definition \ref{defn:tame}
%    \subitem group, Definition \ref{defn:tame}
%    \subitem modulo a subgroup
%      \subsubitem element, Definition \ref{defn:tame}
%      \subsubitem group, Definition \ref{defn:tame}
  \item tame Levi sequence, Definition \ref{defn:tame-Levi-sequence}
  \item tame reductive sequence, Definition \ref{defn:tame-reductive-sqnc}

\end{trivlist}

%%%% Bibliography starts here

\end{document}